\documentclass[11pt,a4paper,reqno]{amsart}
\usepackage{graphicx}
\usepackage{epsfig}
\usepackage{amsmath}
\usepackage{amsfonts}
\usepackage{amssymb}
\usepackage{amscd}
\newtheorem{theorem}{Theorem}[section]

\newtheorem{lemma}[theorem]{Lemma}

\begin{document}


\bigskip
\title [ Approximation in generalized urn models  ]
{ Asymptotic approximation for a certain class of statistics associated with generalized urn models  }

\author{Sherzod M Mirakhmedov}

 \maketitle
\begin{center}
{Institute of Mathematics, Academy of Sciences. Tashkent. Uzbekistan.
 e-mail: shmirakhmedov@yahoo.com}
\end{center}

\bigskip
\bigskip
Abstract. The generalized urn model is defined in terms of a collection of random frequencies, which is treated as a collection of conditionally independent random variables given their sum. We consider a class of random variables that can be represented as series of functions of these frequencies. The general results concern asymptotic normality, estimates of remainder terms, and Edgeworth-type asymptotic expansions for such random variables. We then consider  random variables defined in various versions of the generalized urn model. These models include the random allocation of particles of several types into finite or infinite numbers of cells, random sampling without replacement from a finite population, random coverage of a circle, random forests, and random permutations.





\bigskip
\bigskip
{\footnotesize \noindent \textbf{Keywords: Urn models, occupancy counts, sampling without replacement,random coverage of a circle, random forests, random permutations.} 

\bigskip
\noindent \textbf{AMS Subject Classification:62G20,62G10 60F05 }

\bigskip

\section{	INTRODUCTION.} 

  Urn models (also known as random allocation schemes) are useful tools that allow us to formulate and better understand many combinatorial problems in probability and statistics. Urn models naturally arise in areas such as statistical mechanics, clinical trials, cryptography, computer science, and biology, among others. Consequently, there exists a rich  literature dedicated to studying various random variables (r.v.s) of interest in specific urn models, see e.g.  Ivanov et al (1985), Johnson, N.L., Kotz, S.(1977), Kotz and Balakrishnan (1997), Mikhaylov (1993 ), Mirakhmedov et al (2014) and references therein.\\
\indent
  Loosely speaking, an urn model can be described as follows: an urn (or alternatively, a population) contains a finite or infinite number of different types of items labeled by 1, 2, ...,. According to a predetermined rule (scheme), a sample of size \textit{n} is drawn from the urn. This results in a random vector (r.vec.) $\bar{\eta}(n) = (\eta_{1}(n), \eta_{2}(n), \dots)$, where $\eta_{1}(n) + \eta_{2}(n) + \cdots = n$, and $\eta_{m}(n)$ represents the frequency of the \textit{m}th type of items in the sample. In practical settings, we often deal with species or categories (types) of sampling units (items). The main focus is placed on the sampled array of these types, denoted as $\bar{\eta}(n)$, since it provides insights into the underlying population frequencies. Such species-sampling challenges appear, for instance, in database query optimization. In this context, the sampling units could represent entries within database columns, while the species correspond to all distinct values present in those columns. Similarly, in literature analysis, the sampling units might be words found in an author’s known works, with the species representing all words known to the author. Another application emerges in limiting disclosure risks, where the sampling units may consist of individuals or firms listed in anonymized microdata files, and the types encompass unique combinations of variable values that could potentially lead to indirect identification of these entities. Models involving infinitely many items can serve as approximations for sampling from large finite populations, or from the set of possible values of a discrete random variable.\\
  \indent
   The r.vec. $\bar{\eta}(n)$ is a sufficient statistics, for the model, and common r.v.s of interest in practices are  functionals of this r.vec.. On the other hand a characteristic property of many urn models is the fact that the  distribution of the vector of frequencies  $\bar{\eta}(n)$ can be represented as the joint conditional distribution of independent r.v.s $ \xi_{1}, \xi_{2},... $ given their total $ \xi_{1}+ \xi_{2}+... =n$ , see, e.g.,Holst (1972,1979),  Rao and Schweitze (1985), Kolchin (1985), Holst and Rao (1980), Mirakhmedov (1992), Mirakhmedov et al (2014).  Thus the following definition includes several commonly known urn models as special cases. \\ 
\indent Assume a given integer $\textit{n}>1$ and $\bar{\xi}=( \xi_{1}, \xi_{2},...) $ is a r.vec. with independent and non-negative integer-valued components such that $ P\lbrace \zeta=n\rbrace >0 $, where $ \zeta=\xi_{1}+ \xi_{2}+... , .$ Also let $\bar{\eta}(n) = (\eta_{1}(n), \eta_{2}(n), \dots)$ be a r.vec. whose distribution is defined by 
 \begin{equation}\label{eq1} 
\mathcal{L}(\bar{\eta}(n))=\mathcal{L}(\bar{\xi}\mid\zeta=n),
\end{equation}
where $ \mathcal{L}(X)$ stands for the distribution of the r.vec. \textit{X}. Clearly, the r.v.s $\eta_{m}(n)$ are non-negative integer-valued, and (1) implies that $\eta_{1}(n) + \eta_{2}(n) + \dots = n$ a.s. . Thus, the r.v. $\eta_{m}(n)$ can be interpreted as the number of \textit{m}th type of items appearing in the sample of size \textit{n} drawn from an urn containing either a  finite or infinite number types of items, where the distributions of the r.v.s $  \xi_{m}$ define the sampling scheme via (1). The model defined in (1) is what we call a {\it generalized urn model} (GUM). If we assume $ \xi_{m}=0 $ for all $m>N$, where $N=N(n)\rightarrow\infty$ as $n\rightarrow\infty ,$ we arrive at the GUM introduced by Mirakhmedov et al (2014), where the asymptotic theory and higher-order expansions for a certain class of statistic are presented.\\

Although a single GUM of the form (1) covers many specific urn models, it does not cover several probabilistic models that have been extensively studied in the literature, such as balls-in-infinite-boxes occupancy schemes with the random allocation of balls of several types, the random allocation of particles in sets, and statistics based on several samples from one or more populations. To include such models, we consider $s\geq 1$ independent GUMs. Specifically, let $\bar{\xi}_{k}=(\xi_{k,1}, \xi_{k,2},...) $, $k=1,...,s $, be a collection of independent r.vec.s whose components are independent non-negative integer-valued r.v.s such that $P\lbrace\zeta_{k}=n_{k}\rbrace>0 $  for a given integer $ n_{k}\geq 2 ,$ where $\zeta_{k}= \xi_{k,1}+\xi_{k,2}+.... $, $k=1,...,s .$ The conditional distribution $ \mathcal{L}(\bar{\xi}_{k}\mid\zeta_{k}=n_{k}) $  generates a r.vec. of frequencies  $ \bar{\eta}_{k}(n_{k})=(\eta_{k,1}(n_{k}), \eta_{k,2}(n_{k}),...)  $ via equation (1) forming a collection of s independent GUMs.\\
\indent In this paper we are interested in the following general class of statistics:
\begin{equation} \label{eq2}
       \mathbb{R}^{\eta}(n)= \sum^{\infty}_{m=1}{f_{m}(\eta_{1,m}(n_{1}),...,\eta_{s,m}}(n_{s})) , 
\end{equation}
where $ f_{m}(u_{1},...,u_{s})$, \ $m=1,2,...$, is a sequence of functions (possibly random) defined for  $u_{1}\geq0,...,u_{s}\geq0 $, such that the series (2) is a.s. converges for every $n=(n_{1},...,n_{s})$. These functions may depend on \textit{n}, but to simplify notation, we omit this dependency. \\
  We emphasize that the present setting involves triangular arrays, in which the parameters of the GUMs and the distributions of the r.v.s \( \xi_{k,m} \) may vary as \( \min(n_1,\ldots,n_s) \to \infty \), formally along non-decreasing sequences \( n_{1,\nu}, \ldots, n_{s,\nu} \) as \( \nu \to \infty \).
  
  \indent  This paper serves as a kind of survey extending the results of Mirakhmedov et al. (2014) and updating results on the infinite occupancy scheme by Mirakhmedov (1989, 1992). The results presented in Sections 2 and 3 are of a general nature, proofs of which are given in Section 4. Section 5 focuses on specific versions of the GUM and on the random variable $\mathbb{R}^{\eta}(n)$. The models considered include the random allocation of particles of $s\ge 1$ types into finite or infinite numbers of cells, sampling schemes without replacement from a finite population, random coverage of a circle, random forests, and random permutations.  Most of the results are novel, while others are dispersed throughout the literature.\\
\indent  Throughout the paper, we adhere to the following conventions.   We assume that the considered moments of the r.v.s exist, the series of r.v.s converges almost surely, and the series of numbers converge absolutely. The quantities we consider may depend on $n_{1}, \ldots, n_{s}$, but for simplicity, we will omit this dependence in future notations as well as in those already introduced, unless such omission causes ambiguity. The symbols $c$ and $C$ represent universal positive constants that may vary depending on the context, $\theta $ stands for the quantity such that $\mid\theta\mid\le 1$. All asymptotic relationships are examined under the condition $\min(n_1,\ldots,n_s)\to\infty .$ For  given quantities $a_j$ (which may be fixed or random) we denote $a^s=(a_1,\dots, a_s)$ for short, for instance, $\xi_m^s=(\xi_{1,m},\dots,\xi_{s,m})\;$, $\tau^s=(\tau_1,\dots,\tau_s)$, etc.;  $ \parallel a \parallel$ and $\langle a,b\rangle$ stand for the Euclid norm and scalar product, respectively, of vectors  $a , b \in R^{s}.$   Equality and subtraction of two vectors is defined as component-wise. Finally,    $\Phi(x)$ is a standard normal distribution, $Poi(\lambda),\,\, Bi(k,p)\,\, \text{and} \,\, G(p)$ denotes, respectively, the Poisson, Binomial and Geometrical distributions  with indicated parameters, $ \xi\sim F $ stand for " a r.v. $\xi$ has distribution F " . \\
  
\begin {center}.\\

\section{A FORMAL CONSTRUCTION OF THE ASYMPTOTIC EXPANSION}
\end {center}

We still use notation of Section 1, and  set $n:=n^{s},$ $ x_{k}:=\left( n_{k}-\lambda_{k}\right) /B_{k},$ and hence $x^{s}=(x_{1},...,x_{s}),$ where
 $$   \lambda_{k}=\sum^{\infty}_{m=1}E\xi_{k,m} ,\quad  B^{2}_{k}= \sum^{\infty}_{m=1}var\xi_{k,m} .$$ \\
We adopt also the following notation:
$$\Lambda(n)=\sum^{\infty}_{m=1}Ef_{m}(\xi^{s}_{m}) ,\quad  \gamma_{k}=\dfrac{1}{B^{2}_{k}}\sum^{\infty}_{m=1}cov\left( f_{m}(\xi^{s}_{m}),\xi_{k,m}\right),$$ 

$$g_{m}(y^s)=f_{m}(y^s)-Ef_{m}(\xi^{s}_{m})-\sum^{s}_{k=1}\gamma_{k}(y_{k}-E\xi_{k,m}),$$

$$  G^{\eta}(n)=\sum^{\infty}_{m=1}g_{m}(\eta^{s}_{m}), \ \ \   G^{\xi}(n)=\sum^{\infty}_{m=1}g_{m}(\xi^{s}_{m}),$$ 
\begin{equation} \label{eq3}
\sigma^{2}(n)=varG^{\xi}(n)=\sum^{\infty}_{m=1}varf_{m}(\xi^{s}_{m})- \sum^{s}_{k=1}B^{2}_{k}\gamma^{2}_{k}. 
\end{equation}
Since $\eta_{k,1}+\eta_{k,2}+\dots =n_{k} $ we have 
   $$ G^{\eta}(n)=\mathbb{R}^{\eta}(n)-\Lambda(n)-\sum^{s}_{k=1} x_{k}B_{k}\gamma_{k}.$$
Also, one can see that
\begin{equation} \label{eq5}
  Eg_{m}(\xi^{s}_{m})=0 ,\  m=1,2,... ,\quad \mathrm{and} \ \   \sum^{\infty}_{m=1}cov\left( g_{m}(\xi^{s}_{m}),\xi_{k,m}\right)=0 , \ \  k=1,...,s.
\end{equation}\\
Under some mild conditions, one can show that as $ \min (n_{1},...,n_{s})\rightarrow\infty $
$$  E\mathbb{R}^{\eta}(n)= \Lambda(n)+\sum^{s}_{k=1}\left[ x_kB_k\gamma_k  -\frac{1-x^{2}_k}{2B^{2}_k}\sum^{\infty} _{m=1}E\left(  g_{m}(\xi^{s}_{m})(\xi_{k,m}-E\xi_{k,m})^{2}\right) (1+o(1))\right], $$
 and $$ var \mathbb{R}^{\eta}(n)=\sigma^{2}(n)(1+o(1)).$$    

  Due to (1) and the fact that $ \bar{\eta}_{1},...,\bar{\eta}_{s} $ are a collection of independent r.vec.s, we have
 \begin{equation}\label{eq4}
 E\exp\lbrace itG^{\eta}(n) \rbrace = E\left[\exp\lbrace itG^{\xi}(n) \rbrace\mid\zeta^{s}=n^s\right] . 
\end{equation}
On the other hand by the total expectation formula 
$$ E\exp\left\lbrace  itG^{\xi}(n)+i\langle \tau^{s},(\zeta^{s}-n^s)\rangle\right\rbrace  $$ 
$$ =\sum^{\infty}_{k_{1}=0}...\sum^{\infty}_{k_{s}=0}E\left[  \exp\left\lbrace  itG^{\xi}(n)\right\rbrace \mid\zeta^{s}=(k_{1},...,k_{s})\right]\prod^{s}_{l=1}e^{i\tau_{l}(k_{l}-n_{l})}P\lbrace\zeta_{l}=k_{l}\rbrace .$$
This together with (4) implies, by Fourier inversion 

 $$E\exp\left\lbrace  itG^{\eta}(n)\right\rbrace =\dfrac{1}{(2\pi)^{s}\prod^{s}_{l=1}P\lbrace\zeta_{l}=n_{l}\rbrace}$$
 \begin{equation} \label{eq5}
\cdot\int _{-\pi}^{\pi}...\int _{-\pi}^{\pi}E\exp\left\lbrace  itG^{\xi}(n)+i\langle \tau^{s},(\zeta^{s}(n)-n)\rangle\right\rbrace  d\tau_1\dots d\tau_s .
\end{equation}
 
Next, consider the standardized r.v.s \,
\begin{equation}\label{eq6}
\tilde{\xi}_{l,m}=B^{-1}_{l}(\xi_{l,m}-E\xi_{l,m}), \quad \tilde{\xi}^{s}_{m}=(\tilde{\xi}_{1,m},..., \tilde{\xi}_{s,m})  ,\quad \tilde{g}_{m}(\xi^{s}_{m})=\sigma^{-1}(n)g_{m}(\xi^{s}_{m})
\end{equation}
and set
$$\psi_{m}(t,\tau^{s})= E\exp\left\lbrace  it\tilde{g}_{m}(\xi^{s}_{m})+i\langle\tau^{s},\tilde{\xi}^{s}_{m}\rangle\right\rbrace ,\quad      \Psi(n,t,\tau^{s})=\prod^{\infty}_{m=1}\psi_{m}(t,\tau^{s}),     $$
\begin{equation} \label{eq8}
\Theta(n,x^{s},t)=\dfrac{1}{(2\pi)^{s/2}} \int _{-\pi B_{1}}^{\pi B_{1}}...\int _{-\pi B_{s}}^{\pi B_{s}} e^{-i\langle\tau^{s},x^{s}\rangle}\Psi(n,t,\tau^{s})d\tau_{1}...d\tau_{s} .
\end{equation}
Then equation (5) and the fact that $\prod^{s}_{l=1}B_{l}P\lbrace\zeta_{l}=n_{l}\rbrace=\Theta(n,x^{s},0)$ gives the following \\

{\bf Lemma 1.1. } {\it It holds 
\begin{equation} \label{eq9}
\varphi(n,t):=E\exp\left\lbrace  it\sigma^{-1}(n)G^{\eta}(n)\right\rbrace =\dfrac{\Theta(n,x^{s},t)}{\Theta(n,x^{s},0)} .
\end{equation}}

\indent 
Lemma 1 provides the basis for the formal construction of an asymptotic expansion of the characteristic function  (ch.f.) $\varphi(n,t)$. Indeed: The $ \Psi(n,t,\tau^{s})$   is the ch.f. of the sum (the series) of  independent (\textit{s}+1) -dimensional r.vec.s $(\tilde{g}_{m}(\xi^{s}_{m}),\tilde{\xi}^{s}_{m}) $.  Due to (4) this sum has zero mean, a unit covariance matrix and uncorrelated components. Keeping these facts in mind, we can now derive the asymptotic expansion of  $ \Psi(n,t,\tau^{s}).$  Under appropriate conditions, the terms of this expansion take the form of polynomials, say $P_{k}(n,t,\tau^{s})$ , $k=1,2,...$,  with respect to $t,\tau_{1},...,\tau_{s}$ and contain the  common factor $\exp\lbrace -(t^{2}+\tau^{2}_{1}+...+\tau^{2}_{s})/2\rbrace.$  The degree of $P_{k}(n,t,\tau^{s})$ is $3k$, and the minimal degree is $k+2$. 
Multiplying these terms by $\exp\lbrace -i\langle\tau^{s},x^{s}\rangle\rbrace $  allows us to integrate them with respect to each $\tau_{1},...,\tau_{s}$   over the interval $(-\infty,\infty)$. The result of these integrations, say $V(n,x^{s},t)$, can be regarded as a formally asymptotic expansion of  $\Theta(n,x^{s},t).$  Ultimately, we obtain the asymptotic expansion of $ \varphi(n,t)$   by dividing $V(n,x^{s},t)$ by $V(n,x^{s},0).$ This procedure involves complex but manageable algebra and consists of the following four steps. Let's assume $E|g_{m}(\xi^{s}_{m})|^{d}<\infty  $ and $  E|\xi_{l,m}|^{d}<\infty $, $m=1,2,...$; $l=1,...,s$; some $d\geq3 .$ \\
 
 Step 1. Determine the polynomials $P_{k}(n,t,\tau^{s})$ using, for example,well developed method in theory of sum of independent r.vec.s , see, book by Bhattacharya and Rao (1976, Chapter 2), also Lemma 2 of Section 4 below. 
  Then,  for instance, we obtain: $ \, P_{0}(n,t,\tau^{s})=1 ,$ 
  $$ P_{1}(n,t,\tau^{s})=\dfrac{i^{3}}{6}\sum_{m=1}^{\infty} E\left( t\tilde{g}_{m}(\xi^{s}_{m})+\langle\tau^{s},\tilde{\xi}^{s}_{m}\rangle\right) ^{3},$$
  
  $$ P_{2}(n,t,\tau^{s})=\dfrac{i^{4}}{24}\sum_{m=1}^{\infty} \left[ E\left( t\tilde{g}_{m}(\xi^{s}_{m})+\langle\tau^{s},\tilde{\xi}^{s}_{m}\rangle\right) ^{4}-3 \left( E\left( t\tilde{g}_{m}(\xi^{s}_{m})+\langle\tau^{s},\tilde{\xi}^{s}_{m}\rangle\right) ^{2}\right) ^{2}\right] +\dfrac{1}{2}P_{1}^{2}(n,t,\tau^{s}).$$
 
 Step 2. Define the polynomials $G_{k}(n,x^{s},t)$  wrt \textit{t} for $k=1,2,...$  as
\begin{equation} \label{eq10}
 G_{k}(n,x^{s},t)=\dfrac{1}{(2\pi)^{s/2}}\int_{-\infty}^{\infty}...\int_{-\infty}^{\infty}P_{k}(n,t,\tau_{1},...,\tau_{s})\exp\left\lbrace  -\frac{1}{2}\sum_{l=1}^{s}(\tau_{l}+ix_{l}))^{2}\right\rbrace  d\tau_{1}...d\tau_{s} 
\end{equation}

Then
$$ V(n,x^{s},t)=\exp{\left( -\dfrac{t^{2}+\parallel x^{s} \parallel ^{2}}{2}\right) }\left( 1+G_{1}(n,x^{s},t)+G_{2}(n,x^{s},t)+...+G_{d-3}(n,x^{s},t)\right) + \dots  $$\\
should be considered as asymptotic expansion of $\Theta(n,x^{s},t).$\\

We note that the coefficients of $P_{k}(n,t,\tau^{s})$ , and hence of $G_{k}(n,x^{s},t)$  , only involve the quantities  $\prod_{l=1}^{k}\left( \alpha_{j,j_{1},...,j_{l}}\right) ^{\nu_{l}} $ , where $j+j_{1}+...+j_{k}=l+2$ for the non-negative integers $j, j_{1},...,j_{r}$ and $\nu_{1}+2\nu_{2}+...+k\nu_{k}=k $, 
$\alpha_{j,j_{1},...,j_{l}}=\Sigma_{m=1}^{\infty}E\left( \tilde{g}^{j}_{m}(\xi^{s}_{m})\tilde{\xi} ^{j_{1}}_{1,m}\cdot ... \cdot \tilde{\xi} ^{j_{s}}_{s,m}\right). $ For the details see \cite{14}, Chapter 2, Section 7.
 For instance, having designated for short  
$$ a^{(k,l)}=\sum_{m=1}^{\infty}(Eg^{k}_{m}(\xi^{s}_{m}))^{l} ,\quad     b_{j}^{(k,\nu)}=\sum_{m=1}^{\infty}Eg^{k}_{m}(\xi^{s}_{m})(\xi_{j,m}-E\xi_{j,m})^{\nu} ,$$
$$ c_{l,k}=\sum_{m=1}^{\infty}Eg_{m}(\xi^{s}_{m})(\xi_{l,m}-E\xi_{l,m})(\xi _{k,m}-E\xi _{k,m} ),  $$
 taking into account (4) and the fact that $\tilde{\xi}_{l,m}$  are independent r.v.s whose mean is equal to zero we obtain
 $$ P_{1}(n,t,\tau^{s})=\dfrac{(it)^{3}}{6} \dfrac{a^{(3,1)}}{\sigma^{3}} +\dfrac{(it)^{2}}{2\sigma^{2}}\sum_{l=1}^{s}\dfrac{(i\tau_{l})b_{l}^{(2,1)}}{B_{l}}+\dfrac{it}{2\sigma}\Big(\sum_{l=1}^{s}\dfrac{(i\tau_{l})^{2}b_{l}^{(1,2)}}{B^{2}_{l}}$$
 $$+\sum_{l,k=1,k\neq l}^{s}\dfrac{(i\tau_{l})(i\tau_{k})c_{l,k}}{B_{l}B_{k}}\Big)
  +\sum_{l=1}^{s}\dfrac{(i\tau_{l})^{3}b^{(0,3)}_{l}}{6B^{3}_{l}},  $$
 and
 
$$ G_{1}(n,x^{s},t)=\dfrac{(it)^{3}}{6} \dfrac{a^{(3,1)}}{\sigma^{3}}+\dfrac{(it)^{2}}{2\sigma^{2}}\sum_{l=1}^{s}\dfrac{x_{l}b_{l}^{(2,1)}}{B_{l}}$$
$$-\dfrac{it}{2\sigma)}\left( \sum_{l=1}^{s}\dfrac{(1-x^{2}_{l})b_{l}^{(1,2)}}{B^{2}_{l}}+ \sum_{l,k=1,k\neq l}^{s}\dfrac{x_{l}x_{k}c_{l,k}}{B_{l}B_{k}}\right) -\sum_{l=1}^{s}\dfrac{((x^{3}_{l}-3x_{l})b^{(0,3)}_{l}}{6B^{3}_{l}}.$$
 A similar expression could be written for $P_2(n, t, \tau^s)$ and $G_{2}(n,x^{s},t)$, but it would result in a lengthy formula, so we omit it here.\\

 Step 3. Find the terms of asymptotic expansion of  $[V(n,x^{s},0)]^{-1}.$ To do this, define $Q_{k}(n,x^{s})$ from the equation
 
$$ \sum_{k=0}^{\infty}(-1)^{k}\left( \sum_{\nu =0} ^{d-3}
G_{\nu}(n,x^{s},0)\right) ^{k}=\sum_{m=0}^{\infty}Q_{m}(n,x^{s}).$$
Then

$$ Q_{m}(n,x^{s})=m!\sum\prod_{l=1}^{d-3}\dfrac{1}{m_{l}!}G_{l}^{m_{l}}(n,x^{s},0)$$

where the summation is over all non-negative integer $m_{1},...,m_{d-3}$ such that $m_{1}+2m_{2}+...+(d-3)m_{d-3}=m$. \\

Step 4. Define
\begin{equation} \label{eq11}
W^{(d)}(n,x^{s},t)=\sum_{m=0}^{d-3}\sum_{\nu =0}^{m} G_{\nu}(n,x^{s}),t)Q_{m-\nu}(n,x^{s}).
\end{equation}
Finally we can conclude that
\begin{equation} \label{eq12}
\varphi(n,t)= \dfrac{V(n,x^{s},t)}{V(n,x^{s},0)}=\exp{\left\lbrace  -\frac{t^{2}}{2}\right\rbrace W^{(d)}(n,x^{s},t) } + r(n,t,x^s), 
\end{equation} 
where $r(n,t,x^s)$ is the remainder to be estimated, see Lemma 3 of Sec 4 below.\\
\indent
 Note that $G_{0}(n,x^{s},t)=Q_{0}(n,x^{s})=1$, so that $ W^{(3)}(n,x^{s},t)=1, $ and 
$$ W^{(4)}(n,x^{s},t) =\dfrac{(it)^{3}}{6} \dfrac{a^{(3,1)}}{\sigma^{3}} +\dfrac{(it)^{2}}{2\sigma^{2}}\sum_{l=1}^{s}\dfrac{x_{l}b_{l}^{(2,1)}}{B_{l}}-\dfrac{it}{2\sigma}\sum_{l=1}^{s}\dfrac{(1-x^{2}_{l})b_{l}^{(1,2)}}{B^{2}_{l}} .$$
Also it can be shown that
$$ W^{(5)}(n,x^{s},t)= W^{(4)}(n,x^{s},t)$$ $$+\left[ G_{2}(n,x^{s},t)-G_{2}(n,x^{s},0)-G_{1}(n,x^{s},0)(G_{1}(n,x^{s}(n),t)-G_{1}(n,x^{s}(n),0))\right]. $$

Now, the formal construction of the terms of the asymptotic expansion of the distribution function proceeds in a fairly standard manner, as described, for example, in the book by Petrov (1975), also see Mirakhmedov et al. (2014). Indeed.
 Define $\mathbb{W}^{(d)}(n,x^{s}(n),u)$ so that 
\begin{equation}\label{eq13}
\int_{-\infty}^{\infty}e^{itu}d\mathbb{W}^{(d)}(n,x^{s},u)= W^{(d)}(n,x^{s},t)e^{-{t^{2}/2}} 
\end{equation}
The function $\mathbb{W}^{(d)}(n,x^{s},u)$  can be obtained by formally substituting   

$$(-1)^{v}\dfrac{d^{v}}{du^{v}}\Phi(u)= -e^{-{u^{2}/2}} H_{v-1}(u)/\sqrt{2\pi}, \quad \mathrm{where }\quad  \Phi(u)=\dfrac{1}{\sqrt{2\pi}}\int_{-\infty}^{u}e^{-{t^{2}/2}}dt , $$
instead of $(it)^{v}$  for each $v$ in the expression $W^{(d)}(n,x^{s},t) $, where $H_{v}(u)$ is the v-th order Hermite -Chebuishev polynomial.\  
Note that $\mathbb{W}^{(3)}(n,x^{s},u)=\Phi(u)$ and by (13) 

$$\mathbb{W}^{(4)}(n,x^{s},u)=\Phi(u)$$ $$+\dfrac{1}{\sqrt{2\pi}}e^{-{u^{2}/2}}\left[ \dfrac{(1-u^{2})a^{(3,1)}}{6\sigma^{3}}-\dfrac{u}{2\sigma^{2}}\sum_{l=1}^{s}\dfrac{x_{l}b_{l}^{(2,1)}}{B_{l}}+\dfrac{1}{3\sigma}\sum_{l=1}^{s}\dfrac{(1-x^{2}_{l})b_{l}^{(1,2)}}{B^{2}_{l}}\right] . $$\\   

We emphasis that for various versions of the GUM the  r.v.s $\xi_{k,m}$  are such that $\lambda_{k}=n_{k}$ (see examples of Section 5 below), and hence $x_{k}=0,$ \, $ k=1,...,s $ \, so the terms of asymptotic expansion, i.e. the function  $\mathbb{W}^{(d)}(n,0,u)$ is considerably simplified. In particular, 
\begin{equation}\label{eq14}
\mathbb{W}^{(4)}(n,0,u)=\Phi(u)+\dfrac{1}{\sqrt{2\pi}}e^{-{u^{2}/2}}\Big[ \dfrac{(1-u^{2})a^{(3,1)}}{6\sigma^{3}}+\dfrac{1}{3\sigma}\sum_{l=1}^{s}\dfrac{b_{l}^{(1,2)}}{B^{2}_{l}}\Big] , 
\end{equation}

and
$$\mathbb{W}^{(5)}(n,0,u)=\mathbb{W}^{(4)}(n,0,u)- \dfrac{e^{-{u^{2}/2}}}{\sqrt{2\pi}}\Big[ \dfrac{(a^{(3,1)})^{2}}{72\sigma^{6}}(u^{5}-10u^{3}+15u) $$
$$+\frac{u^{3}-3u}{24\sigma^{4}}\left( a^{(4,1)} - 3 a^{(2,2)}- 3\sum_{l=1}^{s}\frac{(b^{(2,1)}_{l})^{2}}{B_{l}^{2}} -2 a^{(3,1)} \sum_{l=1}^{s}\dfrac{b_{l}^{(1,2)}}{B^{2}_{l}} \right) $$

$$ -\frac{u}{4\sigma^{2}}\left( \sum_{l=1}^{s}\dfrac{b^{(2,2)}_{l}}{B^{2}_{l}}- 2\sum_{l=1}^{s}\frac{b^{(1,2)}_{l}}{B_{l}^{2}}  -  \sum_{l=1}^{s}\dfrac{b_{l}^{(2,1)}b^{(0,3)}_{l}}{B_{l}^{4}} \right)   $$
 
 \begin{equation}\label{eq15}
  + \dfrac{u}{8\sigma^{2}}\left( 3\sum_{l=1}^{s}\dfrac{(b^{(1,2)}_{l})^{2}}{B^{4}_{l}}+\sum_{l,j=1, j\neq l}^{s}\frac{b^{(1,2)}_{l}b^{(1,2)}_{j}}{B_{l}^{2}B_{j}^{2}} + \sum_{l,j=1, j\neq l}^{s}\dfrac{(c^{(l,j)} )^{2}}{B^{2}_{l}B^{2}_{j}}\right)  \Big] . 
\end{equation}

\bigskip
\section{	THE GENERAL RESULTS.} 

As noted above, for common versions of GUM, the r.v.s $\xi_{k,m}$ depend on an arbitrary parameter that can be chosen such that $\lambda_{k}=n_{k}$ and, therefore, $x_{k}=0,$ \, $ k=1,...,s .$  All the versions of GUM that we consider in Section 5 are examples of this kind. Therefore, although the above algorithm for formally constructing the asymptotic expansion does not impose any restrictions on $x^{s}$, throughout the following discussion, we will assume that $x_{k}=0,$ \, $ k=1,...,s .$\\
 \indent
We will use the notation of the previous sections, omitting the dependence on $x^{s}$. Additionally we set 
 $$ \beta_{k}(n)=\sigma^{-k}(n) \sum_{m=1}^{\infty}E\mid g_{m}(\xi_{m}^{s})\mid^{k},\quad \varrho_{k,l}(n_{l})=B_{l}^{-k}\sum_{m=1}^{\infty}E\vert \xi_{l,m}-E\xi_{l,m}\vert ^{k},$$
 
$$ M_l(T)=\inf_{T\leq u \leq\pi}\sum_{m=1}^{\infty}\left( 1-|E\exp{\left( iu\xi_{l,m}\right)} |^{2}\right), \;\; \textrm{if}\;\; 0< T\leq\pi  ,\quad \textrm{else} \quad    M_l(T)=\infty ,$$ 
 $$ T_l=\Big( B_{l}\varrho_{3,l}(n_{l})2s(s+1)^{2}\Big) ^{-1}  ,\quad \mathbf{M}(n)=   \min_{1\le l\le s} M_{l}(T_l), $$
 
  $$ \mathbf{T}(n)=\min \left( (2(s+1)^{2}\beta_{3}(n))^{-1},\sqrt{\mathbf{M}(n)/6s} \right), $$

$$  \rho_{k}(n)=\beta_{k}(n)+\sum_{l=1}^{s}\Bigg[\varrho_{k,l}(n_{l})+B_l\exp\left\lbrace -\frac{1}{4s}\mathbf M(n)\right\rbrace \Bigg], $$ 

$$ \chi(n,a,b)=I\lbrace a<b \rbrace\int_{a}^{b}\dfrac{|\varphi(n,t)|}{|t|}dt \  ,\quad \bar{\chi}(n,a,b)=I\lbrace a<b \rbrace\int_{a}^{b}|\varphi(n,t)|dt ,$$

$$\Delta^{(d)}(n)=\sup_{-\infty < u < \infty}\mid\text{P}\lbrace  \mathbb{R}^{\eta}(n)<u\sigma(n)+\Lambda(n)\rbrace-\mathbb{W}^{(d)}(n,0,u)\mid. $$
We recall that $\mathbb{W}^{(3)}(n,0,u)=\Phi (u).$ 

\begin{theorem}\label{th 3.1}If\\
(i) $\beta_{3}(n)\rightarrow\infty,$\\
(ii) $\varrho_{3}(n)\rightarrow\infty,$\\
(iii) $ \mathbf{M}(n)\to\infty,$ and $\; \ln B_k=o\Big(M_k(T_k) \Big),\; k=1,\dots, s .$\\

then the r.v. $\mathbb{R}^{\eta}(n)$ has asymptotically normal distribution with mean $ \Lambda(n) $ and variance $ \sigma^{2}(n). $
\end{theorem}

\begin{theorem}\label{th 3.2} There exists a constant $C(s)$  such that             
$$ \Delta^{(3)}(n)\leq C(s)\left( \rho_{3}(n)+\frac{1}{\sqrt{ \mathbf{M}(n) }}\right) .$$
\end{theorem}

\begin{theorem}\label{th 3.3}  There exists a  constant  $C(s)$  such that
$$ \Delta^{(d)}(n)\leq C(s)\rho_{d}(n)+\chi(n,\mathbb{T}(n),\beta_{d}^{-1}(n)),\; d=4,5.$$
\end{theorem}
   
\begin{theorem}\label{th4} Let  $\text{R}^{\eta}(n)$ be a lattice r.v. with span $h$  and a set of possible values  $\Re=\lbrace a+kh, k=0,\pm1,\pm2,\dots \rbrace$. There exists a constant  $C(s)$ such that 
$$\sup_{z\in\Re}\mid\dfrac{\sigma(n)}{h}\mathsf{P}\lbrace \mathbb{R}^{\eta}(n)=z\rbrace -\dfrac{d}{du}\mathbb{W}^{(d)}(n,0,u_{z})\mid\leq C(s)\rho_d(n)+\bar{\chi}(n,\mathbb{T}(n),\pi\sigma(n)/h),\; d=4,5, $$  
where $ u_{z}=(z-\Lambda(n))/\sigma(n).  $
\end{theorem}
 We note that Theorem 3 is not valid for the  lattice $\mathbb{R}^{\eta}(n)$. Instead, in this case, we have the following Theorem 5, which follows from Theorem 4 by using  the Euler-Maclaurin summation formula. We state just a 2-term asymptotic expansion to keep the expressions simple. 

\begin{theorem}\label{th5} Let  $\text{R}^{\eta}(n)$ be a lattice r.v. with span $h$  and a set of possible values  $\Re=\lbrace a+kh, k=0,\pm1,\pm2,\dots \rbrace$. There exists a constant  $C(s)$ such that 
 $$\mid \mathsf{P}\lbrace \mathbb{R}^{\eta}(n)<u\sigma(n)+\Lambda(n)\rbrace-\mathbb{W}^{(4)}(n,0,u)-\frac{e^{-u^2/2}}{\sqrt{2\pi}}\frac{h}{\sigma(n)}S_1\Big(\frac{1}{h}( u\sigma(n)-(a-\Lambda(n))\Big)\mid   $$
$$\le C(s)\rho_{4}(n)+\bar{\chi}(n,\mathbb{T}(n),\pi\sigma(n)/h),$$
 where $ S_1(x)=x-[x]+1/2 $ is a well-known periodic function of period one (see for instance [\textbf{hat }4,p.254]),and comes up here due to the Euler - Maclaurin summation formula.
\end{theorem}

{\bf Remark 3.1.} It may be remarked that Theorems 3 , 4  and 5 can be further extended by approximation $\mathbb{W}^{(d)}(n,0,u_{z}) $  for any $d\ge 3$ but at  the expense of added complexity in the proof.\\
  
{\bf Remark 3.2.} For many variants of the GUM, the statements of this section apply and  yield reasonably accurate results; for example, conditions (ii)and (iii) of Theorem 1, which pertain on the parameters of the urn model, are generally met in a very general settings (as elaborated further in Section 5). So that all we need is to check the Lyapunov-type condition (i) for ensuring the asymptotic normality of the r.v. under consideration. We recall also that Theorems  1-5 include the GUMs consisting a finite number $ N$ of items. In this version of the GUM with $s=1$ , our theorems represent a significant refinement of Theorems 2.1 - 2.3 of Mirakhmedov (1989), Theorems 3.2 - 3.4 of Mirakhmedov et al. (2014). Namely, we eliminate the condition (B) of Mirakhmedov (1989) and the term $ \min(B_N,\sqrt{N})/M_N(0.3(B_N\kappa_3)^{-1}) $ in $ \mathcal{E}_N(\delta) $ from equation (2.10) and condition (iv) of  Theorem 3.1 of Mirakhmedov et al (2014).  The existence of such a term (condition) imposes a significant constraint on the parameters of the GUM.  \\
\indent
Nevertheless, for example, when considering a random sample drawn from a discrete distribution with fixed the mass probabilities (this is, in turn, is an example of the " balls-in-infinite boxes " variant of the GUM), condition (ii) does not hold, and hence the results above do not apply, see Section 5. This case requires an alternative analytic framework distinct from that presented in Lemma 1 must be employed (for more details, consult Gnedin et al., 2008).

{\bf Remark 3.3.} In the literature on Edgeworth expansion for sums of random variables, the results are usually presented in a form where integrals similar to $ \chi(n,a,b) $ and $ \bar { \chi }(n,a,b) $ are already estimated, imposing predefined conditions on the structure of the random variables. However, for the problems we are dealing with such conditions lead to restrictions on the parameters in particular urn models. Therefore, it would be more appropriate to leave such integrals unestimated in general statements so that they can later be assessed individually  for each version of the GUM and the r.v. $ \mathbb{R}^{\eta}(n) $, taking into account their specific properties, as demonstrated in Section 5. \\
\indent  
In the context of Remark 3.3 the following general bounds for $ \chi(n,a,b) $ turn out to be useful.
For a given r.v. $  \varsigma$ we define $  \varsigma^{*}= \varsigma- \varsigma'$, where $  \varsigma' $ is an independent copy of $  \varsigma$. Let  $ \langle a\rangle $  stands for the distance between real $a$ and integers. Set $ A(t,\tau^s)=\lbrace(t,\tau^s): \sigma^{-1}(n)a\le\mid t\mid\le b\sigma^{-1}(n),\;\mid\tau_l\mid\le \pi, l=1,\dots, s\rbrace ,$

\begin{equation}\label{16}
H_m(t,\tau^s)=E\Big\langle tg_m^{*}(\xi_m^s)+ \langle \tau^s,\xi_m^{*s}\rangle\Big\rangle^2,\quad H(n,a,b)= \inf_{(t,\tau^s)\in A(t,\tau^s)} \sum_{m=1}^{\infty}H_m\Big(\frac{1}{2\pi}t,\frac{1}{2\pi}\tau^s\Big).
\end{equation}  
Then

$$\chi(n,a,b)\le C\ln(b/a)\prod_{k=1}^s B_k \exp\left\lbrace -\frac{1}{2} H(n,a,b)\right\rbrace ,$$

and 
\begin{equation}\label{17}
\bar{\chi}(n,a,b)\le C\sigma(n)\prod_{k=1}^s B_k \exp\left\lbrace -\frac{1}{2} H(n,a,b)\right\rbrace. 
\end{equation}

These inequalities  follow from the following arguments: From formula (8) with $x^s=0$ it follows that for the ch.f. $ \varphi(n,t) $  one can write in (9) the product $ \prod_{m=1}^{\infty}\psi_m(t,\tau^s) $
instead of $ \Psi(n,t,\tau^s) $. Since $ \Theta(n,0,0)>c>0 $ inequalities  (17) follow by using the facts that $ x\le \exp\lbrace(x^2-1)/2\rbrace $ and by Lemma 4 of [Mukhin, 1991]
$$ 4H_m\Big(\frac{1}{2\pi}t,\frac{1}{2\pi}\tau^s\Big)\le 1-\mid\psi_m(t,\tau^s)\mid^2\le 2\pi^2 H_m\Big(\frac{1}{2\pi}t,\frac{1}{2\pi}\tau^s\Big). $$

\bigskip

\section{ PROOFS OF THE GENERAL RESULTS }

In addition to the notations of the previous sections we set

$$\varpi_{k}(n)=\sum_{m=1}^{\infty}E\Vert(\tilde{g}_{m}(\xi^{s}_{m}),\tilde{\xi}_{m}^{s})\Vert^{k},\quad \mu_k(n)=\rho_k(n)+\sum_{l=1}^s\varrho_{k,l}(n_l),$$ 
 
$$A_{0}(t,\tau^{s}):=\lbrace(t,\tau^{s}): \Vert(t,\tau^{s})\Vert\leq \varpi_{3}^{-1/3}(n)\rbrace,\quad  A_{1}(t,\tau^{s}):=\lbrace(t,\tau^{s}): \Vert(t,\tau^{s})\Vert\leq \varpi_{5}^{-1/5}(n)\rbrace,$$
$$ A_{2}(t,\tau^{s}):=\lbrace (t,\tau^{s}):  \vert t\vert\leq (2(s+1)^{2}\beta_{3}(n))^{-1},\; \vert \tau_{l}\vert\leq (2s(s+1)^{2}\varrho_{3,l}(n_l))^{-1}  , l=1,...,s \rbrace ,$$

$$ \nabla(t,\tau^{s}) := \Psi(n,t,\tau ^{s})-    \exp\left\lbrace-\dfrac{1}{2}\Vert(t,\tau^{s})\Vert^{2} \right\rbrace\left( 1+P_{1} (n,t.\tau^{s})+P_{2} (n,t.\tau^{s})\right) .$$

{\bf{Lemma 4.1.} {\it 

(i) \;\; If \;$\varpi_{3}(n)\leq 1 $ and $\;(t,\tau^{s})\in A_{0}(t,\tau^{s})$ then 
$$ \vert D_{t}^{j}\left( \Psi(n,t,\tau ^{s})-\exp\left\lbrace-\dfrac{1}{2}\Vert(t,\tau^{s})\Vert^{2})\right\rbrace \right) \vert  $$
$$\leq C(s)\varpi_{3}(n)\vert t\vert^{j}(1 +\Vert (t,\tau^{s})\Vert ^{3})\exp\left\lbrace-\dfrac{1}{6}\Vert(t,\tau^{s})\Vert^{2})\right\rbrace, \; j=0 , 1 .$$

 (ii) If \; $\varpi_{5}(n)\leq 1 $ and  $(t,\tau^{s})\in A_{1}(t,\tau^{s})$ then
$$ \vert D_{t}^{j}\nabla(t,\tau^{s})\vert\leq C(s)\varpi_{5}(n)\vert t\vert^{j}\Vert (t,\tau^{s})\Vert ^{5}(1 +\Vert (t,\tau^{s})\Vert ^{2})\exp\left\lbrace-\dfrac{1}{24}\Vert(t,\tau^{s})\Vert^{2})\right\rbrace, \; j=0,1 .$$
\; (iii) If\; $ (t,\tau^{s})\in A_{2}(t,\tau^{s}) $ then
$$ \vert D_{t}^{j}\Psi(n,t,\tau ^{s})\vert \leq (\vert t\vert +\vert \tau_{1}\vert+ \dot +\vert\tau_{s}\vert)e^{1/3})^{j} \exp\left\lbrace -\dfrac{1}{3}\Vert(t,\tau^{s})\Vert^{2}) \right\rbrace ,\; j=0,1 . $$ }} 

{\bf Proof.} Assume an integer $N\geq 2$ and define $\varpi_{l,N}(n)$ similarly to $\varpi_{l}(n)$ except using a finite sum over $ m=1,...,N $ rather than an infinite one. Then by Lemma 6.2 of Bhattacharya  and  Rao (1976), p.49] it holds that $\varpi_{l,N}^{1/(l-2)}\leq \varpi_{j,N}^{1/(j-2)}$ for $3\leq l\leq j ,$ and all values of $N$. Since $N$ can be chosen arbitrarily large, passing to the limit as $N\to\infty$ yields
\begin{equation}\label{eq 18} 
 \varpi_{l}^{1/(l-2)}(n)\leq \varpi_{j}^{1/(j-2)}(n),\quad 3\leq l\leq j .
 \end{equation}
Also, 
\begin{equation}\label{eq 19}
\varpi_{k}(n)\leq (s+1)^{k-1}\mu_{k}(n)\leq (s+1)^{k}\varpi_{k}(n),\;\; (E\vert\varsigma\vert^{j})^{1/j}\leq(E\vert\varsigma\vert^{\nu})^{1/\nu}, \;\; 1\leq j\leq \nu , \text{ for a r.v.}\; \varsigma .
\end{equation}
\textbf{Proof of part (i)} is similar to that of part (ii) with much simpler algebra, so we omit it. \\
 \textbf{Proof of part (ii).} Let $(t,\tau^{s})\in A_{1}(t,\tau^{s})$. Then using the second inequality of (19) we obtain
$$ \vert \psi_{m}(t,\tau^{s})-1\vert \leq \dfrac{1}{2}E\langle(t,\tau^{s}),(\tilde{g}_{m},\tilde{\xi}^{s}_{m})\rangle^{2} \leq \dfrac{1}{2}\Vert(t,\tau^{s})\Vert^{2}\varpi_{5}^{2/5}(n)\leq\dfrac{1}{2}\; ,$$

and hence, 
\begin{equation} \label{eq20}
 \vert \psi_{m}(t,\tau^{s})\vert\geq\dfrac{1}{2}. 
 \end{equation}
Further, by Taylor expansion we have
$$ \psi_{m}(t,\tau^{s})=1-\dfrac{1}{2}E\langle(t,\tau^{s}),(\tilde{g}_{m},\tilde{\xi}^{s}_{m})\rangle^{2}$$

\begin{equation} \label{eq21}
+\dfrac{i^{3}}{6}E\langle(t,\tau^{s}),(\tilde{g}_{m},\tilde{\xi}^{s}_{m})\rangle^{3}+\dfrac{i^{4}}{24}E\langle(t,\tau^{s}),(\tilde{g}_{m},\tilde{\xi}^{s}_{m})\rangle^{4} +\dfrac{\theta}{120}E\vert\langle(t,\tau^{s}),(\tilde{g}_{m},\tilde{\xi}^{s}_{m})\rangle\vert^{5} .  
\end{equation}

On using inequalities (18) we obtain: for $1<i\leq5 $, $1<j\leq5$ and $i+j\geq5 $ 
$$ \vert E\langle(t,\tau^{s}),(\tilde{g}_{m},\tilde{\xi}^{s}_{m})\rangle^{i}E\langle(t,\tau^{s}),(\tilde{g}_{m},\tilde{\xi}^{s}_{m})\rangle^{j}\vert 
\leq \Vert (t,\tau^{s})\Vert ^{i+j}(E\Vert (\tilde{g}_{m},\tilde{\xi}^{s}_{m})\Vert^{5})^{(i+j)/5}$$
$$ \leq \Vert (t,\tau^{s})\Vert ^{5}  E\Vert(\tilde{g}_{m},\tilde{\xi}^{s}_{m})\Vert^{5}. $$

Applying these inequalities we get 
\begin{equation} \label{eq22}
 \left(  \psi_{m}(t,\tau^{s})-1\right)^{2}= \dfrac{1}{4}(E\langle(t,\tau^{s}),(\tilde{g}_{m},\tilde{\xi}^{s}_{m})\rangle^{2} )^{2}+\theta c \Vert (t,\tau^{s})\Vert ^{5}  E\Vert(\tilde{g}_{m},\tilde{\xi}^{s}_{m})\Vert^{5}.  
\end{equation}
Set 
$$ h(t,\tau^{s}):= \dfrac{i^{3}}{6}\sum_{m=1}^{\infty}E\langle(t,\tau^{s}),(\tilde{g}_{m},\tilde{\xi}^{s}_{m})\rangle^{3}$$
$$+\dfrac{i^{4}}{24}\sum_{m=1}^{\infty}\left[ E\langle(t,\tau^{s}),(\tilde{g}_{m},\tilde{\xi}^{s}_{m})\rangle^{4}  -3(E\langle(t,\tau^{s}),(\tilde{g}_{m},\tilde{\xi}^{s}_{m})\rangle^{2} )^{2}\right]  $$
$$=P_{1} (n,t.\tau^{s})+P_{2} (n,t.\tau^{s})-\dfrac{1}{2}P_{1}^{2} (n,t.\tau^{s})\quad \text{(see, Step 1 of Sec.2)} .$$
Now we use (20), (21), (22) and (4) to get
$$\Psi(n,t,\tau ^{s})=\exp\left\lbrace \sum_{m=1}^{\infty}\log \psi_{m}(t,\tau^{s})   \right\rbrace $$
$$=\exp\left\lbrace \sum_{m=1}^{\infty}\left[ \left(  \psi_{m}(t,\tau^{s})-1\right)  -\dfrac{1}{2} \left(  \psi_{m}(t,\tau^{s})-1\right)^{2}+\sum_{j=3} ^{\infty}\dfrac{(-1)^{j-1}}{j}\left(  \psi_{m}(t,\tau{s})-1\right)^{j}\right] \right\rbrace   $$

\[ =\exp\left\lbrace-\dfrac{1}{2}t^{2}-\dfrac{1}{2}\sum_{l=1}^{s}\tau_{l}^{2}\right\rbrace 
 \exp\left\lbrace h(t,\tau^{s})+\theta c_{1} \parallel (t,\tau^{s})\parallel ^{5}\varpi_{5}(n) \right\rbrace .\]
Therefore,  

$$ \vert \nabla(t,\tau^{s})\vert \leq  \exp\left\lbrace-\dfrac{1}{2}t^{2}-\dfrac{1}{2}\sum_{l=1}^{s}\tau_{l}^{2})\right\rbrace  
\vert\exp\left\lbrace h(t,\tau^{s})+\theta c_{1} \parallel (t,\tau^{s})\parallel ^{5}\varpi_{5}(n) \right\rbrace $$
\begin{equation} \label{eq23}
 -1-P_{1} (n,t.\tau^{s})-P_{2} (n,t.\tau^{s})\vert .
\end{equation}
Using inequality (19) we obtain
$$ \vert \sum_{m=1}^{\infty}E\langle(t,\tau^{s}),(\tilde{g}_{m},\tilde{\xi}^{s}_{m})\rangle^{3}\sum_{m=1}^{\infty}E\langle(t,\tau^{s}),(\tilde{g}_{m},\tilde{\xi}^{s}_{m})\rangle^{4}\vert\leq \parallel (t,\tau^{s})\parallel ^{7} \varpi_{3}(n)\varpi_{4}(n)\leq \parallel (t,\tau^{s})\parallel ^{7}\varpi_{5}(n),$$

$$ \vert \sum_{m=1}^{\infty}E\langle(t,\tau^{s}),(\tilde{g}_{m},\tilde{\xi}^{s}_{m})\rangle^{3}\sum_{m=1}^{\infty}\left( E\langle(t,\tau^{s}),(\tilde{g}_{m},\tilde{\xi}^{s}_{m})\rangle^{2}\right) ^{2}\vert\leq \parallel (t,\tau^{s})\parallel ^{7}\varpi_{5}(n),$$
and, since assumption $ \varpi_{5}(n)\leq 1 ,$ for $(t,\tau^{s})\in A_{1}(t,\tau^{s})$ and $k=3,4$
$$ \vert \sum_{m=1}^{\infty}E\langle(t,\tau^{s}),(\tilde{g}_{m},\tilde{\xi}^{s}_{m})\rangle^{k}\vert\leq \Vert (t,\tau^{s})\Vert ^{k} \varpi_{5}(n)^{(k-2)/3}\leq \Vert (t,\tau^{s})\Vert^{2} \varpi_{5}(n)^{ 2(k-2)/15}\leq \Vert (t,\tau^{s})\Vert ^{2}. $$
 Applying in (23) these inequalities and that $\vert e^{u}-1-u-u^{2}/2 \vert\leq \vert u\vert^{3}e^{\vert u \vert}/6  $  we obtain 
\begin{equation} \label{eq24}
\vert\nabla(t,\tau^{s})\vert\leq C\varpi_{5}(n)(\Vert (t,\tau^{s})\Vert ^{5}(1 +\Vert (t,\tau^{s})\Vert ^{2}))\exp\left\lbrace-\dfrac{1}{24}\Vert (t,\tau^{s})\Vert ^{2}\right\rbrace .
\end{equation}
On using Taylor expansion we observe that
\begin{equation} \label{eq25}
\sum_{m=1}^{\infty}D_{t}\log\psi_{m}(n,t,\tau ^{s})=-t+D_{t}h(t,\tau^{s})+c\theta\parallel(t,\tau^{s})\parallel^{5}\varpi_{5}(n) .
\end{equation}
Next, by repeatedly using operations such as subtraction and addition, we obtain
$$ D_{t}\nabla(t,\tau^{s}) =\Psi(n,t,\tau ^{s})\sum_{m=1}^{\infty}D_{t}\log\psi_{m}(n,t,\tau ^{s}) +t \exp\left\lbrace-\dfrac{1}{2}\parallel(t,\tau^{s})\parallel^{2} \right\rbrace $$
$$\cdot\Big( 1+P_{1} (n,t.\tau^{s})+P_{2} (n,t.\tau^{s})\Big) 
-\exp\left\lbrace-\dfrac{1}{2}\parallel(t,\tau^{s})\parallel^{2} \right\rbrace \Big( D_{t}P_{1} (n,t.\tau^{s})+D_{t}P_{2} (n,t.\tau^{s})\Big)$$

$$ =\left( \Psi(n,t,\tau ^{s})-\exp\left\lbrace-\dfrac{1}{2}\parallel(t,\tau^{s})\parallel^{2} \right\rbrace \Big( 1+P_{1} (n,t,\tau^{s})+P_{2} (n,t,\tau^{s})\Big)\right) \sum_{m=1}^{\infty}D_{t}\log\psi_{m}(n,t,\tau ^{s})   $$
$$+\exp\left\lbrace-\dfrac{1}{2}\parallel(t,\tau^{s})\parallel^{2} \right\rbrace \Big( 1+P_{1} (n,t,\tau^{s})+P_{2} (n,t,\tau^{s})\Big) \Big( \sum_{m=1}^{\infty}D_{t}\log\psi_{m}(n,t,\tau ^{s})+t-D_{t}h(n,t.\tau^{s}) \Big)  $$

$$+\exp\left\lbrace-\dfrac{1}{2}\parallel(t,\tau^{s})\parallel^{2} \right\rbrace \Big( P_{1} (n,t.\tau^{s})D_{t}P_{2} (n,t,\tau^{s})+P_{2} (n,t,\tau^{s})D_{t}P_{1} (n,t,\tau^{s})$$ 
\begin{equation} \label{eq26}
+P_{2} (n,t,\tau^{s})D_{t}P_{2}(n,t,\tau^{s})-P_{1}^{2}(n,t,\tau^{s})D_{t}P_{1} (n,t,\tau^{s})- P_{2} (n,t.\tau^{s})P_{1} (n,t,\tau^{s})D_{t}P_{1} (n,t,\tau^{s}) \Big) .
\end{equation}  
 Using inequalities (19) we have
\begin{equation} \label{eq27}
\vert D^j_{t}P_{1} (n,t.\tau^{s})\vert \leq (s+1)^2 \Vert(t,\tau^{s})\Vert^{3-j} \mu_{3}(n)/2\,3^{1-j},\; j=0,1 , 
\end{equation}

\begin{equation} \label{eq28} 
\vert D^j_{t}P_{2} (n,t,\tau^{s})\vert \leq (s+1)^{3}\Vert(t,\tau^{s}\Vert^{4-j}(1+4\Vert(t,\tau^{s}\Vert^{2})\mu_{4}(n)8^{-(1-j)},\; j=0,1. 
\end{equation}
These inequalities imply that  
$ \vert D_{t}h(t,\tau^{s})\vert\leq (s+1)^{3}\Vert(t,\tau^{s}\Vert^{3}\mu_{4}(n)$ and hence by (25) we have
$ \vert\sum_{m=1}^{\infty}D_{t}\log \psi_{m} (n,t,\tau_{s})\vert \leq (\vert t\vert +c(s+1)^{3}\Vert(t,\tau^{s}\Vert^{3}\mu_{4}(n).$ Therefore,applying (19), (25), (27) and (28) to the terms of the RHS of (26) we get
 
$$\vert D_{t}\nabla_{1}(t,\tau^{s})\vert \leq C(s)\varpi_{5}(n)\vert t\vert \Vert (t,\tau^{s})\Vert^{5}(1 +\Vert (t,\tau^{s})\Vert ^{2})\exp\left\lbrace-\dfrac{1}{24}\Vert(t,\tau^{s}\Vert^{2}\right\rbrace .  $$
The calculation details are omitted. Proof of part (ii) is completed. \\

{\bf Proof of part (iii).} Let $ (t,\tau^{s})\in A_{2}(t,\tau^{s}).$ For a given r.v. $ \varsigma $ we set $ \varsigma^{*}=\varsigma-\varsigma'$, where $\varsigma'$ is an independent copy of $\varsigma .$ Noticing that $ \vert \psi_{m}(t,\tau^{s})\vert ^{2}   $ is the characteristic function of the r.vec. $(\tilde{g}_{m}^{*}(\xi^{s}_{m}),\tilde{\xi}_{1,m}^{*}, \dots ,\tilde{\xi}_{s,m}^{*})$  we obtain 
$$ \vert \psi_{m}(t,\tau^{s})\vert ^{2} \leq 1-E\langle(t,\tau ^{s}), (\tilde{g}_{m}(\xi^{s}_{m}),\tilde{\xi}_{m}^{s})\rangle ^{2}+ \dfrac{2}{3}E\vert\langle(t,\tau ^{s}), (\tilde{g}_{m}(\xi^{s}_{m}),\tilde{\xi}_{m}^{s})\rangle\vert ^{3}$$
\begin{equation} \label{eq29}
 \leq\exp\left\lbrace -E\langle(t,\tau ^{s}), (\tilde{g}_{m}(\xi^{s}_{m}),\tilde{\xi}_{m}^{s})\rangle ^{2}+ \dfrac{2(s+1)^{2}}{3}\left( \vert t\vert^{3}E\vert \tilde{g}_{m}\vert^3 +\sum_{l=1}^{s}\vert\tau_{l}\vert^{3} E\vert \tilde{\xi}_{l,m}\vert^3  \right) \right\rbrace .  
\end{equation}
Therefore, taking into account (3 ) and (4) we obtain

$$\vert\Psi(n,t,\tau ^{s})\vert \leq \exp\left\lbrace -\dfrac{t^{2}}{2}-\dfrac{1}{2}\sum_{l=1}
^{s}\tau^{2}_{l}+\dfrac{(s+1)^{2}}{3}\left( \vert t\vert^{3}\beta_{3}(n)+\sum_{l=1}^{s}\vert\tau_{l}\vert^{3} \varrho_{3,l}(n_{l})  \right)  \right\rbrace $$
$$\leq \exp\left\lbrace -\dfrac{t^{2}}{3}-\dfrac{1}{3}\sum_{l=1}^{s}\tau^{2}_{l} 
\right\rbrace .$$
Part (iii)for \textit{j} = 0 follows.\\

We have $ \vert D_{t}\psi_{m}(t,\tau^{s})\vert\leq (\vert t\vert +\vert \tau_{1}\vert E\vert \tilde{g}_{m} \tilde{\xi}_{1,m}\vert +\dots +\vert \tau_{s}\vert E\vert \tilde{g}_{m} \tilde{\xi}_{s,m}\vert ,$ and $ \sum_{m=1}^{ \infty}E\vert \tilde{g}_{m} \tilde{\xi}_{l,m}\vert\leq 1 $, by Hoglund's inequality.  Hence,
\begin{equation} \label{eq30}
\sum_{m=1}^{ \infty}\vert D_{t}\psi_{m}(t,\tau^{s})\vert\leq (\vert t\vert +\vert \tau_{1}\vert +\dots +\vert \tau_{s}\vert) .
\end{equation} 
Let $  d_m$ denote the expression of the exponent of (29). We have

$$e^{-d_m} \leq \exp\left\lbrace (s+1)(t^{2}\beta_{3}^{2/3}(n)+ \sum_{l=1}^{s}\tau_{l}^{2} \varrho _{l,3}^{2/3}(n_{l})- \dfrac{2(s+1)^{2}}{3}\Big( \vert t\vert^{3}\beta_{3}(n)+\sum_{l=1}^{s}\vert\tau_{l}\vert^{3} \varrho _{l,3}(n_{l}  \Big) \right\rbrace \leq e^{1/3},$$
since inequality $ay^{2}-by^{3}\leq 4a^{3}/27b^{2} ,$ where $y\geq 0$, $a>0$ and $b>0 $ . On using (29),(30), the last inequality  and that
\begin{equation} \label{eq31}
D_{t}\Psi(n,t,\tau^{s})=\sum_{m=1}^{\infty}D_{t}\psi_{m}(t,\tau^{s})\prod_{l=1,l\neq m}^{\infty}\psi_{l}(t,\tau^{s})
\end{equation}
we complete the proof of part (iii) for \textit{j} =1 . Lemma 4.1 is proved completely.\\

{\bf Lemma 4.2.} {\it  Let $\vert t\vert \leq\mathbf{T}(n)$ and 
  $j= 0,1 .$ 
$$ (i) \quad  if \;  \rho_3(n)\le 1 \;\; then \;\; \vert D_{t}^{j}\left( \varphi (n,t)- \exp\left\lbrace -\dfrac{t^{2}}{2}\right\rbrace \right) \vert \leq C(s)(1+\vert t\vert^{3})\rho_{3}(n)  \exp\left\lbrace -\dfrac{t^{2}}{24}\right\rbrace ,$$
$$ (ii)\qquad   if \quad   \rho_5(n)\le 1 \quad  then \quad \vert D_{t}^{j}\left( \varphi (n,t)- \exp\left\lbrace -\dfrac{t^{2}}{2}\right\rbrace W^{(5)}(n,x^{s},t)\right) \vert $$
$$\leq C(s)(1+\vert t\vert^{7})\rho_{5}(n)  \exp\left\lbrace -\dfrac{t^{2}}{24}\right\rbrace .$$ }

{\bf Proof.} Let's prove  part (ii). Part (i) can be proved similarly, using simpler algebra. Let $\vert t\vert \leq \mathbf{T}_{1}(n),$ where $ \mathbf{T}_{1}(n)=:\min \left( \varpi_{5}^{-1/5}(n)/2(s+1)^{2},\sqrt{\mathbf{M}(n)/6} \right).$ Assume that $ \mathbf{T}_{1}(n)\leq \mathbf{T}(n).$ Set
$$ A_{1}(\tau^{s})=\lbrace \tau^{s}:\Vert\tau^{s}\Vert\leq \varpi_{5}^{-1/5}/2\rbrace ,\quad \bar{A}_{1}(\tau^{s})=\lbrace \tau^{s}:\Vert\tau^{s}\Vert\geq \varpi_{5}^{-1/5}/2\rbrace, $$
$$ A_{2}(\tau^{s})=\lbrace \tau^{s}:\vert\tau_{l}\vert\leq \varrho_{3,l}^{-1}(n_{l})/2(s+1)^{2}, l=1,\dots, s\rbrace ,$$
$$  \bar{A}_{2}(\tau^{s})=\lbrace \tau^{s}:\vert\tau_{l}\vert\leq \pi B_{l}(n_{l}) , l=1,\dots, s\rbrace \setminus A_{2}(\tau^{s}) ,$$
$$V^{(5)}(n,t) =\exp\left\lbrace-\dfrac{t^{2}}{2}\right\rbrace \left( 1+ G_{1}(n,0,t)+ G_{2}(n,0,t) \right).$$
We have by (8) 

$$ \Delta(n,t)=:\vert D_{t}^{j}\left(  \Theta(n,0,t)-V^{(5)}(n,t) \right) \vert \leq(2\pi)^{-s/2}\int_{A_{1}(\tau^{s})}\vert D_{t}^{j}\nabla(t,\tau^{s})\vert d\tau^{s} $$
$$+ (2\pi)^{-s/2}\int_{\bar{A}_{1}(\tau^{s})}\vert D_{t}^{j}\left( \exp\left\lbrace-\dfrac{1}{2}\Vert(t,\tau^{s})\Vert^{2} \right\rbrace\left( 1+P_{1} (n,t.\tau^{s})+P_{2} (n,t.\tau^{s})\right)\right) \vert d\tau^{s}  $$
$$ + (2\pi)^{-s/2}\int_{A_{2}(\tau^{s})\setminus A_{1}(\tau^{s})}\vert D_{t}^{j}\Psi(n,t,\tau ^{s})\vert d\tau^{s} + (2\pi)^{-s/2}\int_{\bar{A}_{2}(\tau^{s})}\vert D_{t}^{j}\Psi(n,t,\tau ^{s})\vert d\tau^{s}$$
\begin{equation}\label{32}
:=\Delta_1(n,t)+\Delta_2(n,t)+\Delta_3(n,t)+\Delta_4(n,t).
\end{equation}
 To estimate $ \Delta_1(n,t) $ and $ \Delta_3(n,t) $   we apply  Lemma 4.1 part (ii) and (iii), respectively, and for $ \Delta_2(n,t) $ we use inequalities (27), (28), and the fact that $ \varpi_{k}(n)\leq (s+1)^{k-1}\mu_{k}(n). $ Then, after some algebraic manipulations, we obtain 
   
\begin{equation} \label{eq33}
 \Delta_1(n,t)+\Delta_2(n,t)+\Delta_3(n,t) \leq C(s)\mu_{5}(n))(1+\vert t\vert^{7})\exp\left\lbrace -\frac{t^{2}}{24}\right\rbrace  . 
\end{equation}
 Let's estimate $ \Delta_4(n,t). $  We have

$$\vert\psi_{m}(t,\tau^{s})\vert^{2} =\mid \psi_{m}(0,\tau^{s})+ E\left( \exp\left\lbrace  it\tilde{g}_{m}(\xi^{s}_{m})\right\rbrace -1\right) \exp\lbrace i\langle\tau^{s},\tilde{\xi}^{s}_{m}\rangle\rbrace\mid ^2$$
$$\le 2\Big(\mid \psi_{m}(0,\tau^{s})\mid^2+ t^2E\tilde{g}_{m}^2\Big)
\le 2\exp\left\lbrace  -\Big( 1-\mid \psi_{m}(0,\tau^{s})\mid^2 \Big)     + t^2E\tilde{g}_{m}^2\right\rbrace  ,$$
since inequality $x\leq \exp\lbrace x-1\rbrace .$ Whence 
$$\vert\Psi(n,t,\tau^{s})\vert\leq \sqrt{2}\exp\left\lbrace -\dfrac{1}{2s} \sum_{l=1}^s\sum_{m=1}^{\infty}\left( 1- \vert E\exp\lbrace i\tau_{l}\tilde{\xi}_{l,m}\vert^{2}\right) + \frac{t^2}{2}\right\rbrace $$
\begin{equation} \label{eq34}
\leq \sqrt{2}\exp\lbrace-t^2/2\rbrace\exp\left\lbrace -\dfrac{1}{4s} \sum_{l=1}^s\sum_{m=1}^{\infty}\left( 1- \vert E\exp\lbrace i\tau_{l}\tilde{\xi}_{l,m}\vert^{2}\right) \right\rbrace,
\end{equation}
since  in this case there exists $l$ such that $\vert\tau_{l}\vert\geq (2s(s+1)^{2}\varrho_{3,l}(n_{l}))^{-1} ,\;$ $\mid t\mid\le \sqrt{\mathbf{M}(n)/4s}$ and for arbitrary $l=1,\dots,s $
$$ \sum_{m=1}^{\infty}\left( 1- \vert\psi_{m}(0,\tau^{s})\vert^{2}\right)\geq \sum_{m=1}^{\infty}\left( 1- \vert E\exp\lbrace i\tau_{l}\tilde{\xi}_{l,m}\rbrace\vert^{2}\right) .  $$
 Next, if $\mid\tau_l\mid\leq(2s(s+1)^{2}\varrho_{3,l}(n_{l}))^{-1} $ then
\begin{equation} \label{eq35}
 -\sum_{m=1}^{\infty}(1-\mid\psi_{l,m}(0,\tau_l)\mid^2)\leq -\frac{1}{3}\tau_l^2  . 
\end{equation} 
It is clear that an inequality similar to (34) with factor $e$ instead of $\sqrt{e}$ still valid for the $ \prod_{l=1,l\neq m}^{\infty}\psi_{l}(t,\tau^{s}). $
By this fact, (31) , (34) and (35) we obtain  
\begin{equation}\label{36}
 \Delta_4(n,t)\le C(s)\exp\lbrace -t^2/2\rbrace\sum_{k=1}^s \sum_{(l_1,\dots,l_k)}\prod_{j=1}^k B_{l_{j}}\exp\left\lbrace -\frac{1}{4s}M_{l_{j}}(T_{l_{j}})\right\rbrace , 
 \end{equation}
where $\sum_{(l_1,\dots,l_k)}$ stands for the sum over all k-tuples $ (l_1,\dots,l_k)  $ of the set $( 1,\dots, s ).$ By (32), (33), (36) and assumption $ \rho_5(n)\le 1 ,$  we conclude that
  
\begin{equation} \label{eq37}
 \Delta(n,t)\leq C(s)\Big(\mu_{5}(n))(1+\vert t\vert^{7})+\sum_{k=1}^s B_k\exp\left\lbrace -\frac{1}{4s}M_k(T_k)\right\rbrace \Big) \exp\left\lbrace -\frac{t^{2}}{24}\right\rbrace  . 
\end{equation}

Let now $\mathbf{T}_{1}(n)\leq \vert t\vert \leq \mathbf{T}(n).$ Then we have 

$$ \Delta(n,t)\leq (2\pi)^{-s/2}\int_{A_{2}(\tau^{s})}\vert D_{t}^{j}\Psi(n,t,\tau ^{s})\vert d\tau^{s} + (2\pi)^{-s/2}\int_{\bar{A}_{2}(\tau^{s})}\vert D_{t}^{j}\Psi(n,t,\tau ^{s})\vert d\tau^{s}$$
$$+ (2\pi)^{-s/2}\int_{A_{2}(\tau^{s}) \cup \bar{A}_{2}(\tau^{s})}\vert D_{t}^{j}\left( \exp\left\lbrace-\dfrac{1}{2}\Vert(t,\tau^{s})\Vert^{2} \right\rbrace\left( 1+P_{1} (n,t.\tau^{s})+P_{2} (n,t.\tau^{s})\right)\right) \vert d\tau^{s}.  $$
On the RHS of this inequality we apply Lemma 4.1(iii), the facts used to obtain  (36), and inequalities (27),(28) to the first, second, and third terms, respectively. The subsequent calculations show that inequality (37) remains valid in this case as well. 
 We note $ V^{(5)}(n,0)\geq c>0$, some $c>0$. Therefore, from (37) it follows that  
\begin{equation} \label{eq38}
\Theta(n,0,0)=V^{(5)}(n,0)\left( 1+\theta C_{1}(s)\rho_{5}(n)\right).
\end{equation}
  
 Further, by (9), (37) and (38) we have

$$ \vert D_{t}^{j}\left( \varphi (n,t)- \exp\left\lbrace -\dfrac{t^{2}}{2}\right\rbrace W^{(5)}(n,0,t)\right) \vert = \vert D_{t}^{j}\left( \frac{V^{(5)}(n,t)}{V^{(5)}(n,0)}- \exp\left\lbrace -\dfrac{t^{2}}{2}\right\rbrace W^{(5)}(n,0,t)\right) \vert  $$
\begin{equation} \label{eq39}
 + \frac{C_{1}(s)\theta_{1}\rho_{5}(n)}{V^{(5)}(n,0)(1+C_{2}(s)\theta_{2}r_{5}(n))}\left( D_{t}^{j}V^{(5)}(n,t) + (1+\vert t\vert ^{7})\exp\left\lbrace -\frac{t^{2}}{24}\right\rbrace \right).  
\end{equation}
Remind that $ Q_{m}(n,0)$  is actually the result of the expansion of $ \left(  V^{(5)}(n,0)\right) ^{-1} $ ( see Step 3 of Section 2). Using in (39) this fact, (11), (12), (27) and (28) after some algebra we complete the proof of Lemma 4.2.\\

Let's move on to proving the theorems. Note that Theorem 1 follows straightforwardly from Theorem 2. \\
\indent  
  {\bf Proof of Theorem 3.2 and 3.3.} We assume that $ \mathbf{T}(n)\geq 1, $ otherwise the evidence becomes a little simpler. By the well-known Esseen's smoothing inequality(see e.g. Feller [ 1966, v.2, Ch. XVI, eq. (3.13) ] we have for a $\;T>0$ 
$$\Delta^{(d)}(n)\leq \frac{1}{\pi}\int_{\vert t\vert\leq T }  \dfrac{\vert\varphi(n,t)-\exp\left\lbrace -t^2/2\right\rbrace W^{(d)}(n,0,t)\vert}{\vert t\vert}dt +\frac{24}{T\sqrt{2\pi}}.  $$
By this and that $ \vert\varphi(n,t)-\exp\left\lbrace -t^2/2\right\rbrace W^{(d)}(n,0,t)\vert\leq \vert t\vert \max_{\vert u\vert \leq \vert t\vert }\vert D_t ( \varphi(n,u)-\exp\left\lbrace -t^2/2\right\rbrace W^{(d)}(n,0,u))\vert $   we get  

$$\Delta^{(d)}(n)\leq \frac{1}{\pi}\max_{\vert t\vert \leq 1 }\vert D_t ( \varphi(n,t)-\exp\left\lbrace -t^2/2\right\rbrace W^{(d)}(n,0,t))\vert $$
$$ + \frac{1}{\pi}\int_{1\leq\vert t\vert\leq \mathbf{T}(n) }  \vert\varphi(n,t)-\exp\left\lbrace -t^2/2\right\rbrace W^{(d)}(n,0,t)\vert dt 
 +\frac{1}{\pi} \chi(n,\mathbf{T}(n),T)$$
\begin{equation} \label{eq 40}
  + \frac{1}{\pi} \int_{\mathbf{T}(n)\leq\vert t\vert\leq T }  \dfrac{\vert\exp\left\lbrace -t^2/2\right\rbrace W^{(d)}(n,0,t)\vert}{\vert t\vert}dt +\frac{24}{T\sqrt{2\pi}}. 
\end{equation}
\indent
To prove Theorem 3.2 and Theorem 3.3, we consider the non-trivial cases $\rho_3(n)\le 1$ and $\rho_5(n)\le 1$, respectively. Substituting into (40) \textit{d}=3, $T=\mathbf{T}(n)$ and using Lemma 4.2(i), we arrive at the proof of Theorem 3.2. Next, substituting into (40) $ d=5 $, $ T=\beta_5^{-1}(n)$ and then, using Lemma 4.2(ii) after some fairly simple calculations, we complete the proof of Theorem 3.3.\\
\indent {\bf  Proof of Theorem 3.4}  is carried out using the classical method described, for example, in [Petrov, 1995, pp. 204-207], where the inversion formula and Lemma 4.2 are used; details are omitted.\\

\bigskip

\section{ Applications}

\bigskip

In what follows, we will use the notation of the preceding sections, keeping in mind that the distributions of the r.v.s $\xi_{lm}$ as well as the functions $f_m,\;m\ge 1 $  what are relevant for the particular version of GUM and r.v.s $\mathbb{R}^{\eta}(n)$ under consideration. We stress that in all our examples of the GUM, the distributions of the r.v.s $\xi_{lm}$ depend on an arbitrary positive parameter, which can be chosen in a suitably convenient manner. Thus, in what follows, we assume  that the parameter is chosen in advance such that $E\zeta_{l}=n_{l}$, and hence $x_{l}=0 $, $l=1,\dots,s .$ Also, the notations $\text{GUM}(s,\infty)$ and $\text{GUM}(s,N)$ mean that we are dealing with $s$ independent GUMs, having, respectively, an infinite and finite $N$ number of items (cells),, where $N\to\infty$ as $\min(n_1,\dots, n_s)\to\infty.$

 \subsection{ The infinite and finite multinomial occupancy models}

\bigskip
\indent
                                      
Consider a random allocation of particles of $s$ types over an infinite number of cells numbered $1,2,\dots$. The number of particles of the $k$-th type is $n_k$, $k=1,\dots, s$. The particles are thrown randomly one at a time and independently of each other. The probability that a particle of the $k$-th type falling into $m$-th cell is $p_{km}\geq 0 $, \ $p_{k1}+p_{k2}+...=1$. Then $\eta_{km}$ is the number of particles of the $k$-th type in the $m$-th cell, after allocating of all particles, $k= 1,\dots,s, \; m=1,2,\dots .$  This variant of $\text{GUM}(s,\infty)$ is generated by r.v.s $\xi_{km}\sim Poi(n_{k}p_{km}),$ and the r.vec. $\eta_{k}=(\eta_{k,1},\eta_{k,2},\dots )$ has the multinomial distribution with infinitely many outcomes.\\
      
The following statement  is a direct consequence of the general theorem by Mirakhmedov (2010).\\

\begin{theorem}\label{th5.1} Let $ \max_{m\in \mathbb{N}} p_{l,m}\le c <1$ and $\beta_3(n)\to 0$ then
$ \Delta^{(3)}(n)\to 0. $
\end{theorem}
We would like to emphasize that Theorem 6 is valid both for the case when the probability mass distributions $p_{l,m},\; m\in\mathbb{N}, l=1,\dots, s, $ are fixed and strongly positive, and for the case when $\max_m p_{l,m}\to 0$ as $n_l\to\infty$, $l=1,\dots, s .$ \\
\indent
It should be noted that a variant of the scheme with infinite cells and fixed strongly positive probabilities $p_{1,m},\;m\in\mathbb{N} ,$  has been considered in the literature. In this case, most of the papers is devoted to the number of occupied cells and the associated r.v.s; see, e.g., Hwang and Jonson (2008) and references therein. However, this model is not relevant to the results of Section 3. In contrast, we assume that $ \max_m p_{k,m}(n_k)\to 0 $ as $ n_k\to\infty $, $k=1,\dots, s ,$ hence  the classical multinomial models with finite number of cells is also covered. Such a model with $s=1$ was considered in the past, see, e.g. [14,27, 34],   where the number of cells containing a specified number of particles was studied. We also refer to [22, 24, 26] where for the r.v. $\mathbb{R}^{\eta}(n)$ with $s=1$  the asymptotic normality with an estimated  remainder term  and a two-term Edgeworth asymptotic expansion are presented.\\
\indent 
 We have $\;\lambda_{k,m}:= E\xi_{k,m}= var \xi_{k,m} =n_kp_{k,m}$ and $E(\xi_{k,m}-\lambda_{k,m})^4=\lambda_{k,m}+3\lambda_{k,m}^2 ,$  hence  $  \varrho_{3,k}(n_k)\leq\varrho_{4,k}^{1/2}(n_k) = (n_k^{-1}+3P_k(n_k))^{1/2} $ and (see notation of Sec.3) 
 $$ B^2_k=n_k ,\quad  T_k=\Big(2s(s+1)^2\sqrt{1+3n_kP_{2,k}}\Big)^{-1},\; \text{where} \;\; P_{l,k}=P_{l,k}(n_k)=\sum_{m=1}^{\infty}p_{k,m}^l . $$ 

 Next,  $ \vert E\exp\lbrace i\tau_k\xi_{k,m}\rbrace\vert^2=\exp\lbrace -2n_kp_{k,m}\sin^2 (\tau_k/2)\rbrace $. By this and the inequalities 
  \begin{equation}\label{41}  
  \sin^2 (\tau/2)\geq \tau^2/\pi^2 \,\, if \,\, \vert\tau\vert\leq\pi, \quad \text{and}  \quad1-e^{-u}\geq (1-e^{-c})u/c \quad if \quad \vert u\vert \leq c ,
\end{equation}   
 we obtain for $ T>0$
$$ M_k(T)\geq \sum_{m=1}^{\infty}\left( 1-\exp\left\lbrace -\frac{2n_kp_{k,m}}{\pi^2}T^2\right\rbrace \right)  \geq \sum_{m:p_{k,m}\leq\pi^2/2n_kT^2}\left( 1-\exp\left\lbrace -\frac{2n_kp_{k,m}}{\pi^2}T^2\right\rbrace \right)   $$
$$\geq\frac{2(1-e^{-1})}{\pi^2}n_kT^2 \sum_{m:p_{k,m}\leq\pi^2/2n_kT^2}p_{k,m} =\frac{2(1-e^{-1})}{\pi^2}n_kT^2 \left( 1-\sum_{m:p_{k,m}\geq\pi^2/2n_kT^2}p_{k,m}\right) $$ 
$$\geq\frac{2(1-e^{-1})}{\pi^2}n_kT^2 \left( 1-\frac{2n_kT^2}{\pi^2}\sum_{m:p_{k,m}\geq\pi^2/2n_kT^2}p_{k,m}^2\right)\geq\frac{2(e-1)}{e\pi^2}n_kT^2 \left( 1-\frac{2T^2}{\pi^2}n_kP_{2,k}\right).$$

 We let $ T=T_k$ to get  
 \begin{equation}\label{42}
M_k(T_k)\geq c_0(s)n_k/(1+3n_lP_{2,k}),\; \text{where}\;\; c_0(s)=\frac{(e-1)}{4e\pi^2 s^2(s+1)^4 }.
\end{equation}
We will assume that for each $k=1,\dots, s$
\begin{equation}\label{43}
P_{2,k}\le \frac{c_0(s)}{20s\ln n_k} .
\end{equation}
Since (42), (43) from Theorem 3.2 it follows

\begin{theorem}\label{th5.2} Assume condition (43) is fulfilled. Then
$$ \Delta^{(3)}(n)\le C(s)\Big(\beta_3(n)+\max_{k=1,\dots,s}\Big(\frac{1}{\sqrt{n_k}}+\sqrt{P_{2,k}}\Big). $$
\end{theorem}
Let now  
\begin{equation}\label{44}
  \mathbf{T}(n)=\min\Big(\frac{\beta_3^{-1}(n)}{2s(s+1)^2}, \sqrt{\frac{c_0(s)n_k}{6s(1+3n_lP_{2,k})}}\, ,\;k=1,\dots,s\Big). 
  \end{equation}
   Since $ E(\xi_{k,m}-\lambda_{k,})^6\le \lambda_{k,m}+c\lambda_{k,m}^3 $ as consequences of  Theorem 3.3 and Theorem 3.4 we obtain

\begin{theorem}\label{th 5.3} Assume condition (43) is fulfilled. Then for $  d=4\; and \; d=5 $ one has
$$ \Delta^{(d)}(n)\le C(s)\Big(\beta_d(n)+\sum_{k=1}^s \Big(\frac{1}{n_k^{(d-2)/2}}+(P_{d-2,k})^{(d-2)/2(d-3)}\Big)+\chi(n,\mathbf{T}(n),\beta_{d}^{-1}(n)). $$
\end{theorem}

\begin{theorem}\label{th5.4} Assume condition (43) is fulfilled,  $\mathbb{R}^{\eta}(n)$ is a lattice r.v. with span $h$  and a set of possible values  $\Re=\lbrace a+kh, k=0,\pm1,\pm2,\dots \rbrace$. Then for $  d=4\; and \; d=5 $ one has 
$$\sup_{z\in\Re}\mid\dfrac{\sigma(n)}{h}\mathsf{P}\lbrace \mathbb{R}^{\eta}(n)=z\rbrace -\dfrac{d}{du}\mathbb{W}^{(d)}(n,0,u_{z})\mid\leq C(s)\Big(\beta_d(n)$$
$$+\sum_{k=1}^s \Big(\frac{1}{n_k^{(d-2)/2}}+(P_{d-2,k}(n_k))^{(d-2)/2(d-3)}\Big)+\bar{\chi}(\mathbf{T}(n),\pi\sigma(n)/h),$$
where $ u_{z}=(z-\Lambda(n))/\sigma(n).  $
\end{theorem}

Let's consider application of Theorems 5.1-5.4 to particular versions of the r.v. $\mathbb{R}^{\eta}(n).$   

{\bf Example 5.1.1.} We consider the r.v. $\nu_{r}(n)$, where $r=(r_{1},...,r_{s})$, integers $r_{j}\geq0$, $r_{1}+...+r_{s}\geq1$ , which is the number of cells containing exactly $r_{k}$  particles of \textit{k}th type, $k=1,\dots,s, $ after allocation all particles of all types, viz.,
$$\nu_{r}(n)=\sum_{m=1}^{\infty}
\textbf{I}\lbrace\eta_{1,m}=r_{1},...,\eta_{s,m}=r_{s}\rbrace $$
where $\textbf{I}\lbrace A\rbrace$ denotes the indicator of the event \textit{A}. Recall, here   $\lambda_{k,m}=n_kp_{k,m}, \; \lambda_{k}=B_{k}=n_{k}$. Set $\pi_{k}(a)=e^{-a}a^{k}/k!$, $ a>0 $, 
$$\Pi_{r,m}(n)=\prod_{k=1}^{s}\pi_{r_{k}}(\lambda_{k,m}),  \quad \Lambda_{r}(n)=\sum_{m=1}^{\infty}\Pi_{r,m}(n),$$   
  
$$\tilde{\sigma}_r^2(n):=\sum_{m=1}^{\infty}var \textbf{I}\lbrace \xi_m^s=r^s \rbrace =  \sum_{m=1}^{\infty}\Pi_{r,m}(n)(1-\Pi_{r,m}(n),\quad \gamma_{l}(n)=\dfrac{1}{n_{l}}\sum_{m=1}^{\infty}(r_{l}-\lambda_{l,m})\Pi_{r,m}(n),$$

$$g_{r,m}=\textbf{I}\lbrace \xi_m^s=r^s\rbrace -\Pi_{r,m}(n)-\sum_{k=1}^{s}\gamma_{k}(n)(\xi_{k,m}-\lambda_{k,m}),$$   

$$\sigma_{r}^{2}(n):=\sum_{m=1}^{\infty}var g_{r,m}= \tilde{\sigma}_r^2(n)-\sum_{l=1}^{s}n_{l}\gamma_{l}^{2}(n),\quad \beta_{r,k}(n)=\sigma_{r}^{-k}(n)\sum_{m=1}^{\infty}E\mid g_{r,m}\mid^k.$$

As straightforward consequence of Theorem 5.2 we have  

\begin{theorem}\label{th5.5}Assume condition (43) is fulfilled. There exists a constant $C(s)$ such that 
$$\Delta_{r,s}(n):=\sup_{-\infty <u<\infty}\left| \mathsf{P}\lbrace\nu_r(n) < u\sigma_r(n)+\Lambda_r(n)\rbrace-\Phi(u)\right|$$

$$\le C(s)\Big(\beta_{r,3}(n)+\max_{k=1,\dots,s}\Big(\frac{1}{\sqrt{n_k}}+\sqrt{P_{k,2}}\Big)\Big).$$
 
\end{theorem}  
Evidently $\nu_r(n)$ is lattice r.v. whose span is 1. Let $\delta(n)=\mid\lbrace m: \lambda_{1,m}+\dots+\lambda_{s,m}\ge 1\rbrace \mid$ be the number of elements of the indicated set, and  $\mathbb{W}_{r}^{(4)}(n,0,u)$ defined as in (14) where instead of function $g_{m}$ the above function $g_{r,m}$ is putted. From Theorem 5.4 we obtain

\begin{theorem}\label{th 5.6} Assume condition (43) is fulfilled. There exist  constants $c$, $C(s)$ and $0<q<1$ such that 
$$\nabla_{r,s}(n):=\sup_{z\in \mathbb{Z}^{+}}\left|\sigma_r(n)\mathsf{P}\lbrace \nu_r(n)=z\rbrace -
\frac{1}{\sqrt{2\pi}}e^{-u^2/2}\left[ 1 +
\frac{u(3-u^2)}{6}\, \alpha_3(n)+\frac{u}{3}\,\alpha_{1,2}(n)
\right]\right|$$
$$\leq C(s)\Big(\beta_{r,4}(n) +\sum_{k=1}^s \Big(\frac{1}{n_k}+P_{2,k}\Big)\Big)+\sqrt{n_1\dots n_s}\Big(q^{\delta(n)}+\exp\lbrace - \mathbf{T}^2(n)\rbrace\Big),$$
where $ u_{z}=(z-\Lambda_r(n))/\sigma_r(n)  $, $ \mathbb{Z}^{+}=\lbrace0,1,2,\dots\rbrace ,  $ $ \mathbf{T}  $ is defined in (44),
$$\alpha_3(n)=\frac{1}{\sigma_r^3(n)}\sum_{m=1}^{\infty}\Big\{
\Pi_{r,m}(n)(1-\Pi_{r,m}(n))(1-2\Pi_{r,m}(n))-3(1-2\Pi_{r,m}(n))\Pi_{r,m}(n)\sum_{k=1}^s\gamma_k(r_k-\lambda_{k,m})$$
$$+3\Pi_{r,m}(n)\Big[\Big(\sum_{k=1}^s\gamma_k(r_k-\lambda_{k,m})\Big)^2-
\sum_{k=1}^s\gamma_k^2\lambda_{k,m}\Big]-\sum_{k=1}^s\gamma_k^3\lambda_{k,m}
\Big\},$$

$$\alpha_{1,2}(n)=\frac{1}{\sigma_r(n)}\sum_{k=1}^s\left[\frac{1}{n_k}
\sum_{m=1}^{\infty}\Pi_{r,m}(n)\Bigl((r_k-\lambda_{k,m})^2-\lambda_{k,m}\Bigr)
-\gamma_k\right].$$

\end{theorem}      
   
\textbf{Remark 5.1.} The bounds in  Theorem 5.5 and 5.6 should be viewed as universal upper estimates valid under the stated assumptions, although they are, in general, not order-optimal since the second summand in their RHS. Nevertheless, for the regular multinomial models the bounds turned out to be optimal up to constants C(s). In particular, we recall first that $\beta_{r,3}(n)\le \beta_{r,4}^{1/2} .$ In turn 
 
\[\beta_{r,4}(n)\le \frac C{\sigma_r^4(n)}\left[\tilde{\sigma}_r^2(n)+\sum_{k=1}^s\gamma_k^4n_k(1+n_kP_{2,k})\right].\]

Under assumption $\sigma_r^2(n)\asymp \tilde{\sigma}_r^2(n)$ this estimate simplifies to
\[\beta_{r,4}(n)\le C\Big(\frac 1{\sigma_r^2(n)}+\sum_{k=1}^s
P_{2,k}\Big).\]\\

The r.v. $\nu_r(n)$ has been considered previously for the case $s = 1$. In particular, in Mikhaylov (1981), the CLT with a remainder term of order $O\Big((P_{2,1})^{1/12}\Big)$ was proved under the assumption
\begin{equation}\label{45}
\liminf n^{-1}E\nu_r > 0, \quad \text{for some} \quad r \ge 2.
\end{equation}
This condition is referred to (ibid) as the central domain of the parameters $p_{1,m}, , m \ge 1$. In Mirakhmedov (1991), the results of Theorems 5.5 and 5.6 were derived for the case $s = 1$. In this case we set $n=n_1,\; p_m==p_{1,m}$, $P_2(n)=\sum_{m=1}^{\infty}p_m^2 $,  $\Lambda_r(n)=\sum_{m=1}^{\infty}\pi_r(np_m) $, an integer $r\ge 1.$ As a consequence of these results, he presented the following statements.

\begin{theorem}\label{th5.7} (Mirakhmedov (1992)). Let a fixed $r\ge 2.$ \\

 1)$\quad If \quad P_{2}(n)= o\Big(n^{-1}\Lambda_r(n)\Big)$  then 
 $$\Delta_{r,1}(n)=O\Bigg(\sqrt{\frac{1}{\Lambda_r(n)}+P_{2}(n)}\;\Bigg).$$
 2) $\quad If \quad P_{2}(n)\le C \Lambda_r^{-1}(n)$ then there exists a constant $\; 0<q<1\;$ such that
 $$\nabla_{r,1}(n)=O\Big(\frac{1}{\Lambda_r(n)} + q^{\delta(n)}\Big). $$
\end{theorem}

\begin{theorem}\label{th5.8} (Mirakhmedov (1992)). Let condition (45) is fulfilled for some fixed $r\ge 2$. Then for for any $r\ge 1$ \\
$$ 1)\quad\Delta_{r,1}(n)=O\Bigg(\sqrt{\frac{1}{n}+P_{2}(n)}\;\Bigg).$$
$$2)\quad \nabla_{r,1}(n)=O\Big(\frac{1}{n} +P_{2}+ q^{\delta(n)}\Big),\quad some\;\; q\in (0,1). $$
\end{theorem}

\bigskip
\textbf{Proof} of Theorem 5.6. Due to Theorem 5.4 we only need to estimate the $ \chi(n,\mathbf{T}(n),\pi \sigma(n)). $
Similarly to Lemma 1.1 for   $ \varphi _r(n,t):=E\exp\left\lbrace it\nu_{r}(n)/\sigma_{r}(n)\right\rbrace $ it holds  
\begin{equation} \label{eq46}
 \varphi_r (n,t)=\dfrac{1}{(2\pi)^{s/2}\prod_{l=1}^{s}\sqrt{2\pi } \pi_{n_{l}}(n_{l})}\int_{\mid\tau_1\mid\le\pi,\dots,\mid\tau_s\mid\le\pi}e^{i\langle \tau^{s},\sqrt{n}^{s} \rangle }\prod_{m=1}^{\infty}\psi_{m}(t,\tau^{s}) d\tau^{s}, 
 \end{equation}

where $\psi_{m}(t,\tau^{s}):=E\exp\left\lbrace it\sigma_{r}^{-1}(n)\textbf{I}\lbrace \xi_{m}^{s}=r^{s}\rbrace +i\langle \tau^{s},\xi_{m}^{s}\rangle\right\rbrace .$
We have
$$ \psi_{m}(t,\tau^{s})=\Pi_{r,m}(n)\exp\left\lbrace it\sigma_{r}^{-1}(n) +i\langle \tau^{s},r^{s}\rangle\right\rbrace + \prod_{l=1}^{s}\exp\left\lbrace \lambda_{m,l}\left( e^{i\tau_l} -1\right) \right\rbrace $$  
$$- \Pi_{r,m}(n)\exp\left\lbrace i\langle \tau^{s},r^{s}\rangle\right\rbrace $$

$$= \prod_{l=1}^{s}\exp\left\lbrace \lambda_{m,l}\left( e^{i\tau_l} -1\right) \right\rbrace  
+\left( \exp\left\lbrace it\sigma_{r}^{-1}(n)\right\rbrace -1\right)  \Pi_{r,m}(n)\exp\left\lbrace i\langle \tau^{s},r^{s}\rangle\right\rbrace $$

$$= \exp\left\lbrace -\sum_{l=1}^{s}\lambda_{m,l} \right\rbrace   \exp\left\lbrace \sum_{l=1}^{s}\lambda_{m,l} e^{i\tau_{l}} \right\rbrace $$ 
$$+\left( \exp\left\lbrace it\sigma_{r}^{-1}(n)\right\rbrace -1\right)  \prod_{k=1}^{s}\dfrac{\lambda_{k,m}^{r_{k}}}{r_{k}!}\exp\left\lbrace-\sum_{l=1}^{s}\lambda_{m,l}\right\rbrace  \exp\left\lbrace i\langle \tau^{s},r^{s}\rangle\right\rbrace $$

$$= \exp\left\lbrace -\sum_{l=1}^{s}\lambda_{m,l} \right\rbrace  \Big[  \exp\left\lbrace \sum_{l=1}^{s}\lambda_{m,l} e^{i\tau_l} \right\rbrace 
+\left( \exp\left\lbrace it\sigma_{r}^{-1}(n)\right\rbrace -1\right)  \prod_{k=1}^{s}\frac{\lambda_{k,m}^{r_{k}}}{r_{k}!} \exp\left\lbrace i\langle \tau^{s},r^{s}\rangle\right\rbrace \Big] .$$
  
Next, $\vert \psi_{m}(t,\tau^{s})\vert\leq 1$, and our goal is to show that equality is achieved only when $ t=0, \tau_{1}=0, \dots,\tau_{s}=0 $. From the last equality we obtain
$$\exp\left\lbrace \sum_{l=1}^{s}\lambda_{m,l} \right\rbrace -\exp\left\lbrace \sum_{l=1}^{s}\lambda_{m,l} \right\rbrace \vert \psi_{m}(t,\tau^{s})\vert   $$ 
 $$= \exp\left\lbrace \sum_{l=1}^{s}\lambda_{m,l} \right\rbrace - \vert  \exp\left\lbrace \sum_{l=1}^{s}\lambda_{m,l} e^{i\tau_l} \right\rbrace  +\left( \exp\left\lbrace it\sigma_{r}^{-1}(n)\right\rbrace -1\right)  \prod_{k=1}^{s}\frac{\lambda_{k,m}^{r_{k}}}{r_{k}!} \exp\left\lbrace i \tau_lr_l \right\rbrace \vert  $$

$$= \exp\left\lbrace \sum_{l=1}^{s}\lambda_{m,l} \right\rbrace - \vert  \exp\left\lbrace \sum_{l=1}^{s}\lambda_{m,l} e^{i\tau_l} \right\rbrace   -  \prod_{l=1}^{s}\dfrac{\lambda_{l,m}^{r_{k}}}{r_{l}!} \exp\left\lbrace i \tau_lr_l \right\rbrace $$
$$  +  \exp\left\lbrace it\sigma_{r}^{-1}(n)\rangle\right\rbrace \prod_{k=1}^{s}\frac{\lambda_{k,m}^{r_{k}}}{r_{k}!} \exp\left\lbrace i \tau_lr_l \right\rbrace \vert  $$

$$= \prod_{l=1}^{s}\sum_{k=0}^{\infty}\frac{\lambda_{l,m}^{k}}{k!} - \vert  \prod_{l=1}^{s}\sum_{k=0}^{\infty}\frac{\lambda_{l,m}^{k}}{k!} \exp\left\lbrace i \tau_lk \right\rbrace   -  \prod_{l=1}^{s}\frac{\lambda_{l,m}^{r_{k}}}{r_{l}!} \exp\left\lbrace i \tau_lr_l \right\rbrace $$
$$  +  \exp\left\lbrace it\sigma_{r}^{-1}(n)\rangle\right\rbrace \prod_{k=1}^{s}\dfrac{\lambda_{k,m}^{r_{k}}}{r_{k}!} \exp\left\lbrace i \tau_lr_l \right\rbrace \vert  $$
\begin{equation} \label{eq47}
\geq \prod_{l=1}^{s}\sum_{k=0}^{\infty}\dfrac{\lambda_{l,m}^{k}}{k!} - \vert  \prod_{l=1}^{s}\sum_{k=0}^{\infty}\dfrac{\lambda_{l,m}^{k}}{k!} \exp\left\lbrace i \tau_lk \right\rbrace   -  \prod_{l=1}^{s}\frac{\lambda_{l,m}^{r_{k}}}{r_{l}!} \exp\left\lbrace i \tau_lr_l \right\rbrace \vert  - \vert \prod_{k=1}^{s}\frac{\lambda_{k,m}^{r_{k}}}{r_{k}!} \exp\left\lbrace i \tau_lr_l \right\rbrace \vert  
\end{equation}  
$$= \prod_{l=1}^{s}\sum_{k=0}^{\infty}\frac{\lambda_{l,m}^{k}}{k!} - \vert \prod_{k=1}^{s}\frac{\lambda_{k,m}^{r_{k}}}{r_{k}!} \exp\left\lbrace i \tau_lr_l \right\rbrace \vert - \vert  \prod_{l=1}^{s}\sum_{k=0, k\not = r_{l}}^{\infty}\frac{\lambda_{l,m}^{k}}{k!} \exp\left\lbrace i k\tau_l \right\rbrace \vert    $$
\begin{equation} \label{eq48}
= \prod_{l=1}^{s}\sum_{k=0, k\not = r_{l}}^{\infty}\dfrac{\lambda_{l,m}^{k}}{k!}  - \vert  \prod_{l=1}^{s}\sum_{k=0, k\not = r_{l}}^{\infty}\frac{\lambda_{l,m}^{k}}{k!} \exp\left\lbrace i k\tau_l \right\rbrace \vert \geq 0  . 
\end{equation} 

 Note that the inequality $ \vert z_{1}+z_{2}\vert\leq \vert z_{1}\vert+\vert z_{2}\vert $ becomes equality only if $\arg z_{1}=\arg z_{2}.$ Therefore, inequality (47) turns into equality only with a single choice of \textit{t} depending on $\tau_{1},\dots \tau _{s} $. On the other hand for $ \lambda_{1,m}+\dots+\lambda_{s,m}>c>0 $ inequality (48) becomes equality only if $ \tau_{1}=\dots=\tau_{s}=0 $. Then  (46) becomes equality only if $ t=0$. That is $ \vert \psi_{m}(t,\tau^{s})\vert=1  $ only if $\tau_{1}=\dots =\tau _{s}=t=0 .$ Thus, since $ \psi_{m}(t,\tau^{s})$ is continuous function we conclude that if at least one of the inequalities $\vert t\vert \geq c_{0}$, $\vert \tau_{1}\vert \geq c_{1}$, ... ,$\vert \tau_{s}\vert \geq c_{s}$ hold for some positive constants $ c_{0},\dots c_{s} $ then there exists a constant   
 \textit{q}  such that 
 \begin{equation} \label{eq49}
  \vert  \psi_{m}(t,\tau^{s})\vert \leq q<1 .
  \end{equation} 

 On the other hand $ \psi_{m}(t,0)  $ is the ch.f. of a Bernoulli r.v. with the probability of success equal to $ p_{m}:=\Pi_{r,m}(n).$ Therefore, 
 $$ \vert \psi_{m}(t,0)\vert^{2}= 1-p_{m}(1-p_{m})\sin ^{2}\frac{t}{2\sigma_{r}(n)}\leq \exp\left\lbrace -p_{m}(1-p_{m})\sin^{2}\frac{t}{2\sigma_{r}(n)}\right\rbrace \leq \exp\left\lbrace -\frac{p_{m}(1-p_{m})}{\pi^{2}\sigma^{2}_{r}(n)}t^{2}\right\rbrace , $$ 
 for $\vert t\vert\leq\pi\sigma_{r}(n)$, since $\sin^{2}u\geq u^{2}/\pi^{2}$ for $\vert t \vert\leq\pi $. Hence, within interval $\vert t \vert\leq \pi\sigma_{r}(n)$ we have
  $$ \prod_{m=1}^{\infty}\vert \psi_{m}(t,0)\vert \leq \exp\left\lbrace -\dfrac{t^{2}}{\pi^{2}\sigma^{2}_{r}(n)}\sum_{m=1}^{\infty}\Pi_{r,m}(n)(1-\Pi_{r,m}(n))\right\rbrace .$$ 
 In turn, this implies due to continuity of function $\psi_{m}(t,\tau^{s})$ that there exists a small enough $\epsilon >0$ such that for $\Vert \tau^{s}\Vert \leq \epsilon $ it holds
\begin{equation} \label{eq50} 
\prod_{m=1}^{\infty}\vert \psi_{m}(t,\tau^{s})\vert \leq \exp\left\lbrace -\dfrac{t^{2}}{2\pi^{2}\sigma^{2}_{r}(n)}\sum_{m=1}^{\infty}\Pi_{r,m}(n)(1-\Pi_{r,m}(n))\right\rbrace \leq \exp\left\lbrace -\dfrac{t^{2}}{2\pi^{2}}\right\rbrace.
\end{equation}
Further, using Stirling's formula one can observe that $ \sqrt{2\pi } \pi_{n_{l}}(n_{l})\sim \sqrt{n_l} .$ On using this fact from (45)-(49) we get
 \begin{equation} \label{eq51} 
\bar{\chi}(c\mathbf{T}(n),\pi\sigma_r(n))\le C \sqrt{n_l \dots n_s}\Big(q^{\delta(n)}+\exp\lbrace-c\mathbf{T}^2(n)\rbrace \Big) .
\end{equation}
Theorem 5.6 follows by (50), (51) and Theorem 5.4.  \\   

\textbf{Proof of the bound for $\beta_{r,4}(n).$} We have  

\[\sum_{m=1}^{\infty}Eg_{r,m}^4\le 8\sum_{m=1}^{\infty}\Big[E(\textbf{I}\lbrace \xi_m^s=r^s\rbrace-\Pi_{r,m})^4+E\!\Big(\sum_{k=1}^s\gamma_k (\xi_{k,m}-\lambda_{k,m})\Big)^4 \Big]$$
$$\le 8\Big[\sum_{m=1}^{\infty}\Pi_{r,m}(1-\Pi_{r,m})+s^3\sum_{k=1}^s\gamma_k^4\sum_{m=1}^{\infty}(\lambda_{k,m}+3\lambda_{k,m}^2)\Big] =8\Big[\tilde{\sigma}_r^2(n)+s^3\sum_{k=1}^s\gamma_k^4n_k(1+3n_kP_{2,k})\Big].\]  
The bound for the $\beta_{r,4}(n)$ follows by dividing the last inequality by \(\sigma_r^4(n)\). Next, if $\sigma_r^2(n)\asymp \tilde{\sigma}_r^2(n)$ then 
$\sum_{k=1}^sn_k\gamma_k^2=O\left(\sigma_r^2(n)\right),$ also $\mid\gamma_r(n)\mid\le 2,$ whence the second statement of Remark 5.1 follows. 

\bigskip
\textbf{Example 5.1.2.} Let $s = 1$. We rephrase the model as follows. Let a population consist of an infinite number of symbols labeled by $\mathbb{N} = {1,2,\dots}$. From this population, a sample of size $n$ is drawn according to a sampling scheme with replacement. 
Let $p_m$ be the probability that the $m$-th symbol appears in a single draw, where $p_1 + p_2 + \dots = 1$, and the r.v. $\eta_m$ denote the frequency of the $m$-th symbol in the sample. Then the r.vec. $\eta(n) = (\eta_1, \eta_2, \dots)$ defines a multinomial model with infinitely many outcomes (cells); that is, it corresponds to the urn model described in this section for $s = 1$. Several important functionals of $\eta(n)$ that arise in many contexts include:  $\mu_r(n)$- the number of symbols that appear in the sample exactly $r \ge 1$ times, $\mu(n)$ - the number of distinct symbols in the sample, and $C(n)$ - the number of collisions (i.e., the number of draws of symbols already observed in the sample). The r.v. $\mu_r(n)$ already considered in Theorems 5. and 13 above. 
Also, since $C(n) =n-\mu(n) $, findings regarding $\mu(n)$ directly apply to $C(n)$ as well. Therefore, below we concentrate on the r.v.s  $\mu(n).$  

 The r.v. $\mu(n)$ has been extensively studied within the framework of the model with fixed pmf $p_m , m\ge 1$ Specifically, it was proved that $var \mu(n)\to \infty $ as $n\to\infty$ is a sufficient condition for the asymptotic normality of $\mu(n)$. For more information see, [10] and references therein.\\
 \indent
We have the following representation 
$$  \mu(n)=\sum_{m=1}^{\infty}\textbf{I}\lbrace \eta_m >0 \rbrace .$$
Throughout this example we set $\lambda_m =np_m,$ $ \xi_m \sim Poi\lbrace\lambda_m\rbrace, $
$$\Lambda_{\mu}(n)=\sum_{m=1}^{\infty}(1-e^{-\lambda_m}), \quad \gamma_{\mu}(n)=\frac{1}{n}\sum_{m=1}^{\infty}\lambda_m e^{-\lambda_m}, $$
 
 $$\tilde{\sigma}^2_{\mu}(n)=\sum_{m=1}^{\infty}e^{-\lambda_m}(1-e^{-\lambda_m}),\quad \sigma^2_{\mu}(n)=\tilde{\sigma}^2_{\mu}(n) -n\gamma^2_{\mu}(n),$$

\begin{theorem}\label{th 5.9} If $\sigma^2_{\mu}(n)\to \infty$ then  
$$\Delta_{\mu}(n)=\sup_{-\infty<u<\infty}\mid \textbf{P}\lbrace \mu(n)<u\sigma_{\mu}(n)+\Lambda_{\mu}(n)\rbrace -\Phi(u)\mid\to 0.$$
\end{theorem}
Note that this theorem is valid for both cases when the probabilities $p_m, m\ge 1$ are fixed or $p_m=p_m(n)\to 0, m\ge 1.$

\begin{theorem}\label{th 5.10} Let 
$$ P_2(n)=\sum_{m=1}^{\infty}p^2_m \le \frac{e-1}{20\ln n}. $$
 Then
$$\Delta_{\mu}(n)\le c\Big(\frac{1}{\sigma_{\mu}(n)}+\sqrt{P_2(n)}\Big) .$$
\end{theorem}
    
\begin{theorem}\label{th 5.11} Assume $\tilde{\sigma}_{\mu}^2(n)P_2(n)\le c $. There exist positive constants $c$, $C$  such that 
$$\sup_{z\in \mathbb{Z}^{+}}\left|\sigma_\mu(n)\mathsf{P}\lbrace \mu(n)=z\rbrace -
\frac{1}{\sqrt{2\pi}}e^{-u^2/2}\left[ 1 +
\frac{u(3-u^2)}{6}\, \alpha_3(n)+\frac{u}{3}\,\alpha_{1,2}(n)
\right]\right|$$
$$\leq C\Big(\frac{1}{\sigma^2_{\mu}(n)} +\sqrt{n}\sigma_{\mu}(n) \exp\lbrace - c \bar{D}(n)\rbrace\Big),$$
where $ u_{z}=(z-\Lambda_r(n))/\sigma_{\mu}(n)  $, $ \mathbb{Z}^{+}=\lbrace0,1,2,\dots\rbrace ,  $ ,
$$ \bar{D}(n)= \sum_{m=1}^{\infty}\pi_1(\lambda_m)\pi_0(\lambda_m)= \sum_{m=1}^{\infty}\lambda_m\min(1,\lambda_m)e^{-2\lambda_m},  $$

$$\alpha_3(n)= -\frac{1}{\sigma_{\mu}^3(n)}\sum_{m=1}^{\infty}\Big[e^{-\lambda_m}(1-e^{-\lambda_m})(1-2e^{-\lambda_m})$$
$$-3\gamma_{\mu}(n)\lambda_m e^{-\lambda_m}(1-2e^{-\lambda_m})
-3\gamma_{\mu}^2(n)\lambda_m(1-\lambda_m)e^{-\lambda_m}+\lambda_m\gamma_{\mu}^3(n)\Big],$$

$$ \alpha_{1,2}(n):= \frac{1}{n\sigma_{\mu}(n)}\sum_{m=1}^{\infty}\lambda_m(1-\lambda_m)e^{-\lambda_m} -\frac{\gamma_{\mu}(n)}{\sigma_{\mu}(n)}.$$

\end{theorem}    
 
\textbf{Remark 5.2.} A sufficient condition  the exponential term in rhs of Theorem 19 to be small enough  is that  the total probability mass carried by symbols with $p_m\asymp n^{-1}$  should be at least of order $(\ln n)/n.$\\
  
\textbf{Proofs of Theorems 5.9 - 5.11.} Set  $g_m(\xi_m):= \textbf{I}\lbrace \xi_m>0\rbrace-(1-e^{-\lambda_m})-\gamma_{\mu}(n)(\xi_m-\lambda_m)= -[\textbf{I}\lbrace \xi_m=0\rbrace-e^{-\lambda_m}]-\gamma_{\mu}(n)(\xi_m-\lambda_m)$.
We have $\gamma_{\mu}(n)\le 1  $ and
\begin{equation}\label{53 }
\beta_4(n)=\sigma^{-4}_{\mu}(n)\sum_{m=1}^{\infty} E\mid g_m(\xi_m)\mid^4\le 8\sigma^{-4}_{\mu}(n)\Big( \tilde{\sigma}_\mu ^2(n)+\gamma_{\mu}^4(n)n(1+3nP_2(n))\Big).
\end{equation}

In the case  $p_m, m\ge 1$ fixed by Lemma 4 of Gnedin et al (2007)[7]  $\sigma_\mu ^2(n)\sim \tilde{\sigma}_\mu ^2(n),$ which implies $\, n\gamma_\mu^2(n)=o(\sigma_\mu^2(n)). $ Therefore, by (53) we have  $\; \beta_4(n)\le 8(2 \sigma_\mu^{-2}(n)+o(1)).$

Next, in the case $\; p_m=p_m(n)\to 0, m\ge 1\;$  by Proposition 4.3 of Hwang and Jonson (2008)[11]   $\,\sigma_\mu ^2(n)\asymp \tilde{\sigma}_\mu ^2(n)$, and hence, $ n\gamma_\mu^2(n)\le c\sigma_\mu^2(n). $ Using this by (53) we obtain 
\begin{equation}\label{53}
\beta_4(n)\le C\Big(\frac{1}{\sigma_{\mu}^2(n)} + P_2(n)\Big).
\end{equation}
Thus, Theorem 5.9 follows from Theorem 5.1, whereas Theorem 5.10 follows from Theorem 5.2 and inequality $\beta_3(n)\le(\beta_4(n))^{1/2} .$\\
Theorem 5.11 we will derive as consequence of Theorem 5.4 assuming $s=1$ and $d=4$ . To end this we only need to estimate $\bar{\chi}(n,\mathbf{T},\pi\sigma_{\mu}(n)), $ since (53) and the assumption $ \sigma^2_{\mu}(n)P_2(n)\le c $. Here $ \mathbf{T}$ is as in (44) with $s=1.$ 
Note that $\sigma_{\mu}(n) $ can't exceed $n$ and that under the assumption of Theorem 5.11 one has $\mathbf{T}\sigma_{\mu}^{-1}(n)\ge c_0>0.$ \\
 We have ( see, notation in  Remark 3.3)
$$ P\lbrace g_m^{*}(\xi_m)= (1-\gamma_{\mu}(n)), \xi_m^{*}=1\rbrace =\pi_1(\lambda_m)\pi_0(\lambda_m),\quad  P\lbrace g_m^{*}(\xi_m)= -\gamma_{\mu}(n)), \xi_m^{*}=1\rbrace =\pi_2(\lambda_m)\pi_1(\lambda_m).$$
Therefore,
$$ H_m\Big(\frac{t}{2\pi},\frac{\tau}{2\pi}\Big)=E\Big\langle\frac{t}{2\pi}g_m^{*}(\xi_m) +\frac{\tau}{2\pi}\xi_m^{*}\Big\rangle^2 \ge \lambda_m\min(1,\lambda_m^2)e^{-2\lambda_m}\Big(\Big\langle-\frac{t}{2\pi}\gamma_{\mu}(n))+\frac{\tau}{2\pi}\Big\rangle^2 $$

\begin{equation}\label{54}
+\Big\langle\frac{t}{2\pi}(1-\gamma_{\mu}(n))+\frac{\tau}{2\pi}\Big\rangle^2\Big)\ge c \,\lambda_m\min(1,\lambda_m^2)e^{-2\lambda_m},
\end{equation}
since for $ c_0\le\mid t\mid\le\pi ,\;\; \mid\tau\mid\le\pi $  at least one of the summands in the outer brackets is bounded away from zero. 
Finally, by inequality (17) and (55)  we obtain 
$$\bar{\chi}(n,\mathbf{T},\pi\sigma_{\mu}(n))\le  C\exp\left\lbrace -c\sum_{m=1}^{\infty}\lambda_m\min(1,\lambda_m^2)e^{-2\lambda_m} \right\rbrace . $$
Theorem 5.11 follows.\\

In what follows we consider examples of GUMs with a finite number of cells (boxes) $N$, where $N = N(n) \to \infty$ as $\min_k n_k \to \infty$. We use the notation introduced above, adding the subscript $N$ whenever necessary.\\

\textbf{Example 5.1.3.} In the random allocation scheme considered in this subsection we assume $p_{l,m}>0$ for $m=1,...,N$ and $p_{l,m}=0 $ for all $ m>N $, $l=1,...,s $,  $N(n_{1},...,n_{s})\rightarrow\infty$ as $\min(n_{1},...,n_{s})\rightarrow\infty$. That is $n_{l}$ particles of \textit{l}th type are allocated into a finite number \textit{N} of cells by a classical multinomial allocation scheme, see Kolchin et al (1976)[15]. We are interested here in r.v. $\mu_{0}(n)$, the number of empty cells after allocation of particles of all types. In the context of a random sampling design with replacement from a finite population, the r.v. $\mu_{0}(n)$ is the number of elements of the  population that do not appear in any of $s$ independent samples of sizes $ n_1,\dots, n_s .$ \\
 We write 
 $$\mu_{0}(n)=\sum_{m=1}^{N}\textbf{I}\lbrace\eta_{1,m}+...+\eta_{s,m}=0\rbrace$$ We assume that 
\begin{equation}\label{55}
Np_{l,m}\leq C .
\end{equation}
  Recall notation $\lambda_{k,m}= n_kp_{k,m},\;\;\lambda_{m}=\lambda_{i,m}+...+\lambda_{s,m}$ and set
$$\Lambda_{0}(n)=\sum_{m=1}^{N}\exp{\left( -\lambda_{m}\right) }, \quad \gamma_{0,l}(n)=\dfrac{1}{n_{l}}\sum_{m=1}^{N}\lambda_{l,m}\exp{\left( -\lambda_{m}\right) },$$ 
$$\sigma_{0}^{2}(n)=\sum_{m=1}^{N}\left[ \exp{\left( -\lambda_{m}\right) }\left( 1-\exp{\left( -\lambda_{m}\right) }\right) -N^{-1}\sum_{j=1}^{s}n_{j}\gamma_{0,j}^{2}(n)\right]$$ 

\begin{theorem}\label{th5.12} There exist a constant $C(s)$ such that
 $$\sup_{-\infty <u <\infty}|Pr\lbrace\mu_{0}(n)<u\sigma_{0}(n)+\Lambda_{0}(n)\rbrace-\Phi(u)|\leq C(s)\left( \frac{1}{\sigma_{0}(n)}+\sum_{j=1}^{s}\dfrac{1}{\sqrt{n_{j}}}\right). $$
 \end{theorem}
Below we use the above notation  $ \beta_l(n),\; $ $ \alpha_3(n) $ and $ \alpha_{1,2}(n) $ where the function $g_{m}$ is replaced by $g_{0,m}=\textbf{I}\lbrace\xi_{1,m}+...+\xi_{s,m}=0\rbrace-\exp{\left( -\lambda_{m}\right) }-\sum_{j=1}^{s}\gamma_{j}(n)(\xi_{j,m}-\lambda_{j,m}).$ \\

\begin{theorem}\label{th 5.13}There exist a constant $C(s)$ and $ q\in(0,1)$ such that
$$\sup_{z\in \mathbb{Z}^{+}}\left|\sigma_0(n)\mathsf{P}\lbrace \mu_0(n)=z\rbrace - \frac{1}{\sqrt{2\pi}}e^{-u^2/2}\left[ 1 +
\frac{u(3-u^2)}{6}\, \alpha_3(n)+\frac{u}{3}\,\alpha_{1,2}(n)
\right]\right|$$
$$\leq C(s)\Big(\frac{1}{\sigma^2_0(n)} +\sum_{k=1}^s \frac{1}{n_k}+ q^N\Big),$$
where $ u_{z}=(z-\Lambda_0(n))/\sigma_0(n)  ,$ 

$$\alpha_3(n)= \frac{1}{\sigma_0^3(n)}
\Bigg\{ \sum_{m=1}^N e^{-\lambda_m}(1-e^{-\lambda_m})(1-2e^{-\lambda_m})
+3\sum_{j=1}^s n_j\gamma_{0,j}^2(n)
-6\sum_{m=1}^N e^{-2\lambda_m}  \sum_{j=1}^s \gamma_{0,j}(n)\lambda_{j,m}$$

$$+3\sum_{m=1}^N e^{-\lambda_m} \left[\left(\sum_{j=1}^s\gamma_{0,j}(n)\lambda_{j,m}\right)^2 -\sum_{j=1}^s \gamma_{0,j}^2(n)\lambda_{j,m}\right]
-\sum_{j=1}^sn_j\gamma_{0,j}^3(n)\Bigg\},$$

$$ \alpha_{1,2}(n)=\frac{1}{\sigma_0(n)}\left[\sum_{m=1}^N\sum_{k=1}^s
\frac{\lambda_{k,m}^2}{n_k}e^{-\lambda_m}-2\sum_{k=1}^s\gamma_{0,k}(n)
\right]. $$

\end{theorem} 

 \textbf{Proof} of Theorem 5.12 and 5.13. Under the condition (56) we have $ P_{2,k}\le cN^{-1}.\; $ Next, using Holder's inequality one can observe that $ \beta_l(n)\ge N^{-(l-2)/2} ,\; l\ge 3.$ and by inequality (13) of Mirakhmedov and Saidova (1998) $ \;\beta_l(n)\le c(l)\sigma_0^{-(l-2),}(n),\;\; l\ge 3,$ whence we also have $ (\sigma_0(n)\beta_3(n))^{-1}\ge c_0 $ for some $c_0>0.$  Therefore, arguments similar to those used in proving the inequality (49) (see also proof of Lemma 2 of Kolchin et al. (1976, pp. 132 -133)) allow us to get that there exists a constant $ q\in (0,1) $ such that $ \bar{\chi}( n,\mathbf{T},\pi\sigma_0(n))\le  q^N. $ Now, Theorem 5.12 and Theorem 5.13 follows from Theorem 5.2 and Theorem 5.4, respectively.\\

\subsection{ Random sampling scheme without replacement}

\bigskip
\indent
  
\textbf{Example 5.2.1. \;} We consider a random allocation of particles, grouped into sets, among 
$N $ cells indexed by $ 1,\dots,N. $ The sets are indexed by $1,\dots, s $ and $j$th set contains $n_j$ particles. Particles are allocated one set at a time, independently across sets, with the restriction that no two particles from the same set may occupy the same cell. For each $ i =1,2,...,s $ all $C^{n_i}_N$ possible allocations of the $i$th set are equally likely,  with probability $(C^{n_i}_N)^{-1}.$  Finally, we assume $N=N(n_1,\dots,n_s)\to\infty$ as $\min (n_1,\dots,n_s)\to\infty.$\\
The allocation of the \textit{l}th set of particles is determined by the random vector  $ (\eta_{l,1},\dots,\eta_{l,N}) $, where  $ \eta_{l,m}=1$ if a particle of the \textit{l}th set falls into the cell with index \textit{m} , else $ \eta_{l,m}=0$. This scheme is a version of GUM  where for  $l=1,\dots, s,$ 
\begin{equation}\label{57}
\xi_{l.m}\sim Bi(1,p_l), \quad m=1,\dots,N ; \;   \textbf{P}\lbrace \xi_{l,m}=0\rbrace=1 \;\text{if} \; m \ge N+1; \;  p_l=n_l/N ,\; q_l=1-p_l. 
 \end{equation}   
 We are interested in the r.v. $\mu_0(N)$, the number of empty cells left after all the  sets of particles have been allocated. We note that this allocation scheme corresponds to a simple random sampling scheme without replacement from a finite population of size $ N$.  Any samples of size $n_k$  have the same probability $(C^{n_k}_N)^{-1}$ of being drawn. The samples of sizes $n_1, \dots, n_s $, respectively, are drawn independently of each other. Then the r.v. $\mu_0(N)$ is the number of population elements that do not appear in the samples. We have

$$ \mu_0(N)=\sum_{m=1}^N\textbf{I}\lbrace \eta_{1,m}=0,\dots,\eta_{s,m}=0\rbrace .$$

The r.v. $ \mu_0(N) $ has been studied in the past by many authors. Below  the latest results related to the topic of this article are presented.
By Vatutin and Mikhailov (1982), the integral and local central limit theorems  as well as all possible Poisson limit theorems are proved; in all cases, an optimal estimate of the rate of convergence is given also.  In particular they have proved the following
  
\begin {theorem} \label{th 5.14 } Vatutin and Mikhaylov (1982). There exists a constant $ c>0$ such that
$$ \max_{0\le k\le N-n}\mid \sqrt{var\mu_0(N)}\; \textbf{P}\;  \lbrace \mu_0(N)=k\rbrace -\frac{1}{\sqrt{2\pi}}e^{-u_k^2/2}\mid\le \frac{c}{\sqrt{var\mu_0(N)}} ,$$
where $ n=\max(n_1,\dots,n_s), \;$ $ u_k=(k-E\mu_0(N))/\sqrt{var\mu_0(N)}. $\\
\end{theorem}

In [40] an  approach known as the Frobenius - Harper technique have been developed. They show the distributional equality of $\mu_0(N)$  to a sum of $ N-n$ independent Bernoulli r .v.s. whose probability of success is $ (1+d_m)^{-1},\; d_m\ge 0$, where   $\;-d_1,\dots, -d_{N-n}\;$, are only  roots of  the  the probability generating function of $\mu_0(N)$, which, in turn, is a polynomial of degree $ N-n . $ Obviously, this approach can be used to derive the asymptotic expansion of distribution of $\mu_0(N)$, but the result will be impractical, since the expansion terms will depend on the roots $d_m$s, finding which is a separate complex problem. The combination of this technique and the approach described in Section 2 above, with specification due to the r.v. $\mu_0(N)$, in Mirakhmedov and Mirakhmedov (2010) the asymptotic expansion with an explicit formula for the  terms was obtained. A particular result is as follows.\\
     
 Set $ Q_s=q_1\cdot\dots\cdot q_s,\;$ $ g\sim g_m:= \textbf{I}\lbrace \xi_{1,m}+\dots +\xi_{s,m}=0\rbrace - Q_s +Q_s\sum_{k=1}^s q_k^{-1}(\xi_{k,m}-p_k), $
Note that $g_m,\;m=1,\dots,N$ are iid r.v.s, and $E\mu_0(N)=NQ_s $ and $var\mu_0(N)=N\sigma^2(1+o(1))$, where  

$$\sigma^2:=var g=Q_s\Big(1-Q_s\Big(1+\sum_{k=1}^s\frac{p_k}{q_k}\Big)\Big),$$
Put for short $\; \tilde{g}=g/\sigma,\quad \tilde{\xi}_{k,m}=(\xi_{k,m}-p_k)/\sqrt{p_kq_k}, $
$$ A:= \sum_{k=1}^s\Big( Eg_1^2\tilde{\xi}_{k,1}E\tilde{\xi}_{k,1}^3 -Eg_1^2\tilde{\xi}_{k,1}^2 +1\Big),\quad  B:= E\tilde{g}^4-3\sum_{k=1}^s\Big((E\tilde{g}^2\tilde{\xi}_{k,1})^2+1\Big).$$
Let $H_l(u)$ stands for the $l$ th - order Hermite-Chebyshev polynomial:
$ H_2(u)=u^2-1,\; H_3(u)=u^3-3u,\;H_4(u)=u^4-6u^2+3,\; H_6(u)=u^6-15u^4+45u^2-15 ,  $

$$\bar{\mathbb{W}}(u)= \frac{1}{\sqrt{2\pi}}e^{-u^2/2}\Big[ 1+\frac{Eg_1^3}{6\sqrt{N}\sigma^3}H_3(u) + \frac{1}{4N}\Big( \frac{(Eg_1^3)^2}{18\sigma^6}H_6(u) +\frac{1}{8}H_4(u)B+H_2(u)A\Big)\Big] ,    $$

$$ b_N:= \frac{Q_s}{(N-1)\Big(1-Q_s\Big(1+\sum_{k=1}^sp_kq_k^{-1}\Big)\Big)} \sum_{j=2}^s \frac{(-1)^{j-2}}{(N-1)^{j-2}} \sum\;\prod_{l=1}^j\frac{p_{i_l}}{q_{i_l}}\;, $$
where $ \sum $ stands for the summation over all $j $-tuples $(i_1,\dots,i_j)$ where $i_k\neq i_l, \;$ if $k\neq l $  and all $i_l\in (1,2,\dots, s).$ We emphasize that $var\mu_0(N)=N\sigma^2(1+b_N)$.
  
\begin{theorem}\label{th 5.15} Mirakhmedov and Mirakhmedov (2010). Set $u_k=(k-E\mu_0(N))/\sqrt{var\mu_0(N)}$. One has

$$\quad \max_{0\le k\le N-n}\mid \sqrt{var\mu_0(N)}P \lbrace \mu_0(N)=k\rbrace - \bar{\mathbb{W}}(u_k) -  \frac{1}{\sqrt{2\pi}}e^{-u_k^2/2}H_2(u_k)b_N \mid \le\frac{C_1}{(var\mu_0(N))^{3/2}}\;.$$
\end{theorem}

Theorem 5.15 considers the standardized version of $\mu_0(N)$ by its exact mean and variance. If instead we use  the  asymptotic value of the variance then the term containing $ H_2(u_k)b_N  $ will be removed.  

\begin{theorem}\label{th 5.16}Set $x_k=(k-NQ_s)/\sqrt{N}\sigma. $ One has
$$\quad  \max_{0\le k\le N-n}\mid \sqrt{N}\sigma P \lbrace \mu_0(N)=k\rbrace - \bar{\mathbb{W}}(x_k)  \mid \le\frac{C_2}{(N\sigma^2(1+b_N))^{3/2}}\;.$$
\end{theorem}

Theorem 5.16 don't presented in Mirakhmedov and Mirakhmedov (2010), nevertheless it can be easily derived from Theorem 5.15 and the fact that $ var\mu_0(N)=N\sigma^2(1+b_N) .$ Note that if all $p_k, \; k=1,\dots, s $ are bounded away from zero and one, then the rhs of the inequalities of Theorems 23 and 24 are of order $O(N^{-3/2}) .$\\  
  
\textbf{Example 5.2.2.\;} Consider the sample scheme without replacement from a stratified population of size $ \Omega_N ;\; $  the strata are indexed by $ m=1,\dots, N ;\;$  $ \omega_m $ is the size of the $m$th stratum, with $ \Omega_N=\omega_1+\dots+\omega_N ;\;  $ now $\eta_m $ is the number of elements of the \textit{m}th stratum appearing in a sample of size $ n $. This scheme is a version of GUM(1,N) (1) where $\xi_m \sim Bi(\omega_m, p),\; $ $ p=n/\Omega_N ,\; q=1-p , $  $m=1,\dots, N$ and $\xi_m =0$  for $m\ge N+1.$  We consider the case where the strata sizes $\omega_m$ may increase together with $ N, $ but such that $\bar{\omega}_N:=\max_{1\le m\le N}\omega_m=o((nq)^{1/4}).$ \\
Let's consider a "sample-sum" $ S_{n,N} $,that is the sum of elements of the population that appeared in the sample.\\
  \indent Define the indicator r.v.s $ \eta_{m,j} $ which are equal to one if an element
$ X_{m,j}$ (say) of the \textit{m}th stratum appears in the sample, or else it equals to zero, so that $ \eta_{m,1}+\dots + \eta_{m,\omega_m}=\eta_m .$ Then $ S_{n,N}^{(m)}=\sum_{k=1}^{\omega_m}X_{m,k}\eta_{m,k} $ represents the sum of elements
of the \textit{m}th stratum which appear in the sample, and $ S_{n,N}=\sum_{m=1}^N S_{n,N}^{(m)} $ is a particular variant of the r.v. (3).
    
Assume that the elements of \textit{m}th stratum are independent r.v.s $X_{m1},\dots, X_{m,\omega_m} $ which have a common distribution, the same as that of an r.v. $Y_m ,$ $m=1,\dots, N.$ We also assume that $ Y_1,\dots, Y_N$ are independent r.v.s. Then 
$$\mathcal{L}\Big( S_{n,N}\Big)= \mathcal{L} \Bigg(\sum_{m=1}^N \sum_{j=1}^{\eta_m}X_{m,j} \textbf{I} \lbrace\eta_m\ge 1\rbrace \Bigg).$$

The r.v. $ S_{n,N}$ in this general form have been considered in Mirakhmedov et al 2014)[31] for the first time. We will present below just their a particular slightly simplified results. Set 

$$\gamma_N=\frac{1}{\Omega_N}\sum_{k=1}^N \omega_mEY_m ,\quad \sigma_N^2=p\sum_{k=1}^N \omega_m\Big(E
(Y_m-\gamma_N)^2 - p(E(Y_m-\gamma_N) )^2\Big) ,  $$
$$ g_m(\xi_m)=\sum_{k=1}^{\xi_m} \textbf{I}\lbrace\xi_m\ge 1\rbrace \Big(X_{m,k}-\gamma_N)- \omega_m p E(Y_m-\gamma_N)\big) , $$

$$\beta_{k,N}=\frac{1}{\sigma_N^k}\sum_{m=1}^N E\mid g_m(\xi_m)\mid^k\le \frac{2^kp(1+p^{k-1})}{\sigma_N^k} \sum_{m=1}^N\omega_m^kE\mid Y_m-\gamma_N\mid^k .  $$

\begin{theorem}\label{th5.17}[31].  For arbitrary $\delta\in (0,1]$ one has

$$\sup_{-\infty < u <\infty} \mid \textbf{P}\lbrace S_{n,N} <u\sigma_N+n\gamma_N\rbrace -\Phi (u)\mid\le C\Big(\beta_{2+\delta,N}+\frac{1}{(nq)^{\delta/2}}\Big) .$$

\end{theorem}

\begin{theorem}\label{th5.18} [31]). There exist a constants $c>0$ snd $C>0$ such that

$$\sup_{-\infty < u <\infty} \mid \textbf{P}\lbrace S_{n,N} <u\sigma_N+n\gamma_N\rbrace -\Phi (u)
-\frac{(1-u^2)e^{-u^2/2}}{6\sqrt{2\pi}}\alpha_{3,N}\mid $$
$$\le C\Big(\beta_{4,N}+\frac{\bar{\omega}_N}{nq} + \chi_N\Big),$$
where  
 $\mathbf T_N=\min(\beta_{3,N}^{-1},\sqrt{nq}/\bar{\omega}_N^2) ,$

$$\alpha_{3,N}=\frac{p^{3/2}}{ \sigma_N^3}\sum_{m=1}^N\omega_m\Big(E(Y_m-\gamma_N)^3 -3p    E(Y_m-\gamma_N)E(Y_m-\gamma_N)^2 -2p^2 (E(Y_m-\gamma_N)^3  \Big), $$

for the $ \chi_N $ the following bounds hold:
$$ \chi_N\le C_1 \sqrt{n} \ln(\sigma_N \beta_{4,N})^{-1}\exp\left\lbrace -n\Bigg(1-\sup_{c\sigma_N^{-1}\tilde{T}_N\le \mid t\mid\le (\sigma_N \beta_{4,N})^{-1}}   \frac{1}{N}\sum_{m=1}^N \mid Ee^{itY_m}\mid\right\rbrace ,$$

and 

$$ \chi_N\le C_1 \sqrt{n} \ln(\sigma_N \beta_{4,N})^{-1}\exp\left\lbrace -n\Bigg(1-\sup_{c\sigma_N^{-1}\tilde{T}_N\le \mid t\mid\le (\sigma_N \beta_{4,N})^{-1}}   \frac{1}{\Omega_N}\mid\sum_{m=1}^N \omega_m Ee^{itY_m}\mid\Big)\right\rbrace. $$

\end{theorem}
\textbf{Remark 5.3.} The first bound is sharper than the second under certain conditions, for example, when $q$ is close to zero but a sufficiently large number of population elements are non-degenerate r.v.s. In contrast, the second bound is applicable when all population elements are degenerate, i.e., nonrandom, namely, $Y_{mN}=a_{mN},\quad m=1,\ldots,N,$ provided that $q$ is not too close to zero. In the latter case, the second bound is exponentially small. For instance, suppose that, for a given $d_0>0$, there exist constants $d_1>0$, $\varepsilon>0$, and $\delta>0$, independent of $n$, such that
$\sharp\lbrace k:\mid t(a_{kN}-\bar a)/b\sqrt n-x-2\pi\nu\mid >\varepsilon\rbrace\ge \delta n,$ 
where $b=\max_m a_{mN},\quad\bar a=N^{-1}(a_{1N}+\cdots+a_{NN}),$ for every fixed $x$, all $n\ge1$, $t\in(d_0,d_1)$, and every integer $\nu$; see Bickel and van Zwet (1978), Robinson (1978), and Mirakhmedov (1983).\\

Let $ \omega _1=\dots=\omega_N=1 $, that is $ \Omega_N=(Y_{1,N},\dots,Y_{N,N})$ is  a triangular array  of independent r.v.s . Then $\eta_m=1$ if $Y_{m,N}$ appears in the sample, or  else it equal to zero,$m=1,\dots,N :$ also now

\begin{equation}\label{57}
 S_{n,N}=\sum_{m=1}^N\eta_mY_{m,N} .
 \end{equation}

This r.v. is of interest in several contexts, including the estimation of population totals in multistage sampling schemes, for example, in stratified and cluster sampling designs; see Mirakhmedov et al. (2015) [32]. In the particular case where $Y_{m,N}=a_{m,N},\qquad m=1,\dots,N,$ are nonrandom real numbers, $S_{n,N}$ represents a linear two-sample rank statistic, where the values $a_{m,N}$ serve as the scores of the statistic; see Bickel and van Zwet (1978). Consequently, there is an extensive literature devoted to the r.v. (57); see Ibragim and Mirakhmedov (2016) [19] and the references therein. In that work, a Lindeberg-type condition for asymptotic normality, a Berry--Esseen-type bound, an Edgeworth asymptotic expansion with refined remainder terms, and Cramér-type large-deviation results are established. The proofs are based mainly on the approach presented in Section 2, and some of the results may be viewed as consequences of Theorems 3.1-3.4 in Section 3.\\

\textbf{Agreement on the future notations.} In what follows we consider r.v.s of the form $f(\eta_1)+\dots+f(\eta_N) $ ( known as symmetric statistics) defined in the versions of GUM (1,N) generated by a non-negative integer iid r.v.s $\xi_1,\dots,\xi_N  .$ The r.v. $\xi_m$ will correspond to the GUM versions under consideration. To avoid repetition of notations, we use the following general notation: $\langle x\rangle $ is the distance between real $x$ and integers,  $\zeta'$ is independent copy of the r.v. $\zeta$, $\zeta^*=\zeta-\zeta' ,$\\
 $ \xi\equiv\xi_1,\;$ $ v^2=var\xi,\; \gamma=cov(f(\xi),\xi)/v^2,\;\; g(\xi)= f(\xi)-Ef(\xi)- \gamma(\xi-E\xi),\,$ 
  $\sigma^2=var\,g(\xi)=varf(\xi)-v^2\gamma^2,$ 
  $$ \alpha_{k,j}=\frac{1}{\sigma^k v^j}Eg^k(\xi)(\xi-E\xi)^j,\quad M= 1-\sup_{(8v\sqrt{\alpha_{0,4}})^{-1} \le\mid\tau\mid\le\pi}\mid E\exp\lbrace i\tau\xi\rbrace \mid^2 ,\quad \mathbf{T}^2=\min\Big(\frac{1}{64\sigma^2\alpha_{4,0}}, \frac{M}{36\sigma^2}\Big),$$ 
 \begin{equation}\label{58}
   H(t,\tau)=E\langle tg^*(\xi)+\tau\xi^*\rangle^2,\; H_N= N\inf_{\mathbf{T}\le\mid t\mid\le\pi,\,\mid\tau\mid\le\pi}H\Big(\frac{t}{2\pi},\frac{\tau}{2\pi}\Big).
\end{equation}
\medskip

\textbf{Example 5.2.3.\;}  Assume from an urn containing balls of \textit{N} different colors, with $m$  balls of each color, a sample of size $n$ ($0 < n < Nm $) is drawn according to a simple sampling without replacement scheme. That is, we are dealing with scheme from Example 5.2.2, where $\omega_1=\dots=\omega_N= m\ge 2 .$ Then $ \eta_k $ denotes the number of balls of the \textit{k}-th color in the sample, and hence  the independent r.v.s $ \xi_k\sim Bi(m,p),\; k=1,\dots,N $ where $p=n/mN .$  We are interested in the r.v. $\mathbb{H}_N$ which  denote the number of unordered pairs of balls of the same color in the sample:    

$$ \mathbb{H}_N=\frac{1}{2}\sum_{k=1}^N \eta_k(\eta_k -1) .$$    

Set
$q=1-p,\;\;a=m(m-1)p^2/2,$ we have $v^2=mpq, \quad \sigma^2 = m(m-1)p^2q^2/2. $ From the general theorems we obtain the following statements.
 
 \begin{theorem}\label{th 5.19} There exist positive constants $C$ and $c$ such that
 $$\sup_{-\infty <u<\infty} \mid P\lbrace \mathbb{H}_N < u \sigma \sqrt{N}+N a\rbrace - \Phi(u)\mid\le C\Big( \frac{1}{\sqrt{nmpq^2}}+\frac{1}{\sqrt{nq}} +\frac{1}{\sqrt{N}} +\sqrt{nq}\,e^{-cN} \Big) .$$
\end{theorem}
As an obvious consequence of Theorem 5.19 we obtain the following claim:\\
 \begin{theorem}\label{th28}
 If $n,\, N\,\to\infty $ such that $nmpq^2=(nq)^2/N\to\infty$ and $\ln(nq)=o(N)$ (that is $\sqrt{N}\ll nq\ll e^{cN}$ ) then the r.v. $\mathbb{H}_N$ has asymptotically normal distribution whose mean is $Na$ and variance is $N\sigma^2 .$
\end{theorem}
Theorems  5.20 is already considerably refinements of Theorem 1 by Pavlov and  Cherepnova (2002), giving weaker and unified condition of asymptotically normality o the r.v. $\mathbb{H}_N.$\\
\begin{theorem}\label{th5.21} Let $ mpq \le c_0 < \infty $. Then uniformly in $k\in \lbrace0,1,\dots ,n(n-1)/2 \rbrace $  
 $$ \mid \sigma\sqrt{N}\mathbf{P}\lbrace \mathbb{H}_N=k \rbrace -\frac{1}{\sqrt{2\pi}}e^{-u_k^2/2}\Big( 1+\frac{3u_k-u_k^3}{6 \sqrt{N}}\alpha_{3,0} + \frac{u_k}{3\sqrt{N}}\alpha_{1,2} \Big)\mid 
 \le C \Big(\frac{1}{nmpq^2}+\frac{1}{N} +nq\,e^{-cN} \Big),$$  
where $u_k=(k-Na)/\sqrt{N}\,\sigma ,$ 
$$  \alpha_{3,0}=\frac{\sqrt{2}(1+2(m-4)pq)}{pq\sqrt{m(m-1)}} , \quad  \alpha_{1,2}=\sqrt{\frac{2(m-1)}{m} }.$$ 
\end{theorem} 
          
\noindent \textbf{Proof of Theorems 5.19 and 5.21.}  On using moments of the binomial distribution we found the above formula of $\alpha_{3,0}$ and  $\alpha_{1,2}$, also we get
 \begin{equation}\label{59}
 \alpha_{0,4}=\frac{1+3(m-2)pq}{ mpq}, \quad N^{-1} \alpha_{4,0}\le C  \Big(\frac{1}{nmpq^2}+\frac{1}{nq} +\frac{1}{N} \Big) .
\end{equation}
We have $ \mid Ee^{i\tau\xi}\mid^2\le \exp\lbrace-4mpq\sin^2\tau/2\rbrace \le \exp\lbrace-4mpq\tau^2/\pi^2\rbrace, \;\; \mid\tau\mid\le \pi.$ By this and (59), considering separately the cases when $mpq\ge 1 $ and $mpq\le 1 $
using inequality (41) we get 

\begin{equation}\label{60}
 M\ge c\min(1,mpq) , \quad \text{and}\quad \mathbf{M}(n)= c\min(N, nq) .
 \end{equation}
Theorem 5.19 follows from Theorem 3.2 since (59) and (60).
 Further, if $mpq< c_0$ then $nmpq^2\le c_0nq,\; (\sigma^2\alpha_{4,0})^{-1}\ge  c_1 ,$ and $ (N\sigma^2)^{-1}\mathcal{M}(n)\ge  c_2 ,$ some $c_1>0,\; c_2>0$ depend on $c_0$. Therefore, $(N \mathbf{T}(n)\ge c_3>0.$ By this, the statement of Remark 5.4 below, (17) and that  $\xi\sim Bi(m,p)$ we obtain  $ \bar{\chi}(n,c_3,\pi)\le Cnq\, e^{-c N}. $ Thus Theorem 5.21 follows from Theorem 5.3, (59),(60) and the last inequity.\\

 \textbf{Remark 5.4.}  Let $\xi $ be a non-negative integer r.v., $ f(\xi) =\xi(\xi-1)/2 .$  It can be seen that

$$H(t,\tau)= \sum_{k,l}\langle \vartheta_{k,l} \rangle^2 P\lbrace\xi=k\rbrace P\lbrace\xi'=l\rbrace ,$$
where $ \vartheta_{k,l}=(k-l)(k+l-1)t/2 +(k-l)\tau .$ As in Lemma 2 of Pavlov and  Cherepanova (2002),one can prove the following statement:\\
 \indent For all real $t$ and $\tau$ such that $ 0<\varepsilon\le\mid t\mid \le 1/2, \mid\tau\mid\le 1/2,  $ and any non-negative integer $k$ and $l$ one has $ \max (\;\langle\vartheta_{k,l}\rangle,\langle\vartheta_{k+1,l}\rangle, \langle\vartheta_{k+2,l}\rangle \;) \ge c >0,$ here $c$ depend on $\varepsilon.$\\

\subsection{ Random coverage of a circle }

 \bigskip
 \indent
 
Consider $n\ge 2$ boxes arranged clockwise on a circle and numbered
$0,1,\dots,n-1$. Select $N-1$ boxes by simple random sampling
without replacement, where $2\le N\le n$. Let $1\le R_1<\cdots<R_{N-1}\le n-1$ be the numbers of the selected boxes, and set $R_0=0$ and $R_N=n$.
We regard \(R_0,\dots,R_{N-1}\) as \(N\) distinct starting points whose
locations among the \(n\) boxes are chosen at random. Starting from each
of these points, we designate the next $r$, $1\le r\le n$, boxes, including the starting box itself, as covered. Any boxes that remain uncovered are designated as vacant. Let us define the spacings 
$\eta_k=R_k-R_{k-1},\; k=1,\dots,N. $
Then for any integers $k_j\ge 1,$ such that $ k_1+\dots+k_N=n$
\begin{equation}\label{61 }
\mathbf P\{\eta_1=k_1,\dots,\eta_N=k_N\}=\frac{1}{ \binom{n-1}{N-1}}, 
\end{equation}

Therefore, the equality (1) hold, where the iid  r.v.s $\xi_1,\dots,\xi_N$ have a common geometric distribution:
$$\mathbf P\{\xi_m=k\}=(1-p)p^{k-1},\quad k=1,2,\dots,$$
with arbitrary $p\in(0,1)$. It is convenient to choose $p=(n-N)/n,\quad q=1-p=N/n.$

The random circle-coverage problem described above can also be interpreted
as a random allocation model. Namely, it corresponds to allocating $n$
indistinguishable particles among $N$ distinct cells, with each cell
containing at least one particle, in such a way that all $ \binom{n-1}{N-1} $
possible allocations are equally likely. Equivalently, after placing one
particle in each cell, the remaining \(n-N\) particles are allocated according to the Bose--Einstein urn model.

Another interpretation is obtained by considering all ordered partitions
of the integer $n$ into $N$ strictly positive summands, that is, all
solutions of $\eta_1+\cdots+\eta_N=n,\quad\eta_j\ge1,\quad j=1,\dots,N.$
Equipping this set of solutions with the uniform probability distribution
yields a random ordered $N$-partition of the integer $n$.

There is also a Pólya-urn interpretation. Suppose that an urn initially
contains $N$ balls of different colors, one ball of each color. At each
step, a ball is drawn at random and then returned to the urn together with
one additional ball of the same color. This procedure is repeated $n-N$
times. Then $\eta_j$ denotes the number of balls of the $j$th color in
the urn at the end of the procedure. Then for the r.vec. $(\eta_1,\dots,\eta_N)$ equality (61) hold. Thus, this is in fact a multicolor Pólya urn model.\\

 \indent   
 Throughout this sub-section we assume $np\to\infty.$ 
 We have $ E\xi =q^{-1},\; v^2= pq^{-2},$   $\; v^2\alpha_{0,4}= 1+9pq^{-2},\, $
  $\varrho_4(n)=q(np)^{-1}+9N^{-1} .$  Also $\mid Ee^{i\tau\xi_m}\mid^2 = (1+4pq^{-2}\sin^2 \tau /2)^{-1}\le \exp\lbrace - 4p(q^2+4p)^{-1}\sin^2 \tau /2)\rbrace, $ since inequality $(1+x)^{-1}\le e^{-x/(1+x)},\; x\ge 0.$ Therefore, using the inequalities (41) we conclude that 
 \begin{equation}\label{61} 
  M \ge c_0\frac{pq^2}{(q^2+4p)^2} , \;\text{and hence}  \quad\mathbf{M}(n)\ge c_0\frac{Npq^2}{(q^2+4p)^2} .
\end{equation}
\medskip
\textbf{Example 5.3.1.} Let $\mathbb{W}_N(r)=\mathbb{W}(r,N,n)$ be the number of distinct blocks of vacant boxes.  Then  
 $$ \mathbb{W}_N (r)=\sum_{k=1}^N \textbf{I}\lbrace \eta_k > r \rbrace  .$$

We have   $\,\varsigma = \textbf{I}\lbrace\xi > r\rbrace ,\;$ $E\varsigma =p^r\,$ and $\, \sigma^2 = p^r(1-p^r)(1- \theta), $ with $\, \theta=r^2q^2p^{r-1}/(1-p^r). $
 \begin{theorem}\label{th 5.22 }  There exist $C>0$ such that 
  $$\sup_{-\infty <u<\infty} \mid \mathbf{P}\lbrace \mathbb{W}_N (r)< u \sigma \sqrt{N}+N p^r\rbrace - \Phi(u)\mid\le C Q_N, $$
  where $ Q_N =Q_{1,N}+Q_{2,N}+Q_{3,N},$
$$Q_{1,N}= \frac{1}{1-\theta}\frac{1}{ \sqrt{Np^r(1-p^r)}},\,\; Q_{2,N}=\frac{1}{\sqrt{Npq^2}},\,\; Q_{3,N}= 
  \sqrt{Npq^{-2}}\,\exp\left\lbrace -\frac{c_0Npq^2}{(q^2+4p)^2}\right\rbrace.$$
\end{theorem}
 
{\bf Remark 5.5.}   If $q\to 0\;$, $rq\le c_1 <1\;$ and $ q\gg\sqrt{N^{-1}\ln N}\; $ then $Q_N\asymp (q\sqrt{N})^{-1}.$  Let $q\to 1,\;$ then $Q_N\asymp (\sqrt{N}p^{3/2})^{-1}\;$ if  $r=1$, and  $ Q_N\asymp (\sqrt{N}p^{r/2})^{-1} \;$ if $r\ge 2.\;$ Finally, if $q$ is bounded away from 0 and 1, then $Q_N\asymp (\sigma\sqrt{N})^{-1}\sim (p^{r/2}\sqrt{N})^{-1}\;$ if $r\to \infty,$ and else $Q_N\sim N^{-1/2}\;$ if $r$ is fixed.\\
  Thus  as a consequence of Theorem 5.22 we obtain the following
  
 \begin{theorem}\label{th 5.23 }  Let
 $\; (i) \;q\to 0,\; rq\le c_0<1,\;  q\gg\sqrt{N^{-1}\ln N}\;\; and \;\; Nq^2\to\infty,\;$  or
 $ \quad(ii)\; q\to 1,\; and\;  Np^r\to\infty, \; r\ge 2,\;\; or  \quad \;  r=1,\; Np^3\to\infty,\; $ or \; (iii)\;$ q $ is bounded away from 0 and 1 , $\;Np^r\to\infty.$
Then the r.v. $W_N(r)$ has asymptotically normal distribution whose mean and variance are equal to $Np^r$ and $\;N\sigma^2 $, respectively.
 \end{theorem}
   
  \begin{theorem}\label{th 5.24 } Let $ \sqrt{N^{-1}\ln N}\ll q\le c<1 \; and \; rq\le c_0 <1.$ Then  uniformly in $k\in \lbrace0,1,\dots ,N \rbrace$  
 \begin{equation}\label{62} 
  \mid \sigma\sqrt{N}\mathbf{P}\lbrace \mathbb{W}_N(r)=k \rbrace -\frac{1}{\sqrt{2\pi}}e^{-u_k^2/2}\Big( 1+\frac{3u_k-u_k^3}{6 \sqrt{N}}\alpha_{3,0} + \frac{u_k}{3\sqrt{N}}\alpha_{1,2}\Big)\mid \le C \Big( Q_{1,N}^2+ Q_{2,N}^2 \Big), 
\end{equation}
where $u_k=(k-Np^r)/\sqrt{N}\,\sigma ,$ 

$$\alpha_{3,0}=\frac{(1-p^r)(1-2p^r)-3r^2q^2p^{r-1}(1-2p^r)+r^3q^3p^{2r-2}\bigl(3rq-(1+p)\bigr)}{p^{r/2}\left[1-p^r-r^2q^2p^{r-1}\right]^{3/2}},$$
$$\alpha_{1,2}=\frac{r( q^3r-1-p)}{pq(1-p^r-r^2q^2p^{r-1})}.$$
\end{theorem}
 
\textbf{Proof of Theorems 5.22 and 5.24.} 
  One counts 
 $$ Eg^4(\xi)= p^r(1-p^r)(1-3p^r+3p^{2r})+ \frac{r^2q^2}{p}\left(-4p^{2r}+18p^{3r}-18p^{4r}\right)$$
$$+\frac{r^4q^4}{p^3} \left[6p^{3r+1}(1-4p^r)+p^{4r}(1+7p+p^2)\right] 
- 4p^{4r-3}r^6q^6.$$
 Whence (recall $\theta=r^2q^2p^{r-1}/(1-p^r)$ )
 $\;Eg^4(\xi)\le 15\,p^r(1-p^r)(1+\theta+\theta^2)\le 45\,p^r(1-p^r).$
Therefore we obtain
$$\alpha_{4,0}\le \frac{45}{(1-\theta)\sigma^2} . $$
 By this and (61) Theorem 3.2 implies Theorem 5.22.  Further, $\;1-p^r=q(1+p+\cdots+p^{r-1})\ge qrp^{r-1}$ and hence $\theta =rq\, [rp^{r-1}/(1+p+\cdots+p^{r-1})]\le rq.$
  Next, under the assumption of Theorem 5.24 we have $ \mathbf{T}\ge c >0 .$ Therefore,  reasons similar to those used in (54) allow us to get $\bar{\chi}(n,\mathbf{T},\pi\sqrt{N}\sigma)\le c_3 N \sigma \sqrt{p}q^{-1}\exp\lbrace-c_4Nq^2p^{2r}\rbrace ,$ which is negligible relative to RHS terms of (62). Theorem 5.24 follows from Theorem 3.4.\\

\medskip 

\textbf{Example 5.3.2.} Let us consider $ \mathbb{V}_N(r) $, the total number of vacant boxes: 
$$\mathbb{V}_N(r)= \sum_{k=1}^N (\eta_k -r)\textbf{I}\lbrace \eta_k > r\rbrace.$$  
 Let now $\varsigma= (\xi - r)\textbf{I}\lbrace \xi > r\rbrace ,\; r\ge 2.$ We have $E\varsigma =p^rq^{-1},\;$ $  \sigma^2= q^{-2}p^r[1+p-p^r-p^{r-1}( rq+p)^2 ] . $ 

 \begin{theorem}\label{th 5.25 }  There exist  $C>0$ such that 
  $$\sup_{-\infty <u<\infty} \mid \mathbf{P}\lbrace \mathbb{V}_N (r)< u \sigma \sqrt{N}+N p^rq^{-1}\rbrace - \Phi(u)\mid\le C \Big(\bar{Q}_{1,N}+Q_{2,N}+   Q_{3,N}\Big), $$
  where 
$$\bar{Q}_{1,N}= \frac{1}{\sqrt{Np^r}\left[1+p-p^r-p^{r-1}(p+rq)^2\right]}.$$
 \end{theorem}
   
\begin{theorem}\label{th 5.26 }  Uniformly in $k\in \lbrace0,1,\dots ,N \rbrace$  
 $$ \mid \sigma\sqrt{N}\mathbf{P}\lbrace \mathbb{V}_N(r)=k \rbrace -\frac{1}{\sqrt{2\pi}}e^{-u_k^2/2}\Big( 1+\frac{3u_k-u_k^3}{6 \sqrt{N}}\alpha_{3,0} + \frac{u_k}{3\sqrt{N}}\alpha_{1,2}\Big)\mid \le C ( \mathcal{Q}_N + D_N ),  $$

where $u_k=(k-Np^r)/\sqrt{N}\,\sigma ,$ $ \mathcal{Q}_N=\bar{Q}_{1,N}^2+ Q^2_{2,N}+Q_{3,N}, $
$$ D_N= N\sigma q^{-1}\sqrt{p}\exp\left\lbrace -cq^2p^{2r+1}
 \min\left( \frac{q^4\sigma^2}{p^r},\frac{pq^2}{\sigma^2}\right)\right\rbrace          ,$$
$$\alpha_{3,0}=\frac{p^r \mathcal{A}(p,q,r)}{q^3\sigma^3},\quad 
\alpha_{1,2}=\frac{r(r-1)qp^{r-1}}{\sigma},$$
 $$\mathcal{A}(p,q,r)= (1+p)(1+p-3p^r)+2p^{2r}
- 3p^{r-1}(p+rq)\left[(1+p)(1+p+qr)-2p^r(p+qr)-q\right]$$
$$+ p^{2(r-1)}(p+rq)^2\left[3p(1+p)+3qr(2p+qr) -(p+rq)(1+p)\right].$$
\end{theorem}  
  
  A detailed description of the asymptotic behavior of the terms $\mathcal{Q}_N$ and $D_N$ in the cases $q\to 0$ and $q\ge c>0$ is given below at the end of the proof of the theorems. 
  
 \textbf{Proof of Theorems 5.25 and 5.26.} We have $\gamma= p^{r-1}(p+rq),\,$ $g(\xi)=\varsigma -p^rq^{-1} -\gamma (\xi-q^{-1}),\,$ $Eg^4(\xi)\le 8(E( \varsigma-p^rq^{-1})^4+\gamma^4 E(\xi-q^{-1})^4).$ Since $E\xi^4=q^{-4}(1+11p+11p^2+p^3)\le 24q^{-4} ,$ on using the memoryless property of the geometric distribution we obtain $E\varsigma^4=p^rE\xi^4\le 24p^rq^{-4} ,$ whence $E(\varsigma-p^rq^{-1} )^4\le 64 p^rq^{-4}.$ Next, as $E(\xi-q^{-1})^4=q^{-4}p(1+7p+p^2)\le 9pq^{-4}$ we have $\gamma^4 E(\xi-q^{-1})^4\le 9q^{-4}p^{4r-3}(p+rq)^4.\;$ On the other hand  $p+rq=1+(r-1)q,\; p^{r-1}=(1-q)^{r-1}\le \exp\{-(r-1)q\},$  and hence $p^{3(r-1)}(p+rq)^4\le e^{-1}(4/3)^4.$ Therefore $\gamma^4 E(\xi-q^{-1})^4  \le c\,q^{-4}p^r .$  Thus, we have $ Eg^4(\xi)\le cq^{-4} p^r $, that is $\beta_4(n)\le c p^r/Nq^4\sigma^4 .$ By this and (61) Theorem 33 follows from Theorem 2.\\
\indent Let us to get a lower bound for the $ H_N. $  Set for brevity $u=t/2\pi,\; v=\tau/2\pi ,$ and consider  $H(u,v)=E\langle ug^*+v\xi^*\rangle^2,\;$ where  $u g^*+v\xi^*=u(\varsigma-\varsigma')+b(\xi-\xi'),\;\; b= v-u\gamma.$
We now restrict the expectation defining \(H(u,v)\) to several particular values of $(\xi,\xi')$. Consider the two events:
$(\xi,\xi')=(r,r-1),\;\; (\xi,\xi')=(r+1,r)\quad$ and
$\quad (\xi,\xi')=(r-1,r),\;\; (\xi,\xi')=(r,r+1).$
On the event \((\xi,\xi')=(r,r-1)\), we have $u g^*(\xi)+v\xi^*=b,$ whereas on the event \((\xi,\xi')=(r+1,r)\), we have \(u g^*+v\xi^*=u+b.\) On the other hand
$$\mathbf P\{\xi=r,\xi'=r-1\}=q^2p^{2r-3}\quad \text{and}\quad \mathbf P\{\xi=r+1,\xi'=r\}=q^2p^{2r-1}.$$
Taking also the symmetric counterparts into account, and using $\langle -x\rangle=\langle x\rangle,$
we obtain $ H(u,v)\ge 2q^2p^{2r-1}\left[\langle b\rangle^2+p^2\langle u+b\rangle^2\right].$
Consequently,
\begin{equation}\label{63}
\inf_{|v|\le 1/2}H(u,v)\ge 2q^2p^{2r-1}\inf_{b\in\mathbb R}\left[\langle b\rangle^2 + p^2\langle u+b\rangle^2\right]= 2q^2p^{2r-1}\frac{p^2}{1+p^2}\langle u\rangle^2
\end{equation}   

Indeed, since the expression is periodic in $b$ modulo integers, it is enough to minimize $b^2+p^2(u+b)^2$ over the relevant interval. The unconstrained minimum is attained at $b_*=-p^2/(1+p^2)u,$ and the minimum value is
$b_*^2+p^2(u+b_*)^2=u^2p^2/(1+p^2).$ Hence (63) follows.\\
\indent Further, as $|u|\le 1/2,$ we have 
$\langle u\rangle^2=u^2\ge \mathbf{T}^2/4\pi^2.$
Combining this with (63), we obtain 
$$ H_N\ge \frac{q^2p^{2r+1}}{2\pi^2(1+p^2)}\,\mathbf{T}^2 \ge cq^2p^{2r+1}
 \min\left( \frac{q^4\sigma^2}{p^r},\frac{pq^2}{\sigma^2} \right),$$
 since $ Eg^4(\xi)\le cq^{-4}p^r ,\;$ and  $\mathbf{T}^2\ge\min\left(q^4\sigma^2/p^r,pq^2/6\sigma^2\right).$ Thus by (17) $\bar{\chi}(n,\mathbf{T},\pi\sqrt{N}\sigma)\le D_N.$
By this and (61) Theorem 5.26 follows from Theorem 3.4 .

Let us consider the special  regimes of $q$. Set  $\Delta:=1+p-p^r-p^{r-1}(p+rq)^2.$ Then $\sigma^2=p^r\Delta/q^2,\;$ and
 $$D_N = Nq^{-2}p^{(r+1)/2}\sqrt{\Delta}\exp\{-\Lambda_{1,N}\},\; \text{where}\quad \Lambda_{1,N}\asymp N\min\left(q^4p^{2r+1}\Delta\,,\;q^6p^{r+2}\Delta^{-1}\right). $$
Furthermore,
\[\mathcal{Q}_N=\frac1{Np^r\Delta^2}+\frac1{Npq^2}+\frac{Np}{q^2}\exp\left\{-c_1Npq^2\right\}.\]\\

Let $q\to 0.$ Assume $rq\le c_2<1,\; r\ge 2$. Then $p^r\sim 1,\; \Delta \asymp a_rq^3, $ where $d_r=r(r-1)(2r-1).$ Therefore, 
\begin{equation}\label {64}  
\bar{Q}_{1,N}^2\asymp\frac{1}{Nd_rq^6},\;\;\mathcal{Q}_N\asymp \frac{1}{Nd_r^2q^6}+\frac{1}{Nq^2}+\frac{N}{q^2}e^{-cNq^2}\asymp \frac{1}{Nd_r^2q^6}+\frac{1}{Nq^2},
\end{equation}
provided $Nq^2\gg\ln N.$ Also

$$D_N\asymp N\sqrt{d_r}q^{-1/2}\exp\left\lbrace -cN\min\Big(d_rq^7,q^3d_r^{-1}\Big)\right\rbrace. $$
In particular, if $r$ is fixed, 

$$\mathcal{Q}_N\asymp \frac{1}{Nq^6},\quad D_N=o(\mathcal{Q}_N),\; \text{if}\;\; Nq^7\gg\ln N .$$
If $r\to\infty$, then $\Delta \asymp r^3q^3,\; \sigma^2\asymp r^3q ,\,$ and 

$$D_N\asymp Nr^{3/2}q^{-1/2}\exp\left\lbrace -cN\min\Big(r^3q^7,q^3r^{-3}\Big)\right\rbrace. $$
Next, assume $r^3q^2\le C. $ Then 
$$\mathcal{Q}_N\asymp \frac{1}{Nr^6q^6},\quad D_N=o(\bar{Q}^2_{N}),\;\; \text{if}\;\; Nr^3q^7\gg\ln N. $$
Let now $r^3q^2\ge C>0. $ Then   
$$\mathcal{Q}_N\asymp \frac{1}{Nq^2},\quad D_N=o(\bar{Q}^2_{N}),\;\; \text{if}\;\; Nq^3 \gg r^3 \ln N. $$
That is $D_N=o(\mathcal{Q}_N),\;$ if $\; (\ln N)/Nq^7\ll r^3\ll q^{-2}\;$ (hence $Nq^7\gg \ln N\;$) or $\; q^{-2}\ll r^3\ll Nq^3/\ln N, $ (hence $Nq^3\gg \ln N, $) and finally at the transition scale $ r^3\asymp q^2$ we should assume $Nq^5\gg \ln N.$\\
\indent Let $q\to 1,\;r\ge 2.$ Then $ \Delta\sim 1,\; \sigma^2\sim p^r.$ Therefore,
$\mathcal{Q}_N\asymp (Np^r)^{-1},\; D_N\asymp Np^{(r+2)/2}\exp\left\lbrace -cNp^{2r+1}\right\rbrace ,\;\; \text{hence }\;\;D_N=o(\mathcal{Q}_N)\;\;\text{if}\;\; Np^{2r+1}\gg \ln N .$\\
\indent Finally, let $q$ is bounded away from 0 and 1. Then for the fixed $r\ge 2$ we have $D_N\ll \mathcal{Q}_N\asymp N^{-1}.$ Whereas if $r\to\infty $ the results are similar to the case $q\to 1.$ \\
\indent
 The problem of random coverage of a circle has been studied previously. Holst (1985)[11] derived the joint distribution of the $(\eta_1,\dots,\eta_N)$ and explored the connection between this r.vec. and the Bose-Einstein distribution. Ivchenko (1994)[22] established Poisson limit theorems and the central limit theorem for r.v.s $\mathbb W(r)$ and $\mathbb V(r).$ Barlevy and Nagaraja (2015)[2] obtained an explicit formula of the distribution and the first two moments of $\mathbb V(r)$ and described their limiting properties. He also discuss applications to a study of contagion in networks.\\

\subsection{Random forests .}

\bigskip
\indent

 We consider the set $\mathcal{F}(n,N)$ of all distinct forests with $N$ roots and $n$ non-root vertices. The roots are labeled $1,\dots,N$  and the non-root vertices are labeled 
$1,\dots,n.$ We equip this set with the uniform probability distribution, thereby obtaining random forests. The number of non-root vertices of a tree is called its size.  Let $\eta_k$ denotes the size of  a tree of \textit{k}-th  root, then (see Kolchin (1986)) for an arbitrary integers $k_j\ge 0$ such that $k_1+\dots+k_N=n$ :

$$\mathbf P\lbrace \eta_=k_1,\dots,\eta_N=k_N\rbrace =\frac{n!}{N(N+n)^{n-1}} \prod_{m=1}^N\frac{(k_m+1)^{k_m}}{(k_m+1)!} \,. $$
Therefore, the r.vec. $(\eta_1,\dots,\eta_N)$ defines a particular version of GUM(N,1) generated by iid r.v.s  $\,\xi_1,\dots,\xi_N$  having common shifted Borel distribution:

$$ a_k:=\textbf{P}\{\xi = k\}=\frac{(k+1)^k}{(k+1)!}\,p^{\,k}\,e^{-(k+1)p},\quad k=0,1,2,\dots,$$
where arbitrary $p\in (0,1)$. This easily can be seen by using the fact that 
$$\mathbf P\lbrace\xi_1+\dots+\xi_N=n\rbrace=\frac{N}{N+n}\frac{(p(N+n))^n}{n!}E^{-(N+n)p}.  $$ 
   
 It is convenient to choose $p= n/(N+n)=\mu/(1+\mu), $ with $\mu=n/N. $  We assume that $\mu=o(N).$\\
  
 We have 
\begin{equation}\label{67}
  E\xi=p/(1-p)=\mu,\; v^2=var\xi=\mu(1+\mu)^2, \quad\alpha_{04}=(1+11p+3p^2)/p(1-p)\le 15(1+\mu)^2/\mu.
 \end{equation} 
  For every $k\ge 0$ the two events  $\;(\xi,\xi')=(k+1,k),\;\; (\xi,\xi')=(k,k+1) $ give $\xi^*=\pm 1.$ Therefore 
 \begin{equation}\label{68} 
 1-\mid Ee^{i\tau\xi}\mid^2=E(1-\cos\tau\xi^*)\ge 4\sin^2 \frac{\tau}{2}\sum_{k=0}^{\infty}a_ka_{k+1}\ge 4a_0a_1\sin^2 \frac{\tau}{2}.
\end{equation}  
Whence and since $\,v^2\alpha_{0,4}= (1+11p+3p^2)/(1-p)^4$ we get ( see notation (58))
 \begin{equation}\label{69}
  M \ge \frac{e^{-3p}p(1-p)^4}{2\pi^2(1+11p+3p^2)}\ge \frac{1}{30e^3\pi^2}p(1-p)^4=\frac{1}{30e^3\pi^2}\mu(1+\mu)^{-5}. 
  \end{equation}
 
\textbf{Example 5.4.1.} Consider the r.v.  $ \mathbb{S}(r)$ - the number of trees of size $r\ge 1$ of the forest, viz.,

$$ \mathbb{S}_N(r)=\sum_{m=1}^N \textbf{I} \lbrace \eta_m=r  \rbrace. $$ 
   
 Set $\varsigma=\textbf{I}\lbrace\xi=r\rbrace$.We have 
 \begin{equation}\label{70}
 E\varsigma=a_r,\quad \sigma^2=a_r(1-a_r)-\mu(1+\mu)^2\gamma^2  ,\;\quad\gamma=\frac{(r-\mu) a_r}{\mu(1+\mu)^2}.
 \end{equation}
 
 \begin{theorem}\label{th 5.27}
 There exist $C>0$ such that 
  $$\sup_{-\infty <u<\infty} \mid \mathbf P\lbrace \mathbb{S}(r)< u \sigma \sqrt{N}+N a_r\rbrace - \Phi(u)\mid\le C \mathcal{L}_N, $$
  where
$\;\mathcal{L}_N ;= \mathcal{L}_{1,N}+ \mathcal{L}_{2,N} + \mathcal{L}_{3,N},\quad \Delta_r(\mu)= 1- a_r-\mu(1+\mu)^2\gamma^2/a_r,$
$$\mathcal{L}_{1,N}=\frac{1}{\sqrt{Na_r}\,\Delta_r(\mu)},\quad \mathcal{L}_{2,N}=\frac{1+\mu}{\sqrt{n}},\quad \mathcal{L}_{3,N}=\sqrt{n}(1+\mu)\exp\left\lbrace -c_0n(1+\mu)^{-5} \right\rbrace   .$$
\end{theorem}
    
\begin{theorem}\label{th 5.28 }  Uniformly in $k\in \lbrace0,1,\dots ,n \rbrace $  
 $$ \mid \sigma\sqrt{N}\mathbf P\lbrace \mathbb{S}(r)=k \rbrace -\frac{1}{\sqrt{2\pi}}e^{-u_k^2/2}\Bigg( 1+\frac{3u_k-u_k^3}{6 \sqrt{N}}\alpha_{3,0} + \frac{u_k}{3\sqrt{N}}\alpha_{1,2}   \Bigg)\mid \le C \bar{\mathcal{L}}_N, $$
  
where $u_k=(k-Np^r)/\sqrt{N}\,\sigma ,\;$ $\; \bar{\mathcal{L}}_N = \mathcal{L}^2_{1,N} + \mathcal{L}^2_{2,N}  + \mathcal{L}_{3,N}+\mathcal{L}_{4,N} , $  $\; \alpha_{3,0}=a_r(1-a_r)(1-2a_r)/\sigma^3, $

$$ \alpha_{1,2}=\frac{a_r}{\sigma}\Big( \frac{(r-\mu)^2}{\mu(1+\mu)}-1\Big),
\quad
  \mathcal{L}_{4,N}=  N\sigma\sqrt{\mu}(1+\mu)\exp\left\lbrace -CN K_r(\mu,a_r)T(\mu, a_r)\right\rbrace  ,$$
 $$  K_1(\mu,a_1)=\frac{\mu^3}{(1+\mu)^3},\quad K_r(\mu,a_r)=\frac{\mu}{1+\mu}\min\Big(1,\,\frac{a_r^2(1+\mu)^2}{\mu^2}\,\Big),\;\; r\ge 2 ,$$

$$T(\mu,a_r)=\min\left(\frac{\sigma^2}{a_r},\frac{\mu}{\sigma^2(1+\mu)^5}
\right).$$
\end{theorem}  

\medskip

{\bf Remark 5.6 .} (i)\; Let  $r\ge 1$ is fixed. Assume $\mu\le c<\infty .$ Then 
\begin{equation}\label{70 }
\mathcal L_{1,N}^2\asymp\frac{1}{N\mu^3},\;\;\mathcal L_{4,N}\asymp N\mu^{3/2}e^{-cN\mu^4} \;\; \text{if}\; r=1;\; \;\mathcal L_{1,N}^2\asymp\frac{1}{N\mu^r},\;\;\mathcal L_{4,N}\asymp N\mu^{(r+1)/2}e^{-cN\mu^{2r-1}}  \;\text{if}\; r\ge 2, 
\end{equation}
else 
\begin{equation}\label{71}
 \mathcal L_{2,N}^2\sim\frac{1}{n},\;\; \mathcal L_{3,N}\asymp\sqrt{n}e^{-cn},\;\text{for all}\; r\ge 1.
 \end{equation}
Thus $\bar{\mathcal L}_N\asymp(N\mu^3)^{-1}$ \; if \; $N\mu^4\gg\log N ,\; r=1, $ and $\bar{\mathcal L}_N\asymp(N\mu^r)^{-1}$ provided $N\mu^{2r-1}\gg\log N ,\; \text{for}\;\; r\ge 2.$\\
Assume now $\mu\to\infty.$ Then
\begin{equation}\label{72}
\mathcal L_{1,N}^2\asymp\frac{1}{N},\;\;\mathcal L_{2,N}^2\sim\frac{\mu}{N},\;\;\mathcal L_{3,N}\asymp\sqrt{N}\mu^{3/2}e^{-cN\mu^{-4}},\;\;\mathcal L_{4,N}\asymp N\mu^{3/2}e^{-cN\mu^{-4}},
\end{equation}
whence $\bar{\mathcal L}_N\asymp(\mu/N)$ if\; $\mu\ll (N/\log N)^{1/4} .$\\

(ii)\; Let now $r=r(n,N)\to\infty.$ Set $ \psi(\mu):=\log(1+\mu^{-1})-(1+\mu)^{-1}. $ We have
\begin{equation}\label{73}
\mathcal L_{1,N}^2+\mathcal L_{2,N}^2\asymp \frac{r^{3/2}}{N\mu^r}\; \;\text{if}\;\; 0<\mu\le c<\infty,\;\; \text{and else}\;\; \mathcal L_{1,N}^2+\mathcal L_{2,N}^2\asymp \frac{1}{N}\Big[r^{3/2}\exp \{r/\mu^2\}+\mu\Big]\;\;\text{if}\;\; \mu\to\infty. 
\end{equation} 
Further, 
\begin{equation}\label{74}
\bar{\mathcal L}_N\asymp\frac{1}{Na_r}\asymp r^{3/2}\exp \{r\psi(\mu)\},\;\;\; \;\text{if}\;\; 0<\mu\le c<\infty,\;\;\text{provided}\;\; r^3\exp \{2r\psi(\mu)\} \ll\frac{N}{\log N}.
\end{equation}

If $\mu\to\infty$ then
\begin{equation}\label{75}
\bar{\mathcal L}_N\asymp\frac{\mu}{N}\;\; \text{if}\;\; r\ll\mu^{3/2},\;\;
\bar{\mathcal L}_N\asymp\frac{r^{3/2}}{N}\;\; \text{if}\;\; \mu^{3/2}\ll r\ll \mu^2,\;\;\text{and}\;\; \bar{\mathcal L}_N\asymp\frac{\mu^3}{N}\; \;\text{if}\;\;r\sim c\mu^2, \;\; \mu\ll (N/\log N)^{1/7} .
\end{equation}

\textbf{ Proof of Theorem 5.27 and 5.28.} 
We have $E(\xi-\mu)^4=\mu(1+\mu)^4\left(1+13\mu+15\mu^2\right)\le 15\mu(1+\mu)^6.$ We shall use the following two elementary inequalities: $a_r|r-\mu|/\mu\le1\,$ and $a_r(r-\mu)^2/\mu(1+\mu)^2\le1-a_r\le1.\,$ Indeed, if $r\ge\mu$, then $ra_r\le E\xi=\mu,$ whereas, if $r<\mu$, then
$a_r(\mu-r)/\mu\le a_r\le1.$ Furthermore, $cov^2(\varsigma,\xi)\le var\varsigma\,var\xi,$  and $cov(\varsigma,\xi)=a_r(r-\mu)$ we obtain $a_r^2(r-\mu)^2\le a_r(1-a_r)\mu(1+\mu)^2,$ which gives the second inequality.\\
\indent We have $E\varsigma^4=a_r(1-a_r)^4+(1-a_r)a_r^4\le a_r(1-a_r)\le a_r.$ Also,
$$\gamma^4E(\xi-\mu)^4\le \frac{a_r^4(r-\mu)^4}{\mu^4(1+\mu)^8}15\mu(1+\mu)^6=
15a_r\left(\frac{a_r\mid r-\mu\mid}{\mu}\right)^2\left(
\frac{a_r(r-\mu)^2}{\mu(1+\mu)^2}\right)\le 15a_r.$$
Therefore, $Eg^4(\xi)\le 8(E\varsigma^4+\gamma^4E(\xi-\mu)^4)\le 128a_r .$ Thus
$\beta_{4,N}\le 128 a_r/N\sigma ^4=128/Na_r(1-a_r)^2\theta ^2=\mathcal{L}_{1,N}^2. $ By this, (68) and (69) Theorem 5.27 follows from Theorem 3.2. \\    
   
Let's proof Theorem 2.28. Since $\sigma^2/Eg^4(\xi)\ge \,\sigma^2/128a_r$ we have,   
\begin{equation}\label{76}  
 \mathbf{T}^2\ge
C\min\left(\frac{\sigma^2}{a_r},\frac{\mu}{\sigma^2(1+\mu)^5}\right).
\end{equation} 

Next, $g(k)-g(l)=\mathbf 1\{k=r\}-\mathbf 1\{l=r\}
-\gamma(k-l).$ Therefore 
\[t(g(k)-g(l))+\tau(k-l)=t\bigl(\mathbf 1\{k=r\}-\mathbf 1\{l=r\}\bigr)+
u(k-l),\; \text{where} \;\; u=\tau-\gamma t.\]
We now use only adjacent pairs \((k+1,k)\). For such pairs,
\[H\left(\frac{t}{2\pi},\frac{\tau}{2\pi}\right)\ge 2\sum_{k=0}^{\infty}a_ka_{k+1}\left\langle\frac{u+\delta_k t}{2\pi}\right\rangle^2,\]  
where $\delta_k =\mathbf 1\{k+1=r\}-\mathbf 1\{k=r\}.$
That is, $\delta_{r-1}=1,\; \delta_r=-1,\;$ $\delta_k=0$ if $\,k\ne r-1,r.$
Define $A_r^+ =\,a_{r-1}a_r,\; A_r^-=a_ra_{r+1},$
and $A_r=A_r^++A_r^-.$ Also define
\[\bar{A}_r=\sum_{k\ge0,\; k\ne r-1,r}a_ka_{k+1}.\]
Then   
\[H\left(\frac{t}{2\pi},\frac{\tau}{2\pi}\right)\ge 2A_r^+\left\langle\frac{u+t}{2\pi}\right\rangle^2+2A_r^-\left\langle\frac{u-t}{2\pi}\right\rangle^2+
2\bar{A}_r\left\langle\frac{u}{2\pi}\right\rangle^2.\]
If \(m,n\in\mathbb Z\), then $\langle u+v\rangle\le |u-m+v-n|\le |u-m|+|v-n|.$ Taking the infimum over \(m,n\in\mathbb Z\) gives 
$\langle u+v\rangle\le \langle u\rangle+\langle v\rangle .$ 
Now, $\langle y\rangle=\langle (x+y)-x\rangle\le\langle x+y\rangle+\langle x \rangle .$ Put $a=\langle x+y\rangle,\quad b=\langle x\rangle.$ By Cauchy's inequality,
\[(a+b)^2=\left(\frac{\sqrt{\alpha}a}{\sqrt{\alpha}}+\frac{\sqrt{\beta}b}{\sqrt{\beta}}\right)^2\le\left(\alpha a^2+\beta b^2\right)\left(\frac1\alpha+\frac1\beta\right).\]

Hence $\alpha a^2+\beta b^2\ge\alpha\beta(\alpha+\beta)^{-1}(a+b)^2\ge\alpha\beta(\alpha+\beta)^{-1}\langle y\rangle^2,$ since $a+b\ge \langle y\rangle.$ Returning to $a$ and $b$, this gives
\[\alpha\langle x+y\rangle^2+\beta\langle x\rangle^2\ge\frac{\alpha\beta}{\alpha+\beta}\,\langle y\rangle^2.\]

Applying this inequality to the terms above gives
\[H\left(\frac{t}{2\pi},\frac{\tau}{2\pi}\right)\ge \frac{A_r\bar{A}_r}{4\pi^2(A_r+\bar{A}_r)}t^2,\]
since $|t|\le \pi,$ and hence $\left\langle t/2\pi\right\rangle=|t|/2\pi.$ 
Consequently,
\begin{equation}\label{77}
H_N=\inf_{\mathbf{T}\le |t|\le \pi,\; |\tau|\le \pi}\;H\left(\frac{t}{2\pi},\frac{\tau}{2\pi}\right)\ge \frac{A_r\bar{A}_r}{4\pi^2(A_r+\bar{A}_r)}
\mathbf{T}^2.
\end{equation}
We now simplify $A_r.$ We have $\; a_{r-1}/a_r=(r/(r+1))^{r-1}p^{-1}e^p,\; $  $\; a_{r+1}/a_r=(1+(r+1)^{-1})^rpe^{-p},\;\;r\ge 1.$  
Therefore $a_{r-1}\asymp a_r/p,\quad a_{r+1}\asymp p\,a_r.$
Consequently, for \(r\ge1\),
$$A_r=a_r(a_{r-1}+a_{r+1})\asymp a_r^2\left(\frac1p+p\right)\asymp \frac{a_r^2}{p},\quad r\ge1.$$

If \(r\ge2\), then \(k=0\) is included in the definition of \(\bar{A}_r\). Hence
$\bar{A}_r\ge a_0a_1.$ But $a_0=e^{-p},\quad a_1=p e^{-2p}.$ Therefore $\bar{A}_r\ge p e^{-3p}\ge e^{-3}p,\;\ r\ge 2.$ Thus
\begin{equation}\label{78}
\frac{A_r\bar{A}_r}{A_r+\bar{A}_r}\asymp\min(A_r,\bar{A}_r)\asymp \min\left(\frac{a_r^2}{p},p\right)\;\; r\ge 2, \;\; \text{and else},\;\; \frac{A_1\bar{A}_1}{A_1+\bar{A}_1}\asymp p^3.
\end{equation}

Combining (76),(77) and (78) we obtain
$ N\sigma\sqrt{\mu}(1+\mu) \exp\left\lbrace -cN H_N\right\rbrace= \mathcal{L}_{4.N}.  $  
By this and (17) Theorem 5.28 follows from Theorem 3.4.\\

\medskip
{\bf Proof of the assertions of Remark 5.6.} \; The case of fixed $r.$
Recall $\Delta_r(\mu):=1-a_r-a_r(r-\mu)^2/\mu(1+\mu)^2 , $ 
and $\sigma^2=a_r\Delta_r(\mu).$ We rewrite
$$T(\mu,a_r)=\min\left(\Delta_r(\mu),\frac{\mu}{a_r\Delta_r(\mu)(1+\mu)^5}\right).$$
Let $\mu\to0 ,\; r=1.$ 
 We have $a_1=\mu-3\mu^2+O(\mu^3).$ A direct expansion gives
$\sigma^2=6\mu^2+O(\mu^3).$ Consequently, $\Delta_1(\mu)\sim6\mu.$ Therefore,
$$\mathcal L_{1,N}^2\sim \frac{1}{36N\mu^3},\quad \mathcal L_{2,N}^2\sim\frac1{N\mu},\;\text{hence}\quad \mathcal L_{2,N}^2=o\left(\mathcal L_{1,N}^2\right), \quad \mathcal L_{3,N}\asymp\sqrt{N\mu}\,e^{-cN\mu}=o(\mathcal L_{1,N}^2).$$
Furthermore, $K_1(\mu,a_1)\asymp\mu^3.$ Also, since $\sigma^2/a_1\asymp\mu $ and
$\mu/\sigma^2(1+\mu)^5\asymp\mu^{-1},$ we obtain $T(\mu,a_1)\asymp\mu.$ Thus
$K_1(\mu,a_1)T(\mu,a_1)\asymp\mu^4.$ Therefore, since $\sigma\asymp\mu,$ we conclude that 
$\mathcal L_{4,N}\asymp N\mu^{3/2}e^{-cN\mu^4}.$ 

Assume $\mu\to0$, fixed $r\ge2.$ Put $d_r:=(r+1)^r/(r+1)!.$ Then $a_r\sim d_r\mu^r,\; \Delta_r(\mu)\to1 $ and $ \sigma^2\sim a_r. $ It follows that
$$\mathcal L_{1,N}^2\sim\frac1{Na_r}\sim\frac1{Nd_r\mu^r},\quad \mathcal L_{2,N}^2\sim\frac1{N\mu}=o\left(\mathcal L_{1,N}^2\right).$$
Also,
$\,\mathcal L_{3,N}\asymp\sqrt{N\mu}\,e^{-cN\mu}=o(\mathcal L_{2,N}^2).$
Further, we have $a_r^2(1+\mu)/\mu\asymp d_r^2\mu^{2r-1},$ whereas $\mu/(1+\mu)\sim\mu.$ Therefore $ K_r(\mu,a_r)\asymp\mu^{2r-1} $ since $r\ge2 .$ Furthermore,
$\sigma^2a_r\to1,$ while $\mu/\sigma^2(1+\mu)^5\asymp \mu^{1-r}\to\infty.$
Hence $T(\mu,a_r)\asymp1.$ 
Consequently, $K_r(\mu,a_r)T(\mu,a_r)\asymp\mu^{2r-1},$ and since $\sigma\asymp\mu^{r/2},$ we obtain
$\mathcal L_{4,N}\asymp N\mu^{(r+1)/2}e^{-cN\mu^{2r-1}}.$ Therefore,
$$\bar{\mathcal L}_N\sim\frac1{Nd_r\mu^r},\quad\mu\to0,\quad r\ge2,\quad \mu^{2r-1}\gg N^{-1}\log N.$$

Assume now $ 0<c_1\le\mu\le c_2<\infty,$ fixed $r\ge 1$. Then 
$a_r\asymp1,\quad \Delta_r(\mu)\asymp1,\quad \sigma^2\asymp1.$ Also both $K_r(\mu,a_r)$ and $T(\mu,a_r)$ are bounded away from 0 and $\infty$. Therefore,
$$\mathcal L_{1,N}^2\asymp\mathcal L_{2,N}^2\asymp\frac1N,\quad \mathcal L_{3,N}
\asymp\mathcal L_{4,N}\asymp\sqrt N\,e^{-cN}.$$
 Combining the above results, we obtain (70) and (71) .\\
 \indent
Let $\mu\to\infty.$ Put $\hat{d}_r=d_re^{-(r+1)}.$ Then $a_r\rightarrow \hat{d}_r
\quad \sigma^2\rightarrow \hat{d}_r(1-\hat{d}_r),\quad \Delta_r(\mu)\rightarrow1-\hat{d}_r.$
It follows that
$$\mathcal L_{1,N}^2\sim\frac1{N\hat{d}_r(1-\hat{d}_r)^2},\quad \mathcal L_{2,N}^2\sim\frac{\mu}{N},\quad \mathcal L_{3,N}\asymp\sqrt N\,\mu^{3/2}\exp\left\{-cN\mu^{-4}\right\}.$$

Furthermore, $\sigma^2/a_r\rightarrow1-\hat{d}_r\,$ whereas $\mu/\sigma^2(1+\mu)^5\asymp \mu^{-4}.$ Therefore, $\,T(\mu,a_r)\asymp\mu^{-4} $ and $K_r(\mu,a_r)\asymp 1.$ Hence $K_r(\mu,a_r)T(\mu,a_r)\asymp\mu^{-4}.$ Also $N\sigma\sqrt\mu(1+\mu)\asymp N\mu^{3/2}.$ we obtain
$$\mathcal L_{4,N}\asymp N\mu^{3/2}\exp\left\{-cN\mu^{-4}\right\}, \;\text{and hence}\;\; 
\bar{\mathcal L}_N\sim\frac{\mu}{N},\quad\mu\to\infty, \quad  \mu^4\ll N/log N .$$
The case fixed $r\ge 1$ is proved completely.\\

\textbf{The case $r=r(n,N)\to\infty. $} By Stirling's formula $a_r\asymp r^{-3/2}e^{-r\psi(\mu)}\to 0,$ therefore the first quantity in the definition of $K_r$ is eventually the smaller one. Thus
$$K_r(\mu,a_r)\asymp\frac{a_r^2(1+\mu)}{\mu},\quad \text{ also}\quad T(\mu,a_r)\asymp \min\left(1,\frac{\mu}{a_r(1+\mu)^5}\right).$$
Therefore, $K_r(\mu,a_r)T(\mu,a_r)\asymp K_r^*(\mu),$ and hence 
\begin{equation}\label{79}
 \mathcal L_{4,N}\asymp N(1+\mu)\sqrt{\mu a_r}e^{-cNK_r^*(\mu)},\quad\text{where}\quad   K_r^*(\mu):=\min\left(\frac{a_r^2(1+\mu)}{\mu},\frac{a_r}{(1+\mu)^4}\right).
 \end{equation}

As $r\to\infty$ one can see that $ \sigma ^2\sim a_r\asymp r^{-3/2}e^{-r\psi(\mu)}\to0 .$  It follows that
$$\mathcal L_{1,N}^2+\mathcal L_{2,N}^2\asymp\frac1N\left[\frac1{a_r}+\frac{(1+\mu)^2}{\mu}\right]\asymp\frac1N\left[r^{3/2}e^{r\psi(\mu)}+\frac{(1+\mu)^2}{\mu} \right].$$
Note that $\psi(\mu)\asymp -\log\mu,$ if $0<\mu\le c<\infty$, and $ \psi(\mu)\asymp \mu^{-2}$ if $\mu\to\infty.$ Applying this we obtain (73). 
 We next prove that the exponential terms are negligible. Since
$a_r\le1,\;ra_r\le  E\xi=\mu,$ we have $a_r\le\min(1,\mu)\le2\mu/(1+\mu).$ Therefore,
$$K_r^*(\mu)\le\frac{a_r}{(1+\mu)^4}\le\frac{2\mu}{(1+\mu)^5}.$$ 
Hence $NK_r^*(\mu)\gg\log N $ implies 
\begin{equation}\label{80}
\frac{N\mu}{(1+\mu)^5}\gg\log N.
\end{equation}
Set $A_N:=\mathcal L_{1,N}^2+\mathcal L_{2,N}^2.$ Since $A_N\ge\mathcal L_{2,N}^2=
(1+\mu)^2/N\mu,$
we obtain
$$\frac{\mathcal L_{3,N}}{A_N}\le N^{3/2}\frac{\mu^{3/2}}{1+\mu}\exp\left\{
-c_0\frac{N\mu}{(1+\mu)^5}\right\}.$$
The condition  (80) implies that $\mathcal L_{3,N}=o(A_N).$ Similarly,
$$\frac{\mathcal L_{4,N}}{A_N}\le N^2\sigma\frac{\mu^{3/2}}{1+\mu}\exp\{-cNk_r^*(\mu)\}.$$
Since $\sigma^2\sim a_r\le1,$ and $K_r^*(\mu)\le (1+\mu)^{-4}$ the condition $NK_r^*(\mu)\gg\log N$ implies $N(1+\mu)^{-4}\gg\log N,$ and therefore
$1+\mu=O(N^{1/4}).$ So, it follows that the prefactor $N^2\sigma\mu^{3/2}/(1+\mu)$ is bounded by a fixed power of $N$. Hence $\mathcal L_{4,N}=o(A_N).$ Thus
\begin{equation}\label{81}
\bar{\mathcal L}_N\asymp A_N\asymp\frac1N\left[\frac1{a_r}+\frac{(1+\mu)^2}{\mu}\right].
\end{equation}
It remains to simplify this expression in the individual regimes of $\mu.$\\
\indent If $\mu\to0$, then $ra_r\le\mu$ implies $a_r/\mu\le r^{-1}\to0,$ that is $a_r^{-1}\gg\mu^{-1}$,and hence $\mathcal L_{2,N}^2=o(\mathcal L_{1,N}^2).$ Moreover,
$K_r^*(\mu)\asymp a_r^2/\mu.$ 
 If $\mu$ is bounded away from 0 and $\infty$, then $(1+\mu)^2/\mu\asymp1,\quad
a_r^{-1}\gg 1,$ also $K_r^*(\mu)\asymp a_r^2.$  The relations (74) follows.\\

\indent Suppose now  $\mu\to\infty$. Then 
$$\bar{\mathcal L}_N\asymp\frac1N\left(\frac1{a_r}+\mu\right).$$
If $r=o(\mu^2)$, then $r\psi(\mu)\to0,$ because $\psi(\mu)\sim (2\mu^2)^{-1}.$
Hence $ a_r\asymp r^{-3/2}.$ Moreover, $a_r\mu^4\asymp \frac{\mu^4}{r^{3/2}}
\rightarrow\infty.$ 
Therefore, $K_r^*(\mu)\asymp a_r/\mu^4\asymp \mu^{-4}r^{-3/2},$ and
$$\bar{\mathcal L}_N \asymp \frac{r^{3/2}+\mu}{N}.$$
The transition between the two terms of the numerator is determined by $r\asymp\mu^{2/3}.$ If $r\sim c\mu^2,$ then $r\psi(\mu)\to c/2,$
and consequently $a_r\asymp\mu^{-3}.$ Thus $K_r^*(\mu)\asymp a_r/\mu^4\asymp
\mu^{-7}.$
Therefore,
$$\bar{\mathcal L}_N\asymp\frac{\mu^3}{N}.$$
Finally, if $r\gg\mu^2$, then $r\psi(\mu)\to\infty.$ Furthermore, 
$a_r\le Cr^{-3/2}=o(\mu^{-3}),$ so that $a_r\mu\to0.$ Hence $a_r^{-1}\gg\mu,$
and therefore
$$\bar{\mathcal L}_N\asymp\frac1{Na_r}\asymp\frac{r^{3/2}}{N}e^{r\psi(\mu)}.$$
In this regime,
$\, K_r^*(\mu)\asymp\min\left(a_r^2,a_r/\mu^4\right).$
This completes the proof of (75).\\

 \bigskip
 {\bf Example 5.4.2.} Consider $\mathbb{V}_N $- the total number of pairs of vertices in the forest; where a pair defined as two non-root vertices connected by a chain.  We note that the total number of distinct chains in the forest is $\mathbb{V}_N +n. $  We have
 
 $$\mathbb{V}_N=\frac{1}{2}\sum_{m=1}^N \eta_m(\eta_m -1).$$
    Set $a:= (3+\mu)\mu^2/2,$  and $\sigma^2:=3\mu^2(1+\mu)^5/2.$  

\begin{theorem}\label{th 5.29}
 There exist  $C>0$ such that 
  $$\sup_{-\infty <u<\infty} \mid \mathbf P\lbrace \mathbb{V}_N< u \sigma \sqrt{N}+N a\rbrace - \Phi(u)\mid\le C \Big(\bar{\mathcal L}_{1,N}+\mathcal{L}_{3,N}\Big), $$
  where 
$\,\bar{\mathcal L}_{1,N} =(1+\mu)^{3/2}/\sqrt{N}\mu.$ 
\end{theorem}

\begin{theorem}\label{th 5.30 }  Uniformly in $k\in \lbrace0,1,\dots ,n \rbrace $  
 $$ \mid \sigma\sqrt{N}\,\mathbf P \lbrace \mathbb{V}_N=k \rbrace -\frac{1}{\sqrt{2\pi}}e^{-u_k^2/2}\Bigg[ 1+\sqrt{\frac{2}{3N}}\Bigg((3u_k-u_k^3)\frac{(1+18\mu+27\mu^2)}{6\mu\sqrt{1+\mu}} + u_k\sqrt{1+\mu}\Bigg) \Bigg]\mid \le C \mathcal{D}_N, $$
  
where $u_k=(k-Na)/\sqrt{N}\,\sigma ,\;$ $\; \mathcal{D}_N = \bar{\mathcal L}_{1,N}^2 +  + \mathcal{L}_{3,N}+\bar{\mathcal L}_{4,N} , $

$$   \bar{\mathcal L}_{4,N}=  N\mu^{3/2}(1+\mu)^{7/2}\exp\left\{-cN\frac{\mu^2}{(1+\mu)^{11}}\min\left(\mu,\frac{1}{(1+\mu)^2}\right)\right\}  ,$$

\end{theorem}  
 
 \textbf{Proofs of Theorems 5.29 and 5.30.} Here $\varsigma=\xi(\xi-1)/2,\,$   $ g(\xi)=\varsigma - a- \gamma(\xi-\mu),\, $  $ \gamma = 3\mu (2+\mu)/2. $ One counts 
 \begin{equation}\label{82}
 \alpha_{3,0}=\frac{\sqrt6(1+18\mu+27\mu^2)}{3\mu\sqrt{1+\mu}},\quad \alpha_{1,2}=\sqrt{6(1+\mu)}. 
 \end{equation}
For $\alpha_{4,0}$ we get
\begin{equation}\label{83}
 \alpha_{4,0}=\frac{2+204p+1617p^2+2032p^3+420p^4}{3p^2(1-p)}\le\frac{c}{p^2(1-p)}= \frac{c(1+\mu)^3}{\mu^2}.
\end{equation}
Thus Theorem 5.29 follows from Theorem 3.2 , since (66),(68),  (83) and the fact that $\mathcal L_{2,N}\le \bar{\mathcal L}_{1,N}.$
 Further, by (68) and (83) we have $ \, M/\sigma^2\ge 2/3\mu(1+\mu)^{10}, $ and $\, \sigma^2/Eg^4(\xi)\ge =c/(1+\mu)^8. \,$ Therefore we define 

\begin{equation}\label{84} 
   \mathbf{T}^2=\frac{1}{(1+\mu)^8}\min\Big(1,\frac{1}{\mu(1+\mu)^2}\Big).
  \end{equation}
 Next, we define a lower bound for $ H_N $ using the neighboring atoms $0,1,2 .$ We use the events $(\xi,\xi')=(0,1) $ and $(\xi,\xi')=(1,2), $ then including the reverse pairs gives
 $$H\left(\frac{t}{2\pi},\frac{\tau}{2\pi}\right)\ge 2a_1\min(a_0,a_2)\Big[\Big\langle\frac{1}{2\pi}(t\gamma-\tau)\Big\rangle^2 + \Big\langle\frac{t}{2\pi}+ \frac{1}{2\pi}(t\gamma-\tau)\Big\rangle^2\Big]\ge cp^3\Big\langle \frac{t}{2\pi}\Big\rangle^2 = cp^3\frac{t^2}{4\pi^2}.$$
 since $|t|\le \pi$ and $\,a_0=e^{-p},\; a_1=pe^{-2p},\; a_2=3p^2e^{-3p}/2.$
 Thus, for $\mathbf{T}\le |t|\le \pi$,putting $p=\mu/(1+\mu)$ we get 
$$H_N\ge c\frac{N\mu^3}{(1+\mu)^3}\mathbf{T}^2.$$
Combining above results with (84),  (17) and Theorem 5 we complete the proof of  Theorem 5.30 . \\

\subsection{Random  permutations.}

\bigskip
\indent
 
 Consider the set $\mathcal{P}(n,N)$ of all permutations of degree $n$
having exactly $N$ cycles. A permutation may be viewed as a directed
graph whose connected components are cycles. We equip
$\mathcal{P}(n,N)$ with the uniform probability distribution and thus
obtain a random permutation with $N$ cycles.
To define an ordered vector of cycle lengths, we additionally label the
$N$ cycles by $1,\dots,N$, independently and uniformly over all $N!$
possible labelings. The number of elements in a cycle is called its
length. Let $\eta_m$ denote the length of the cycle labeled $m$,
$m=1,\dots,N$. Then $(\eta_1,\dots,\eta_N)$ is an exchangeable r.vec. satisfying $\eta_m\ge 1,\quad\eta_1+\cdots+\eta_N=n$  and we arrive at the GUM (1,N) generated by iid r.v.s $\xi_1,\dots,\xi_N$ having the common logarithmic-series distribution ( see Kolchin (1986)):  
$$ b_k:=\mathbf{P}\{\xi = k\} =\frac{p^k}{kL},\quad p \in (0,1), \; q=1-p, \quad L = -\ln q,\quad k=1,2,\dots.$$
We assume that $p$ is chosen so that $E\xi=p/(1-p)L =n/N=: \mu.$ \\
 \indent
We have
\begin{equation}\label{85}
 v^2:=var\xi=\frac{\mu(L-p)}{Lq}=\frac{p(L-p)}{q^2L^2},\;  \; \alpha_{0,4}\le c\Big( \frac{1}{p}+L\Big).
 \end{equation}

The last inequality follows from the exact formula:
\[\alpha_{0,4}=\frac{L^3(1+4p+p^2)-4L^2p(1+p)+6Lp^2-3p^3 }{p(L-p)^2}.\]
To get the bound we consider two regimes of $p.$ For $0<p\le 1/2$, we have $\,L-p\asymp p^2.$ Therefore the numerator of $\alpha_{0,4}$ does not exceed $C p^4,\,$ whereas denominator is $c\,p^3.\,$ That is $\alpha_{0,4}\le C\, p^{-1}.$ For $1/2<p<1$ straightforwardly we have $\alpha_{0,4}\le C\,L. $\\
Further, note that equation (68) still hold for any integer r.v., so using it we have
$$\;1-\mid Ee^{i\tau\xi}\mid^2\ge 4b_1b_2\sin^2 \frac{\tau}{2}\ge \frac{p^3}{\pi^2L^2}\tau^2.$$
Next, we have $v^2\alpha_{0,4}\le c\,(L-p)(1+pL)/ q^2L^2 .$ Apply this to get (see notation () )
\begin{equation}\label{86}
M=1-\sup_{(8v\sqrt{\alpha_{0,4}})^{-1}\le |\tau|\le\pi}\left|Ee^{i\tau\xi}\right|^2 \ge c\,\frac{p^3q^2}{(L-p)(1+pL)}.  
\end{equation}   

{\bf Example 5.5.1.} Let  $\mathbb{L}(r)$ be the number of cycles of length $r\ge 1$  in a random permutation of $\mathcal{P}(n,N).\,$ It is clear that 
$$ \mathbb{L}(r)=\sum_{k=1}^N \textbf{I}\lbrace\eta=r\rbrace,$$

 Set $ \varsigma=\textbf{I} \lbrace \xi=r\rbrace.\,$ We have $ E\varsigma =b_r,\, $ and 
\begin{equation}\label{87}
 \gamma = \frac{q b_r(r-\mu)}{\mu(1-\mu q)},\quad        \sigma^2 =  b_r(1-b_r)\theta,\quad 
\theta=1-\frac{b_r\left(rqL-p\right)^2}{(1-b_r)p(L-p)}\,.
\end{equation}

\begin{theorem}\label{th 5.31}
 There exist $C>0$ such that 
  $$\sup_{-\infty <u<\infty} \mid P\lbrace \mathbb{L}(r)< u \sigma \sqrt{N}+N b_r\rbrace - \Phi(u)\mid\le C\,\Big(\mathcal{G}_{1,N}+  \mathcal{G}_{2,N}\Big), $$
where
$$ \mathcal{G}_{1,N}=\frac{1}{\theta\sqrt{Nb_r(1-b_r)}}\sqrt{\frac{1}{p}+L},\quad \mathcal{G}_{2,N}=\sqrt{N\frac{p(L-p)}{q^2L^2}}\exp\left\lbrace -c\frac{Np^3q^2}{(L-p)(1+pL)}\right\rbrace .$$

\end{theorem}

\begin{theorem}\label{th 5.32 }  Uniformly in $k\in \lbrace 1,\dots ,N \rbrace $  
 $$ \mid \sigma\sqrt{N}P\lbrace \mathbb{L}(r)=k \rbrace -\frac{1}{\sqrt{2\pi}}e^{-u_k^2/2}\Bigg( 1+\frac{3u_k-u_k^3}{6 \sqrt{N}}\alpha_{3,0}+ \frac{u_k}{3\sqrt{N}}\alpha_{1,2}   \Bigg)\mid \le C\, \mathcal G_N, $$
  
where $u_k=(k-Nb_r)/\sqrt{N}\,\sigma,\;$ $\; \mathcal G_N =\mathcal G^2_{1,N}+ \mathcal G_{2,N} +\mathcal G_{3,N} , $  \\

$$\alpha_{3,0}=\frac{1}{\sigma^3}\Big[ b_r(1-b_r)(1-2b_r)-3\gamma b_r(r-\mu)(1-2b_r)+3\gamma^2 b_r\left((r-\mu)^2-\dfrac{\mu(1-\mu q)}{q}\right)$$
$$-\gamma^3 \mu\Big(\dfrac{1+p}{q^2}-\dfrac{3\mu}{q}+2\mu^2\Big)\Big],$$

$$ \alpha_{1,2}=\frac{p^{r-1}q}{r\sigma (1-\mu q)}\left[q(r-\mu)^2-\mu (1-\mu q)-\frac{q^2(r-\mu)}{1-\mu q}\left(\frac{1+p}{q^2}-\frac{3\mu}{q}+2\mu^2\right)\right].$$

$$\mathcal{G}_{3,N}= N\sigma\,v\exp\left\{-cN\kappa_r(p^{-1}+L) \min\left(
\frac{\sigma^2}{b_r(1-b_r)},\frac{p^2q^2}{(L-p)\sigma^2}\right)\right\},$$
$$ \kappa_1=p^5/L^2,\;\; \kappa_2=p^7/L^2 \;\; \text{and else}\;\; \kappa_r=p^{2r-1}/r^2L^2, \;\text{if}\;\; r\ge 3.  $$

\end{theorem}

{\bf Remark 5.7.}
If $p\to0,\; r\ge3;\quad\text{or}\quad p\to1,\; r\ge1;\quad\text{or}\quad 0<p_0\le p\le p_1<1,\; r\to\infty,$ then $\alpha_{3,0}\sim b_r^{-1/2}.$ If $p\to0$ and $r=1,2$, then
$$\alpha_{3,0}\sim(-1)^{r-1}\frac{\sqrt3}{p},\quad\alpha_{1,2}\sim
(-1)^{r-1}\frac{4\sqrt3}{3}.$$

If $p\to0$ and $r\ge3$, then
$$\alpha_{1,2}\sim\frac{2(r-1)(r-2)\sqrt{b_r}}{p}.$$
If $0<p_0\le p\le p_1<1,\quad r\to\infty,$ then
$$\alpha_{1,2}\sim\frac{r^2\sqrt{b_r}}{v^2}.$$

If $p\to1$, while $r\ge1$ is arbitrary.
Then
$$\alpha_{3,0}\sim b_r^{-1/2},\quad \alpha_{1,2}=\sqrt{b_r}\,Q_r(p)(1+o(1)),$$
where
$$Q_r(p)=\frac{q^2L^2}{p(L-p)}\,r^2-\frac{qL^2\bigl((1+p)L-p\bigr)}     {p(L-p)^2}\,r+\frac{pL^2}{(L-p)^2}.$$
In particular, if no cancellation occurs in the leading terms, namely if $\,qLr \nrightarrow 1/2 \,$ or $\,qr \nrightarrow 2 \,$, then
$$ Q_r(p)\sim 1-2qLr+q^2Lr^2,\quad \text{and hence}\quad \alpha_{1,2}\sim\sqrt{b_r}\bigl(1-2qLr+q^2Lr^2\bigr).$$
In the cancellation regimes, such as 
$qLr\to\frac12\quad\text{or}\quad qr\to2,$
the exact expression for $Q_r(p)$, or a higher-order expansion, must be used. For instance, if   $\,qLr=1/2 \,+o(L^{-1})$ then $\,Q_r(p)\sim 3/4L\,$ and $\,\alpha_{1,2}\sim 3\sqrt{2q}/4. $ Also, if $\,qr=2+o(L^{-1})\, $ then $Q_r(p)\sim -1$ and $\alpha_{1,2}\sim -e^{-1}\sqrt{q/2L} .$\\

{\bf Remark 5.8.} Suppose that \(r\ge1\) is fixed.     

\begin{enumerate}
\item If $p\to0$, then 
$$ \mathcal G_{1,N}+\mathcal G_{2,N}=O\Big(\frac{1}{\sqrt{N}p^2}\Big),\quad     \mathcal G_N=O\Big(\frac{1}{Np^4}\Big),\quad \text{if}\quad r=1,\; p^3\gg \frac{\log N}{N},\quad \text{or}\quad r=2,\; p^5\gg \frac{\log N}{N}.$$

For every fixed $r\ge3$,

$$ \mathcal G_{1,N}+\mathcal G_{2,N}=O\Big(\frac{1}{\sqrt{Np^r}}\Big),\quad     \mathcal G_N=O\Big(\frac{1}{Np^r}\Big),\quad \text{if}\quad  p^{2r-4}\gg \frac{\log N}{N}.$$
\item If $0<p_0\le p\le p_1<1,$ then 
$$ \mathcal G_{1,N}+\mathcal G_{2,N}=O\Big(\frac{1}{\sqrt{N}}\Big),\quad     \mathcal G_N=O\Big(\frac{1}{N}\Big).$$

\item If \(p\to1\), equivalently \(q\to0\) and \(L\to\infty\), then
$$\mathcal G_{1,N}+ \mathcal G_{2,N}=O\Big(\frac{L}{\sqrt{N}}\Big),\quad \mathcal G_N=O\Big(\frac{L^2}{N}\Big),\quad\text{if}\quad \frac{q}{L}\gg\sqrt{\frac{\log N}{N}}.$$

\end{enumerate}

{\bf Remark 5.9.} Suppose that $r=r(N)\to\infty$. \\
(1)\; If $p\le c<1$, then
$$\mathcal G_{1,N}+\mathcal G_{2,N}\sim\sqrt{\frac{r}{Np^r}},\quad \mathcal G_N\sim\frac{r}{Np^r},\quad\text{if}\quad  Np^{2r}{r^2}\gg\log N .$$

(2)\;If $p\to1$ and hence $q\to 0$, $L\to\infty $ , then
$$\mathcal G_{1,N}^2\sim\frac{rL^2}{Np^r},\quad \mathcal G_{2,N}=\frac{\sqrt N}{q\sqrt L}\exp\left\{-c\frac{Nq^2}{L^2}\right\},\quad \mathcal G_{3,N}\asymp\frac{Np^{r/2}}{q\sqrt r\,L}\exp\left\{-cN\frac{p^r}{rL}\min\left(\frac{p^r}{r},q^2\right)\right\}. $$
Thus, this case is more complicated. A restriction on the growth
rate of $r$ is needed in order to ensure that the exponential
terms are negligible compared with $\mathcal G_{1,N}^2$.  In particular, suppose that $r=O(q^{-1}).$ Then $p^r\asymp1.$ Moreover,$r^{-1}p^r\gtrsim q\gg q^2,$ so the \(q^2\)-branch is the smaller branch. Consequently,
$$\mathcal G_N\sim\frac{rL^2}{Np^r},\quad\text{provided that}\quad
\frac{q^2}{L^2}\min\left(1,\frac{Lp^r}{r}\right)\gg\frac{1}{N}\log\frac{N}{qL}.$$
Equivalently, the displayed condition can be written as
$$\frac{q}{L}\gg\sqrt{\frac{1}{N}\log\frac{N}{qL}},$$
together with one of the following two alternatives:
$r^{-1}p^r\ge L^{-1},$ or $ L^{-1}>r^{-1}p^r\gg (L/Nq^2)\log(N/qL).$
Obviously, if, in addition, $ -\log (qL)=O(\log N),$ then these conditions  may be slightly simplified.

If $rq\to\infty$, then one can see that the exponent in $\mathcal G_{3,N}$
is the smaller than the exponent in $\mathcal G_{2,N}.$ Therefore,

$$\mathcal G_N\sim\frac{rL^2}{Np^r},\quad\text{provided that}\quad
\frac{Np^r}{rL}\min\left(\frac{p^r}{r},q^2\right)\gg \log\frac{N}{qL}.$$

The last condition reduces to $r^{-1}p^r\gg (l/Nq^2)\log(n/qL)$ if $r^{-1}p^r\ge q^2 $ whereas, if $r^{-1}p^r\le q^2)$, it reduces to
$r^{-2}p^{2r}\gg (L/N)\log (N/qL).$\\

\textbf{Proof of Theorems 5.31 and 5.32}. Let us prove the bound for the $\alpha_{4,0} .$ Write $g(\xi)=\varsigma-b_r-\gamma (\xi-\mu),$ then $Eg^4(\xi)\le 8(E(\varsigma-b_r)^4+\gamma^4E(\xi-\mu)^4).$ We have $E(\varsigma-b_r)^4\le b_r(1-b_r).$
Next, by Cauchy's inequality, $\operatorname{cov}(\varsigma,\xi)^2\le b_r(1-b_r)v^2.$ Whence $\gamma^2\le b_r(1-b_r)/v^2.$ Consequently,  
$$\gamma^4E(\xi-\mu)^4\le b_r^2(1-b_r)^2 v^{-4}E(\xi-\mu)^4\le  b_r(1-b_r)\alpha_{0,4}/4.$$
Therefore
\begin{equation}\label{88}
\alpha_{4,0}=\sigma^{-4}Eg^4(\xi)\le \frac{1}{4\sigma^4} b_r(1-b_r)\left(1+\frac{1}{p}+L \right).
\end{equation}
Thus Theorem 5.31 follows from Theorem 3.2 by using (85), (86), (88) and that $\,\sqrt{b_r(1-b_r)}\gg \sigma^2$ and hence $ \varrho_3(n)\le\sqrt{N^{-1}\alpha_{0,4}}\le \mathcal G_{1,N} $.\\
\indent
Let us find a lower bound for $H_N.$  From (87) and (88) it follows that 
$$ \sigma^2/Eg^4(\xi)\ge c\frac{\sigma^2p}{b_r(1-b_r)(1+Lp)}.$$
By this and (86) we have
\begin{equation}\label{89}
\mathbf{T}^2  \ge c\frac{p}{1+Lp}\min\Big(\frac{\sigma^2}{b_r(1-b_r)},\frac{p^2q^2}{(L-p)\sigma^2}\Big).
\end{equation}
Next, $ \Delta_k:= g(k+1)-g(k)=\textbf{I}\lbrace k+1=r\rbrace -\textbf{I}\lbrace k=r\rbrace  -\gamma.$
The event $(\xi,\xi')=(k+1,k)$ has  probability to occur  $w_k:=b_{k+1}b_k.$ Therefore for two adjacent indices $k,\ell$
\[H\left(\frac{t}{2\pi},\frac{\tau}{2\pi}\right)\ge \sum_{j=1}^{\infty}w_j
\left\langle\frac{t\Delta_j+\tau}{2\pi}\right\rangle^2\ge w_k\left\langle
\frac{t\Delta_k+\tau}{2\pi}\right\rangle^2+w_\ell\left\langle\frac{t\Delta_\ell+\tau}{2\pi}\right\rangle^2 \ge\frac{w_kw_\ell}{w_k+w_\ell}\left\langle\frac{t(\Delta_k-\Delta_\ell)}{2\pi}\right\rangle^2 ,\]
where  the inequality $a\langle x\rangle^2+b\langle y\rangle^2\ge
ab(a+b)^{-1}\langle x-y\rangle^2,\quad a,b>0$ is used  with
$x=(t\Delta_k+\tau)/2\pi,$ and $y=(t\Delta_\ell+\tau)/2\pi.$ 
Since $w_kw_\ell/(w_k+w_\ell)\asymp\min(w_k,w_\ell),$ we get
\[H\left(\frac{t}{2\pi},\frac{\tau}{2\pi}\right)\ge c\min(w_k,w_\ell)
\left\langle\frac{t(\Delta_k-\Delta_\ell)}{2\pi}\right\rangle^2.\]
If $|\Delta_k-\Delta_\ell|=1,$ then for $\mathbf{T}\le |t|\le\pi$, we have
$ \left\langle t/2\pi\right\rangle ^2\ge c\mathbf{T}^2.$ 
Consequently,

$$H\left(\frac{t}{2\pi},\frac{\tau}{2\pi}\right)\ge c\min(w_k,w_\ell)\mathbf{T}^2.$$

We now choose $k,\ell$ so that $|\Delta_k-\Delta_\ell|=1$.
First, let $r=1$.  Then we chose $(k,\ell)=(2,1)$ and $(k,\ell)=(3,2).$ 
The corresponding weights are $w_1=b_2b_1,\; w_2=b_3b_2,$ and $\min(w_1,w_2)\asymp p^5/L^2.$ Thus 
\[H_N\ge cN\frac{p^5}{L^2}\mathbf{T}^2,\quad r=1.\]

Next, let $r=2$. Here use $(k,\ell)=(2,3)$. Then the corresponding weights are
$w_2=b_3b_2,\quad w_3=b_4b_3$ with  $\min(w_2,w_3)\asymp p^7/L^2.$
Therefore the uniform bound is
\[H_N\ge cN\frac{p^7}{L^2}\mathbf{T}^2,\quad r=2.\]

Finally, let $r\ge3$. Use the special adjacent pair $(k,\ell)=(r-1,1)$ . Then
the weights are $w_{r-1}=b_rb_{r-1},$ and $w_1=b_2b_1,$ with 
Thus
$$\min(w_{r-1},w_1)=\min\left(\frac{p^{2r-1}}{r(r-1)L^2},\frac{p^3}{2L^2}
\right)=\frac{p^{2r-1}}{r(r-1)L^2}.$$
 Consequently,
\[H_N\ge cN\frac{p^{2r-1}}{r(r-1)L^2}\mathbf T^2,\quad r\ge3.\]
  
Combining the cases, we obtain $H_N\ge cN\kappa_r \mathbf T^2.$ By this, (17), (85), (86), (88),(89) and Theorem 3.4 we conclude Theorem 5.32.\\

The asymptotic results presented in above remarks can be obtained by routine asymptotic analysis of the parameters of $b_r,$ the calculations are omitted.

\bigskip
{\bf Example 5.5.2.} We say that two distinct elements of a permutation form a pair if they belong to the same cycle. Let $\mathbb{W}_N $ stands for the total number of  pairs.   It is clear that 
$$ \mathbb{W}_N=\frac{1}{2}\sum_{k=1}^N\eta(\eta-1) . $$
 Note that $2\mathbb{W}_N$ is the number of simple (non-closed) directed paths in the permutation graph.\\

 As above we set $\varsigma=\xi(\xi-1)/2$. We get $ E\varsigma=p^2/2q^2L=L\mu^2/2,$ and  $\operatorname{cov}(\varsigma,\xi)=p^2(2L-p)/2q^3L^2,$ whence 
$$\gamma=\frac{p(2L-p)}{2q(L-p)}, \quad \sigma^2=\frac{p^2B}{4q^4L(L-p)},\quad B:=2(L-p)-p^2.$$

\begin{theorem}\label{th 5.33}
 There exist  $C>0$ such that 
  $$\sup_{-\infty <u<\infty} \mid P\lbrace \mathbb{W}_N< u \sigma \sqrt{N}+ n\mu L/2\rbrace - \Phi(u)\mid\le C\,\Big(\sqrt{\frac{1+p^2L}{Np^2} }\,+  \mathcal{G}_{2,N}\Big), $$

\end{theorem}

\begin{theorem}\label{th 5.34 }  Uniformly in $k\in \lbrace 1,\dots ,N \rbrace $  
 $$ \mid \sigma\sqrt{N}P\lbrace \mathbb{W}_N=k \rbrace -\frac{1}{\sqrt{2\pi}}e^{-u_k^2/2}\Bigg( 1+\frac{3u_k-u_k^3}{6 \sqrt{N}}\alpha_{3,0}+ \frac{u_k}{3\sqrt{N}}\alpha_{1,2}   \Bigg)\mid \le C \widehat{\mathcal{G}}_N, $$
  
where $u_k=(k-Nb_r)/\sqrt{N}\,\sigma ,\;$ 
$$\widehat{\mathcal{G}}_N =\frac{1+p^2L}{Np^2} + \mathcal{G}_{2,N}+\widehat{\mathcal{G}}_{3,N} ,\quad \alpha_{1,2}=\sqrt{\frac{ L^3\,B }{(L-p)^3}}, $$  

$$\alpha_{3,0}=\frac{B}{2pL}\,\alpha_{1,2}\left[p^2+6p+1+\frac{3p^3(p+4)}{B}+\frac{p^4(3p^2+6p+1)}{B^2}+\frac{p^8}{B^3}\right],$$
$$\widehat{\mathcal G}_{3,N}\le
\begin{cases}
\displaystyle
C N p^{3/2}\exp\{-cNp^3\},
& 0<p\le c_1<1,
\\[10pt]
\displaystyle
\frac{CN}{q^3L}
\exp\left\{-c\frac{Nq^6}{L^3}\right\},
& p\to1.
\end{cases}$$
\end{theorem}

\textbf{Remark 5.10.} We have
$$\widehat{\mathcal G}_N\asymp\frac{1+p^2L}{Np^2}\quad\text{ if} \quad 0<p\le c_1<1,\; p^3\gg \frac{\log N}{N},\quad \text{or \; if} \quad p\to 1,\;\; qL^{-1}\gg\sqrt{N^{-1}\log(N/qL)}, $$ 
since the exponent of $\mathcal G_{2,N}$ is smaller than the exponent of $ \widehat{\mathcal G}_{3,N}$ in the case $p\to 1.$\\

\textbf{Proof of Theorems 5.33 and 5.34.}
The  arguments similar to those used above  give 
$$\alpha_{4,0}\le C\left(\frac1{p^2}+L\right).$$
Therefore, Theorem 5.33 follows  from Theorem 3.2 since (87), (88).

For $\alpha_{1,2}$, we use the following exact identity. Recall $g(\xi)=\varsigma-E\varsigma-\gamma(\xi-\mu).$

\begin{lemma}
 One has $Eg(\xi)(\xi-\mu)^2=2\sigma^2.$ 
\end{lemma}
\textbf{Proof.} Put for brevity $Y:=\xi-\mu,\; v^2=EY^2,\; m_3:=EY^3,\; m_4:=EY^4.$ 
Since $\xi=\mu+Y $  we have
\[\varsigma=\frac{Y^2}{2}+\left(\mu-\frac12\right)Y+\frac{\mu^2-\mu}{2}\]
and we write $E\varsigma=v^2/2+ (\mu^2-\mu)/2.$ Therefore
$\varsigma-E\varsigma=(Y^2-v^2)/2+(\mu-  1/2)Y.$
Whence $\operatorname{cov}(\zeta,\xi)=E(\zeta-E\zeta)Y=EY^3/2+(\mu- 1/2)EY^2.$ Hence $\gamma=m_3/2v^2+\mu-1/2.$
Thus combining above equations after simplification we have 
$$g(\xi)=\zeta-E\zeta-\gamma Y=\frac12\left(Y^2-v^2-\frac{m_3}{v^2}Y\right).$$
On using this we get
\[Eg(\xi)Y^2=\frac12\left(m_4-v^4-\frac{m_3^2}{v^2}\right).\]
On the other hand, $\sigma^2=Eg^2(\xi),$
and
\[\sigma^2=\frac14E\left(Y^2-v^2-\frac{m_3}{v^2}Y\right)^2=\frac14\left\{E(Y^2-v^2)^2-2\frac{m_3}{v^2}E(Y^2-v^2)Y+\frac{m_3^2}{v^4}EY^2\right\}.\]
But $E(Y^2-v^2)^2=m_4-v^4,$ and $E(Y^2-v^2)Y=m_3,$ and $EY^2=v^2.$
Hence
\[\sigma^2=\frac14\left(m_4-v^4-2\frac{m_3^2}{v^2}+\frac{m_3^2}{v^2}\right)
=\frac14\left(m_4-v^4-\frac{m_3^2}{v^2}\right).\]
Consequently,
\[Eg(\xi)Y^2=\frac12\left(m_4-v^4-\frac{m_3^2}{v^2}\right),\]
whereas
\[2\sigma^2=\frac12\left(m_4-v^4-\frac{m_3^2}{v^2}\right).\]
Lemma follows.\\
It follows from lemma $\alpha_{1,2}=2\sigma/v^2.$
Equivalently,
$\alpha_{1,2}=\sqrt{L^3B/(L-p)^3}.$

Further, a direct expansion gives $Eg^3(\xi)=p^2P_3(p,L)/8q^6L(L-p)^3,$
where
$$P_3(p,L)=\frac12\left[B^3(p^2+6p+1)+3B^2p^3(p+4)+Bp^4(3p^2+6p+1)+p^8\right].$$
So we obtain
$$\alpha_{3,0}=\frac{\sqrt L\,P_3(p,L)}{p(L-p)^{3/2}B^{3/2}}.$$
which coincides with formula given above.\\
\indent
On using (86) and above obtained bound for $Eg^4(\xi)$ we obtain
\begin{equation}\label{90}
\mathbf{T}^2\ge\frac{c}{\sigma^2}\min\left(\frac1{1/p^2+L},\frac{p^2q^2}{(L-p)(1/p+L)}\right).
\end{equation}
Assume that \(0<\mathbf T\le \pi\). Then will prove that
$$H_N\ge\frac{Np^4}{2\pi^2L^2(3+8p)}\,\mathbf T^2\ge\frac{Np^4}{22\pi^2L^2}\,\mathbf T^2.$$
Put $u:=t/2\pi,\quad w:=(\tau-\gamma t)/2\pi.$ Then  
$tg^*/2\pi +\tau\xi^*/2\pi=u(\varsigma-\varsigma')+w(\xi-\xi').$
Consider first the event $\{(\xi,\xi')=(2,1)\}$ on which 
$\varsigma-\varsigma'=1,\quad\xi-\xi'=1,$ and therefore
$u(\varsigma-\varsigma')+w(\xi-\xi')=u+w.\,$ On the event $\{(\xi,\xi')=(3,1)\},$ we have $u(\varsigma-\varsigma')+w(\xi-\xi')=3u+2w.$
Taking into account the reversed pairs $(\xi,\xi')=(1,2)\quad\text{and}\quad
(\xi,\xi')=(1,3),$ and using \(\langle-x\rangle=\langle x\rangle\), we obtain
$$H\left(\frac{t}{2\pi},\frac{\tau}{2\pi}\right)\ge 2b_1b_2\langle u+w\rangle^2+2b_1b_3\langle3u+2w\rangle^2\ge b_1\min(b_2,2b_3)(2\langle u+w\rangle^2 +\langle 3u+2w\rangle^2)$$
$$\ge b_1\min(b_2,2b_3)(\langle 2(u+w)\rangle^2 +\langle 3u+2w\rangle^2)=b_1\min(b_2,2b_3)(\langle 2(u+w)\rangle^2 +\langle 2(u+w)+u\rangle^2)$$
$$\ge b_1\min(b_2,2b_3)\langle u \rangle^2, $$
since $2(u+w)+u=(3u+2w),$ the elementary inequalities
$\langle x+y\rangle\le\langle x\rangle+\langle y\rangle,\quad \langle mx\rangle\le|m|\langle x\rangle,\quad m\in\mathbb Z,$
give $\langle u\rangle\le 2\langle u+w\rangle+\langle3u+2w\rangle.$
Thus, since $\mid u\mid\le 1/2$ , substituting instead of $b_k$s their exact form we obtain
$$H\left(\frac{t}{2\pi},\frac{\tau}{2\pi}\right)\ge \frac{p^4}{9L^2}t^2. $$ 

Taking the infimum over $\mathbf T\le |t|\le\pi,\quad |\tau|\le\pi,$
gives $H_N\ge Np^4\,\mathbf T^2/9L^2.$ Equivalently, since (90),
$$H_N\ge cN\min\left\{\frac{p^4q^4(L-p)}{LB(1+p^2L)},\frac{p^5q^6}
{LB(1+pL)}\right\}.$$

Consequently, the following asymptotic estimates hold:

Suppose first that $p\to0$. Then  
\[L=p+\frac{p^2}{2}+\frac{p^3}{3}+O(p^4),\]
and hence $L-p\sim p^2/2,\quad B=2L-2p-p^2\sim 2p^3/3.$
The two quantities in the minimum have orders $p^2$ and $p$,
respectively. Therefore, $H_N\gtrsim Np^2.$ 
As $0<\mathbf T\le\pi$ and $\sigma^2\asymp p^2,$ we obtain $H_N\lesssim Np^2.
$ Thus $H_N\asymp p^2.$ 

If $p$ is bounded away from both zero and one, all the quantities
appearing in the lower bound are bounded above and below by positive
constants. Hence $H_N\ge c N.$ \\
Finally, suppose that $p\to1$. Then $L-p\sim L,\quad 2L-2p-p^2\sim2L,\quad
1+p^2L\sim L,\quad 1+pL\sim L.$ Therefore, the first quantity in the minimum is of order $q^4/L^2,$ whereas the second one is of order $q^6/L^3.$
Therefore we conclude that $H_N\ge c q^6/L^3 .$ Theorem 5.34 now follows from Theorem 3.4.


\bigskip
\bigskip
\begin{thebibliography}{99}

\bibitem{1} Bartlett,M.S.:The characteristic function of a conditional statistic. J.Lond.Math.13,62-67(1938)

\bibitem{2}  Barlevy Gadi and  Nagaraja H. N. Properties of the Vacancy Statistic in the Discrete Circle Covering Problem . Springer Proceedings in Mathematics and Statistics ((PROMS,v. 149)), 121 - 146. (2015)

\bibitem{3}Bhattacharya R.N. and R. Ranga Rao (1976). Normal approximation and asymptotic Expansions. John Wiley and Sons, New-York.	

\bibitem{3}    von Bahr, b. (1972). On sampling from a finite set of independent random variables. Z. Wahrsch. Verw. Geb. 24, 279–286.

\bibitem{4} Bickel P.J.    bickel, p. j. and von Zwet, W.R. (1978). Asymptotic expansions for the power of distribution-free tests in the two-sample problem. Ann. Statist. 6, 937–1004.

\bibitem{4} Cressie N.A.C. and Read T.R. C. (1984). Multinomial goodness-of-fit tests. J. Roy. Statist. Soc. Ser B. 46, pp. 440-464.

\bibitem{5} Devroye L., A limit theory for random skip lists. Ann.Appl. Probabl. v.2,597-609 (1992). 

\bibitem{6} Feller  V., Introduction to  probability theory and its applications., John Wiley and Sons, Inc. New York, v.2,(1966).

\bibitem{7} Gnedin A.,Hansen, B. and Pitman, J. (2007). Notes on the occupancy problem with
infinitely many boxes: General asymptotics and power laws. Probab. Surveys 4, 146 - 171.

\bibitem{8}Holst L., (1979). A unified approach to limit theorems for urn models. J. Appl., Probab., 16,Issue 1, 154-162.

\bibitem{9} Holst, L.( 1985). On discrete spacings and the Bose–Einstein distribution. In: Lanke, J., G. Lindgren (Eds.), Contributions to probability and statistics in honour of Gunnar Blom, pp. 169–177. Lund.

 \bibitem{10} Holst, L., Rao, J.S.: Asymptotic theory for families of two-sample non-parametric statistics. Sankhya, 42(Ser. A), 19-52 (1980)

\bibitem{11} Hsien-Kuei Hwang and Svante Jonson.(2008) Local limit theorems for finite and infinite  urn models. Ann.Probabl. 36,,pp. 992–1022.
 
\bibitem{11}  Jammalamadaka S.R., Schweitzer, R.L.: On tests for the two-sample problem based on higher order spacing
frequencies. In: Matusita, K. (ed.) Statistical Theory and Data Analysis, pp. 583-618. Nort-Holland,
Amsterdam (1985)


\bibitem{12} Johnson, N.L., Kotz, S.: Urn Models and Their Applications. Wiley, New York (1977)

\bibitem{13} Kotz S. and Balakrishnan N., (1997). Advances in urn models during the past two decades. In Advances in Combinatorial Methods and Appl. to Probab. and Statist., 203-257. Birkhauser, Boston. MA.

\bibitem{14}Kolchin V.F., Sevast'yanov, B.A., and Chistyakov, V.P.  (1978). Random Allocations,  V.H. Winston and Sons, Washington, DC.
\bibitem{15} Kolchin V.F., (1986) Random Maps. Springer.(Translations Series in Mathematics and Engineering)
\bibitem{16} Kolchin,V.F.(2010). Random graphs, 2nd ed., Cambridge Univ. Press, Cambridge, 2010 

\bibitem{17} Ibrahim B.M. and Mirakhmedo S.M., (2016), Approximation by Normal Distribution for a Sample Sum in Sampling Without Replacement from a Finite Population. Sankhya ,  V. 78-A, No 2, pp. 188-220.

\bibitem{17}   Ismatullaev Sh.A. and Mirakhmedov Sh.A.(1993), Asymptotic expansion for permutation tests with several samples. Kolchin V.F. (ed) et al., Probabilistic methods in discrete mathematics. Proceedings of the Petrozavodsk Conference. Petrozavodsk (Russia), May 11-16, 1992. Moscow, TVP, Prog. Pure and Appl. Discrete Math. 1, p.226-239 (1993)


\bibitem{18} Ivanov, V.A., Ivchenko, G.I., Medvedev Yu, I.: Discrete problems of the probability theory (a survey).
J. Sov. Math. 31(2), 3-60 (1985)


\bibitem{19} Ivchenko G.I. (1994). On the random covering of a circle: a discrete model (in Russian). Diskretnaya Matematika 6(3): 94–109.

\bibitem{20} Ivchenko G.I. and  Mirakhmedov S.M.(1995) Large deviations and intermediate efficiency of the decomposable statistics in multinomial scheme. Math. Methods Statist. 4 (1995)p. 294-311.

 \bibitem{21}  L'ecuyer P., Simard R. and Wegenkittl S. ,  Sparse serial tests of uniformity for random number generators. SIAM J. Sci. Comput.  24,  pp. 652– 668, (2002). 
\bibitem{22} Mirakhmedov S.M. ,  Jammalamadaka S.R. and Ibrahim B. M (2014), On Edgeworth Expansions in Generalized Urn Models.. J Theor Probab (2014) 27:725–753 DOI 10.1007/s10959-012-0454-z

 \bibitem{23}	Mirakhmedov Sh.A. (1983). Asymptotical expansion of the distribution of the sample sum from a finite population. Theory Probabl.and Appl., v.28, p. 468-477. 

\bibitem{24}	Mirakhmedov S.M (1996), Limit theorems on decomposable statistics in generalized allocation schemes. Discrete Math.Appl.,6, №4, p. 379-404. 

\bibitem{25}	 Mirakhmedov Sh.A, (1992) Probability of moderate deviations in limit theorems for conditional distribution of sums of independent random variables. Discrete Math. Appl. 4,  p.22-32.(In Russian)
	

\bibitem{26}	Mirakhmedov S. M.(1992).Randomized decomposable statistics in the scheme of independent allocating particles into boxes. Discrete Math. Appl. 2 , 91–108.

 \bibitem{27} Mirakhmedov S.M and Saidova O.A.,(1998) Estimates for remainder term on CLT for empty boxes statistics. Turkish J. Math., 22,  47-52.
\bibitem{16}	Mirakhmedov S. M.( 2007). Asymptotic normality associated with generalized occupancy problem. Statist. Probab. Letters, 77, 1549-1558.

\bibitem{28}	Mirakhmedov S.S and Mirakhmedov S.M. (2010) On Asymptotic Expansion in the Random Allocation of Particles by Sets. J Theor Probab., v.23: pp. 904–919 DOI 10.1007/s10959-009-0225-7 

\bibitem{29}	Mirakhmedov,S.M.,Jammalamadaka S.R. and Ibragim Bin Muhamad.(2014), On Edgeworth Expansion in Generalized Urn Models.. J. Theor. Probab.  27: 725–753, 
DOI 10.1007/s10959-012-0454-z.


\bibitem{30}	Mirakhmedov , S.M., Jammalamadaka S.R.and  Ekstrom M.(2015). Edgeworth expansions for two-stage sampling with applications to stratified and cluster sampling. Canadian J. Statist. Vol. 43, No. 4,  Pages 578–599.  DOI: 10.1002/cjs.11266	

\bibitem{31}  Mirakhmedov S.M. (2020). The Probabilities of Large Deviations for a Certain Class of Statistics Associated with Multinomial Distributions. ESAIM , Probability and Statistics . v.24, p. 581-606 (Correction, 25 (2021) 53-54),http://arxiv.org/abs/2205.03317 [math.PR](2022).  .

\bibitem{32} Mikhaylov V.G., (1981).A central limit theorem for the independent allocation scheme. Proc. Steklov Inst.Math., 157, p.138-152 (in Rassion)

\bibitem{33} Mikhaylov V.G., (1991). Polinomial and polinomial like allocations: recent developments. Probabl. methods in discrete Math. pp. 40-59, V.F. Kolchin et al Eds.

\bibitem{34} Mukhin A.B.,  Local limit theorems for lattice random variables. Theory Probab. Appl. 36, 660–674,(1991) .

\bibitem{35} Pavlov Yu. L. (1977).The asymptotic distribution of maximum tree size in a random forest. Probabl. Theor. Appl., v.22, pp.523-533.

\bibitem{36} Pavlov Yu., Cherepanova Ye., Limit distributions of the umber of pairs in a generalized allocation scheme. Discrete Math. Appl., 14, 149 - 159.
\bibitem{37}, Petrov V.V.: Sums of Independent Random Variables. Springer, New York (1995)
\bibitem{38} Vatutin, V.A., Mikhailov, V.G.(1982), Limit theorems for the number of empty cells in an equiprobable scheme for group allocation of particles. Theory Probab. Appl. 27, 734–743 


\end{thebibliography}
\end{document}